\documentclass{ndjflart_igpl}
\usepackage[T1]{fontenc}
\usepackage{tgtermes}
\usepackage{mathptmx}  
\usepackage[scaled=.92]{helvet}
\usepackage{amsthm,amsmath,amssymb}
\usepackage{mathrsfs}
\usepackage[numbers]{natbib}  %% numbers is required.
\usepackage[colorlinks,citecolor=blue,urlcolor=blue]{hyperref}  %%check

\usepackage{caption}

\usepackage{graphicx}
\usepackage{xspace}

\usepackage{tikz}
\usepackage{bpextra}

\usepackage{txfonts}

\usepackage{stmaryrd}

\usepackage{multicol}
\usepackage{multirow}

\usepackage{amscd}
\usepackage{mathptmx}
\usepackage{mathrsfs}
\usepackage{color}
\usepackage{xspace}
\usepackage{bussproofs}

\usepackage{cleveref}
\renewcommand{\cref}{\Cref}

\usepackage{array}
\newcolumntype{C}[1]{>{\centering\arraybackslash}p{#1}}
\newcolumntype{L}[1]{>{\arraybackslash}p{#1}}

\usepackage{enumerate}

\artstatus{am} %%% leave this alone!!  That means you, too!!
\newtheorem{theorem}{Theorem}[section]
\newtheorem{lemma}[theorem]{Lemma}

\newtheorem{corollary}[theorem]{Corollary} 
\theoremstyle{definition}
\newtheorem{definition}[theorem]{Definition}

\theoremstyle{remark}

\newtheorem{proposition}[theorem]{Proposition}

\newcommand{\MD}{\ensuremath{\mathbf{MD}}\xspace}
\newcommand{\MDe}{\ensuremath{\mathbf{MD}^+}\xspace}
\newcommand{\MID}{\ensuremath{\mathbf{MID}}\xspace}
\newcommand{\MIDz}{\ensuremath{\mathbf{MID^-}}\xspace}
\newcommand{\MDor}{\ensuremath{\mathbf{MD}^\vee}\xspace}

\newcommand{\PDor}{\ensuremath{\mathbf{PD}^\vee}\xspace}

\newcommand{\IPC}{\ensuremath{\mathbf{IPC}}\xspace}

\newcommand{\IK}{\ensuremath{\mathbf{IK}}\xspace}
\newcommand{\K}{\ensuremath{\mathbf{K}}\xspace}

\newcommand{\LL}{\ensuremath{\mathsf{L}}\xspace}

\newcommand{\KP}{\ensuremath{\mathbf{KP}}\xspace}

\newcommand{\MT}{\ensuremath{\mathbf{MT_0}}\xspace}
\newcommand{\MTz}{\ensuremath{\mathbf{MT_0^-}}\xspace}

\newcommand{\Inql}{\ensuremath{\mathbf{InqL}}\xspace} 

\newcommand{\PROP}{\ensuremath{\mathsf{Prop}}\xspace}

\newcommand{\FF}{\ensuremath{\mathfrak{F}}\xspace}
\newcommand{\MM}{\ensuremath{\mathfrak{M}}\xspace}

\newcommand{\dep}{\ensuremath{\mathop{\,=\!}}}
\newcommand{\cmor}{\ensuremath{\otimes}}
\newcommand{\ior}{\ensuremath{\vee}}

\renewcommand{\qed}{\hfill $\clubsuit$}

\newcommand{\ci}{\ensuremath{\wedge\textsf{I}}\xspace}

\newcommand{\ce}{\ensuremath{\wedge\textsf{E}}\xspace}
\newcommand{\sori}{\ensuremath{\sor\textsf{I}}\xspace}
\newcommand{\sore}{\ensuremath{\sor\textsf{E}}\xspace}
\newcommand{\sorwe}{\ensuremath{\sor\textsf{E}}\xspace}

\newcommand{\sors}{\ensuremath{\sor\textsf{Sub}}\xspace}
\newcommand{\bori}{\ensuremath{\bor\textsf{I}}\xspace}
\newcommand{\bore}{\ensuremath{\bor\textsf{E}}\xspace}

\newcommand{\com}{\ensuremath{\textsf{Com}}\xspace}
\newcommand{\ass}{\ensuremath{\textsf{Ass}}\xspace}
\newcommand{\dstr}{\ensuremath{\textsf{Dstr}}\xspace}

\newcommand{\boti}{\ensuremath{\neg\textsf{E}}\xspace}

\newcommand{\negi}{\ensuremath{\neg\textsf{I}}\xspace}
\newcommand{\dnege}{\ensuremath{\neg\neg\textsf{E}}\xspace}
\newcommand{\diamon}{\ensuremath{\Diamond\textsf{Mon}}\xspace}
\newcommand{\boxmon}{\ensuremath{\Box\textsf{Mon}}\xspace}
\newcommand{\diaior}{\ensuremath{\textsf{Distr}\Diamond\vee}\xspace}
\newcommand{\boxior}{\ensuremath{\textsf{Distr}\Box\vee}\xspace}
\newcommand{\depdf}{\ensuremath{\dep(\cdot)\textsf{df}}\xspace}
\newcommand{\boxdiainter}{\ensuremath{\textsf{Inter}\Box\Diamond}\xspace}
\newcommand{\diaboxinter}{\ensuremath{\textsf{Inter}\Diamond\Box}\xspace}

\newcommand{\impi}{\ensuremath{\to\textsf{I}}\xspace}
\newcommand{\impe}{\ensuremath{\to\textsf{E}}\xspace}

\newcommand{\depiz}{\ensuremath{\textsf{DepI}_0}\xspace}
\newcommand{\depez}{\ensuremath{\textsf{DepE}_0}\xspace}

\newcommand{\depik}{\ensuremath{\textsf{DepI}_k}\xspace}
\newcommand{\depek}{\ensuremath{\textsf{DepE}_k}\xspace}
\newcommand{\se}{\ensuremath{\textsf{SE}}\xspace}

\newcommand{\sor}{\ensuremath{\otimes}\xspace}
\newcommand{\bor}{\ensuremath{\vee}\xspace}

\startlocaldefs
\endlocaldefs

\begin{document}

\begin{frontmatter}

  %% TITLE OF YOUR PAPER%%%
  %% Words in title should begin with uppercase, except%%%
  %%% articles (and, the, a), conjunctions (and, for, nor, but),
  %%% prepositions (by, with, for, over, and so on)
  \title{Modal Dependence Logics: Axiomatizations and Model-theoretic Properties}
  %%% Choose a short title to be used as the running head on
  %%% odd-numbered pages%%%
  %%% Use this same title as "short title" when you submit MS to
  %%% EM%%%%
%  \runtitle{Modal Dependence Logics: Axiomatizations and Model-theoretic Properties
%  }

  \author{\fnms{Fan} %first name
    \snm{Yang}%last name
    \corref{}%to denote who is the corresponding author
    \ead[label=e1]{fan.yang.c@gmail.com}%author email, leave this [label] as is
%    \ead[label=u1,url]{http://ndjfl.nd.edu/}%%web page, leave this [label] as is
  }
  %%% ADDRESS NOTES
  %%% Dept listed first, University/company second, street or PO box
  %%% third%%%
  %%% U.S. Postal Service guidelines request no punctuation in the
  %%% street...country lines%%%
  %%% Country should be all uppercase%%%%
\vspace{-1.8\baselineskip}

  \address{%Department of Values, Technology and Innovation\\
  Delft University of Technology\\ 
  Jaffalaan 5, 2628 BX Delft\\
  The Netherlands\\
    \printead{e1}\\
    %\printead{u1}
     }%
 
  % \and \author{\fnms{???} \snm{???}\ead[label=e3]{???}}
  % \address{\printead{e3}} \affiliation{???}

  %%%% INSERT EACH AUTHOR'S FIRST INITIAL AND SURNAME%%%%
%  \runauthor{F. Yang}

\begin{abstract}
Modal dependence logics are modal logics defined on the basis of team semantics and have the downward closure property.
In this paper, we introduce sound and complete deduction systems for the major modal dependence logics, especially those with intuitionistic connectives in their languages. We also establish a concrete connection between team semantics and single-world semantics, and show that modal dependence logics can be interpreted as variants  of  intuitionistic modal logics.
\end{abstract}

\begin{keyword}[class=AMS]
  \kwd{03B45} \kwd{03B60} \kwd{03B55}
\end{keyword}

\begin{keyword}
  \kwd{dependence logic} \kwd{team semantics} \kwd{modal logic} \kwd{intuitionistic modal logic} \kwd{intermediate logics} 
\end{keyword}

\end{frontmatter}

%\newpage

\emph{Dependence logic} is a logical formalism, introduced by V\"{a}\"{a}n\"{a}nen \cite{Van07dl}, that captures the notion of \emph{dependence} in social and natural sciences. The modal version of the logic is called \emph{modal dependence logic} and was introduced in \cite{VaMDL08}. Modal dependence logic extends the usual modal logic by adding a new type of atomic formulas $\dep(p_1,\dots,p_1,q)$, called \emph{dependence atoms}, to express dependencies between propositions, and by lifting the usual single-world semantics to the so-called \emph{team semantics}, introduced by Hodges \cite{Hodges1997a,Hodges1997b}. Formulas of modal dependence logic are evaluated on sets of possible worlds of Kripke models, called \emph{teams}. Intuitively, a dependence atom $\dep(p_1,\dots,p_1,q)$ is true if within a team  the truth value of the proposition $q$ is \emph{functionally determined} by the truth values of the propositions $p_1, \dots, p_n$. 
%Closely related is the framework of \emph{inquisitive logic} \cite{InquiLog}, where team semantics is also (independently) adopted (see e.g. \cite{Ciardelli2015,VY_PD} for further discussions on the connection).

Research on modal dependence logic and its variants has been active in recent years. Basic model-theoretic properties of the logics were studied in e.g., \cite{sevenster09}, a van Benthem Theorem for the logics was proved  in \cite{KontinenMulleerSchnoorVollmer15}, and the frame definability of the logics was studied in \cite{SanoVirtema:2015,SanoVirtema:2016}. The expressive power and the relevant computation complexity problems of the logics were investigated extensively  in e.g., \cite{EHMMVV2013,eblo11,MID_mc,HLSV14,lovo10,MullerVollmer2013,sevenster09}. In this paper, we study two problems that received less attention in the literature, namely the axiomatization problem and the comparison between team semantics and the single-world  semantics.

For the axiomatization problem, Sano and Virtema gave in  \cite{SanoVirtema2014} Hilbert-style systems and label tableau calculi for modal dependence logic and its extended version, and Hannula defined in \cite{Hannula16entail} natural deduction systems for the same logics. 
However, these axiomatizations did not cover the modal dependence logics with intuitionistic connectives of team semantics, especially with  intuitionistic implication. The intuitionistic connectives of team semantics are crucial connectives of  \emph{inquisitive logic} \cite{InquiLog}, a closely related logic to dependence logic that adopts (independently) team semantics too (see e.g. \cite{Ciardelli2015,VY_PD} for further discussions on the connection). While inquisitive  modal logic has been axiomatized in \cite{Ciardelli_PhD}, the logic has different modalities and slightly different Kripke models than those of modal dependence logics.
%propositional base of modal dependence logic with intuitionistic connectives forms propositional \emph{inquisitive logic} \cite{InquiLog}, a closely related framework to dependence logic (see e.g. \cite{Ciardelli2015,VY_PD} for further discussions on the connection).
%The intuitionistic connectives of team semantics play a crucial role in inquisitive logic, which is  the propositional base of modal dependence logic.
 In this paper, we define Hilbert and natural deduction systems for modal dependence logic extended with intuitionistic disjunction and implication,  called \emph{full modal downward closed team logic} (\MT). These systems are extensions of the systems of Fischer Servi's intuitionistic modal logic (\IK) \cite{FS81_IK} and  inquisitive (propositional) logic \Inql \cite{InquiLog}. We also introduce  deduction systems for  modal intuitionistic dependence logic, modal dependence logic with intuitionistic disjunction, and (extended) modal dependence logic (\MD) as fragments or variants of the system of \MT.  We adopt the rules for (extended) dependence atoms introduced in \cite{VY_PD}. These rules are simpler than those in the systems of \cite{Hannula16entail,SanoVirtema2014}. We also point out that for the logic \MD, which does not have implication in its language, the deduction system that enjoys \emph{(weak) completeness} can have less rules than the system that enjoys \emph{strong completeness}. This interesting difference between the \emph{validity problem} (i.e., determining whether $\models\phi$) and the \emph{entailment problem} (i.e., determining whether $\phi\models\psi$) in \MD was noted also in \cite{Hannula16entail}. 

% for axiomatizing \emph{theoremhood} (i.e., $\vdash\phi$) can have less rules than the system for axiomatizing \emph{conquence relation} (i.e., $\psi\vdash\phi$)

% an interesting difference that was pointed out also in \cite{Hannula16entail}.

Our axiomatizations make  use of the disjunctive normal form of modal dependence logics, which is essentially known in the literature. We apply this normal form to prove a characterization theorem for flat formulas, the Interpolation Theorem and the Finite Model Property of modal dependence logics.

%we will define Hilbert systems for these intuitionistic variants of modal dependence logic. These systems are extensions of the systems of Fischer Servi's intuitionistic modal logic \IK \cite{FS81_IK} and inquisitive logic  \cite{InquiLog}. We also introduce a natural deduction system for (extended) modal dependence logic by extending the  system of propositional dependence logic given in \cite{VY_PD}, which has more refined rules for dependence atoms than those in the systems of \cite{Hannula16entail,SanoVirtema2014}. Moreover, we will discuss 
%%The arguments of our proofs of the Completeness Theorem of the deduction systems make essential use of 
%the disjunctive normal form of modal dependence logics, which is essentially known in the literature, and apply the normal form to prove the Interpolation Theorem and the finite model property of these logics. \todo{A unified system for comparing the logics}

For the second topic of this paper, it is well-known  that the team-based first-order dependence logic can be translated into the single-assignment-based existential second-order logic \cite{Van07dl,KontVan09}. In a similar fashion, we show in this paper that the team-based modal dependence logics can be  interpreted as  certain  single-world-based intermediate modal logics.
% team semantics is often considered (without further justification) in the literature as a natural generalization of the usual single-world (or single-assignment) semantics. In this paper, we propose a precise understanding of this commonly made claim. 
We  first provide a rigorous proof for a seemingly folklore observation in the field that clarifies the natural connection between team semantics and single-world semantics in the modal  case, namely,  the team semantics of modal dependence logics over a (classical) modal Kripke model $\mathfrak{M}$ coincides with the usual single-world semantics over an intuitionistic Kripke model whose domain consists of all teams of $\mathfrak{M}$ (i.e., the domain is the powerset of  $\mathfrak{M}$) and whose partial order is the superset relation between teams. The tensor (disjunction) connective of team semantics will be interpreted in this setting as a binary diamond modality under the single-world semantics, the idea of which is developed from \cite{AbVan09}, where tensor is understood as a multiplicative conjunction. 
%We will demonstrate that the powerset models construction also provides a natural justification for the perfect information semantic set games of (modal) dependence logic developed in \cite{Van07dl,VaMDL08}.

On the basis of the powerset models, we establish a comparison between modal dependence logics and  familiar single-world-based non-classical logics, especially intuitionistic modal logic and intermediate logics.
%we will build a more general connection between the team-based modal  logics and the single-world-based modal logics by showing 
We show that modal dependence logics  are complete (in the usual single-world semantics sense) with respect to a class of bi-relation or tri-relation intuitionistic Kripke models that generalise the powerset models. The bi-relation models are special bi-relation intuitionistic Kripke models of Fischer Servi's intuitionistic modal logic \IK, and the tri-relation intuitionistic Kripke models are endowed with an extra ternary relation interpreting the binary diamond that corresponds to the tensor. 
Our results generalise the results in \cite{ivano_msc} that inquisitive logic  \Inql  can be viewed as a variant of the Kreisel-Putnam intermediate logic (\KP) \cite{KrsPutnamLog57}, %the team semantics of \Inql coincides with the usual single-world semantics over negative saturated intuitionistic Kripke models, 
and \Inql is complete (in the usual single-world semantics sense) with respect to the class of negative intuitionistic Kripke models of \KP.

This paper is structured as follows. Section 1 recalls the basics of modal dependence logics. 
In particular, we sketch the standard translation from modal dependence logics into first-order dependence logics and derive the Compactness Theorem for modal dependence logics without intuitionistic implication as a corollary. 
In Section 2 we study the axiomatization problem for modal dependence logics, and also prove a few metalogical properties of the logics, including the Interpolation Theorem and the Finite Model Property.  Section 3 provides single-world semantics interpretation of modal dependence  logics. In Section 4 we make concluding remarks. 

Preliminary results of this paper were included in the author's dissertation \cite{Yang_dissertation}.

%We also prove that modal dependence logics with intuitionistic implication admit a disjunctive normal form whose similar 

%It was discussed in the literature (see e.g. \cite{lovo10,HLSV14}) that formulas of modal dependence logics without intuitionistic implication have disjunctive normal form. We prove that   modal dependence logics with intuitionistic implication admit a similar disjunctive normal form. Using this normal form, we prove the (uniform) Interpolation Theorem and the finite model property of these logics.

\section{Preliminaries}\label{sec:prel}

%\todo{need to introduce more background on team semantics}

In this section, we recall the basics of modal dependence logics, which are modal logics defined on the basis of team semantics. 

Though  team semantics is intended for the extension of (classical) modal logic obtained by adding dependence atoms, for the sake of comparison, we  start by defining the team semantics and fixing  notations for the usual (classical) modal logic.
%language of the usual (classical) modal logic and its team semantics.
Fix a set \PROP of propositional variables and denote its elements by $p,q,r,\dots$ (possibly with subscripts).
Formulas of \emph{(classical) modal logic}, also called \emph{classical (modal) formulas}, are defined recursively as:\vspace{-0.2\baselineskip}%built from the following grammar:
\[\alpha::=p\mid\bot\mid\neg\alpha\mid\alpha\wedge\alpha\mid \alpha\cmor\alpha\mid\alpha\to\alpha\mid \Box \alpha\mid\Diamond\alpha\]
where $\cmor$ (called \emph{tensor}) denotes the  disjunction of classical modal logic, and the implication $\to$ is called \emph{intuitionistic implication} for  reasons that will become clear in the sequel.

%Note that in this paper we denote the disjunction of classical modal logic by $\cmor$ (called \emph{tensor}). %and we define $\neg\alpha:=\alpha\to\bot$. %$\alpha\cmor\beta:=\neg\alpha\to\beta$ and $\Diamond\alpha:=\neg\Box\neg\alpha$. 

A (modal) \emph{Kripke frame} is a couple
$\mathfrak{F}=(W,R)$  consisting of a nonempty set $W$ and a binary relation $R\subseteq W\times W$. Elements of $W$ are called \emph{possible worlds} or \emph{nodes} or \emph{points}. A (modal) \emph{Kripke model} is a triple $\MM=(W,R,V)$ such that $(W,R)$ is a Kripke frame and $V:\PROP\to\wp(W)$ is a valuation function.  A set $X\subseteq W$ of possible worlds is called a \emph{team}. For any team $X$,  define 
\(R(X)=\{w\in W \mid \exists v\in X~\text{s.t. }vRw\}\)
and write $R(w)$ for $R(\{w\})$. %Elements in $R(w)$ are called \emph{successors} of $w$. 
A team $Y$ is called a \emph{successor team} of $X$, written $XRY$, if $Y\subseteq R(X)$ and $Y\cap R(w)\neq\emptyset$ for every $w\in X$.

%\begin{definition}
%We  define the notion of a classical  modal formula $\alpha$ being \emph{satisfied} in a Kripke model $\mathfrak{M}=(W,R,V)$ on a team $T\subseteq W$, denoted  $\mathfrak{M},T\models \phi$, as
%\begin{itemize}
%\item $\mathfrak{M},T\models\alpha$\quad iff \quad for all $w\in T$, $\mathfrak{M},w\models\alpha$ holds in the usual sense 
%\end{itemize}
%\end{definition}

\begin{definition}
We  define inductively the notion of a classical modal formula $\alpha$ being \emph{satisfied} in a Kripke model $\mathfrak{M}=(W,R,V)$ on a team $X\subseteq W$, denoted  $\mathfrak{M},X\models \alpha$,  as follows:
\begin{itemize}
\item $\mathfrak{M},X\models  p$ iff $X\subseteq V(p)$
\item $\mathfrak{M},X\models  \bot$ iff $X=\emptyset$
\item $\mathfrak{M},X\models  \neg  \phi$ iff $\mathfrak{M},\{w\}\not\models \phi$ for all $w\in X$.

\noindent In particular, $\mathfrak{M},X\models  \neg p$ iff $X\cap V(p)=\emptyset$
\item $\mathfrak{M},X\models  \phi\wedge\psi$ iff $\mathfrak{M},X\models  \phi$ and $\mathfrak{M},X\models  \psi$
\item $\mathfrak{M},X\models   \phi\cmor\psi$ iff there exist $Y,Z\subseteq X$ such that $X=Y\cup Z$, $\mathfrak{M},Y\models   \phi$ and $\mathfrak{M},Z\models  \psi$
\item $\mathfrak{M},X\models   \phi\to\psi$ iff for all $Y\subseteq X$,  $\mathfrak{M},Y\models \phi$ implies $\mathfrak{M},Y\models \psi$
\item $\mathfrak{M},X\models \Box \phi$ iff $\mathfrak{M},R(X)\models \phi$
\item $\mathfrak{M},X\models  \Diamond \phi$ iff there exists $Y\subseteq W$ such that $XRY$ and $\mathfrak{M},Y\models \phi$
\end{itemize}

%It is easy to see that the semantic clauses for negation, atomic negation, disjunction and diamond reduce to
%\begin{itemize}
%\item $\mathfrak{M},X\models  \neg  \phi$ iff $\mathfrak{M},\{w\}\not\models \phi$ for all $w\in X$
%\item $\mathfrak{M},X\models  \neg p$ iff $X\cap V(p)=\emptyset$
%\item $\mathfrak{M},X\models   \phi\cmor\psi$ iff there exist $Y,Z\subseteq X$ such that $X=Y\cup Z$, $\mathfrak{M},Y\models   \phi$ and $\mathfrak{M},Z\models  \psi$
%\item $\mathfrak{M},X\models  \Diamond \phi$ iff there exists $Y\subseteq W$ such that $XRY$ and $\mathfrak{M},Y\models \phi$
%\end{itemize}

For any Kripke model $\mathfrak{M}=(W,R,V)$, if $\mathfrak{M},X\models\phi$ holds for all $X\subseteq W$, then we say that $\phi$ is \emph{true} on $\mathfrak{M}$ and write $\mathfrak{M}\models\phi$. For any Kripke frame $\FF$, if $(\FF,V)\models\phi$ holds for all valuations $V$ on $\FF$, then we say that $\phi$ is \emph{valid on} $\FF$ and write $\FF\models\phi$. If $\FF\models\phi$ holds for all frames $\FF$, then we say that $\phi$ is \emph{valid} and write $\models\phi$.
We write $\Gamma\models\phi$ if for all Kripke models $\mathfrak{M}$ and all teams $X$,  $\mathfrak{M},X\models\gamma$ for all $\gamma\in \Gamma$ implies $\mathfrak{M}, X\models\phi$. We write simply $\phi\models\psi$ for $\{\phi\}\models\psi$, and write $\phi\equiv\psi$ if $\phi\models\psi$ and $\psi\models\phi$.
\end{definition}

%\noindent By the flatness property, $\phi\to\bot$ has the same semantics as the formula $\neg\phi$ obtained from 
%In view of this, lifting the single-possible world semantics to team semantics does not add new features to classical modal formulas. %It is only when dependence atoms or intuitionistic connectives are added the flatness property is violated and the environment of team semantics gives rise to interesting new properties. 

It is easy to check that classical modal formulas $\alpha$ satisfy the \emph{flatness property}:%, that is,
\begin{description}
\item[Flatness Property] 
$\begin{array}[t]{rl}
\mathfrak{M},X\models\alpha\!\!\!\!\!\!\!\!\!\!&\iff\mathfrak{M},\{w\}\models\alpha\text{ for all }w\in X\\
&\iff \mathfrak{M},w\models\alpha\text{ in the usual sense for all }w\in X\end{array}$
\end{description}
As a consequence, a few usual equivalences, such as the following ones, hold for classical formulas:\vspace{-0.35\baselineskip}
\[\neg\alpha\equiv\alpha\to\bot, \quad\alpha\sor\beta\equiv\neg\alpha\to\beta\quad\text{and}\quad\Diamond\alpha\equiv\neg\Box\neg\alpha.\vspace{-0.35\baselineskip}\]

%\begin{definition}\label{Hlibert_K}
Recall that the Hilbert-style system of classical modal logic \K consists of the following axioms and rules:
\begin{description}
\item[Axioms]\

%\vspace{-1.85\baselineskip}
\begin{enumerate}
%\addtolength{\itemindent}{0.2cm}
\item all axioms of classical propositional logic 
\item $\Box(\alpha\to \beta)\to(\Box \alpha\to\Box \beta)$
\item $\Diamond \alpha\leftrightarrow\neg\Box\neg  \alpha$
\end{enumerate}
%\begin{center}
%\begin{tabular}{C{0.05\linewidth}L{0.4\linewidth}C{0.05\linewidth}L{0.4\linewidth}}
%(i)& all axioms of \CPC & (iii)&  $\Diamond \phi\leftrightarrow\neg\Box\neg  \phi$\\
%(ii) &\textsf{K}: $\Box(\phi\to \psi)\to(\Box \phi\to\Box \psi)$& \\
%\end{tabular}
%\end{center}

%\vspace{-0.5\baselineskip}
\item[Rules]\

%\vspace{-1.85\baselineskip}

 \begin{enumerate}
\item Modus Ponens: $\alpha,\alpha\to\beta/\beta$
\item Necessitation: $\alpha/\Box\alpha$
\item Uniform Substitution: $\alpha/\alpha(\beta/p)$
\end{enumerate}
%\begin{center}
%\begin{tabular}{C{0.05\linewidth}L{0.4\linewidth}C{0.05\linewidth}L{0.4\linewidth}}
%(i)& Modus Ponens: $\phi,\phi\to\psi/\psi$ & (iii)& Necessitation: $\phi/\Box\phi$\\
%(ii) &Substituion: $\phi/\phi(\psi/p)$& \\
%\end{tabular}
%\end{center}
\end{description}
%\end{definition}
For classical formulas, the system of \K is sound and complete with respect to all Kripke frames in the sense of team semantics too. To see why, by the Completeness Theorem of \K with respect to the usual single-world semantics and the flatness property of classical formulas, we have
\begin{align*}
\vdash_{\K}\alpha&\iff\mathfrak{M},w\models\alpha\text{ for any model $\mathfrak{M}$ and any possible world $w$ in $\mathfrak{M}$}\nonumber\\
&\iff\mathfrak{M},\{w\}\models\alpha\text{ for any model $\mathfrak{M}$ and any possible world $w$ in $\mathfrak{M}$}\nonumber\\
&\iff\mathfrak{M},X\models\alpha\text{ for any model $\mathfrak{M}$ and any team $X$ of $\mathfrak{M}$}\nonumber\\
&\iff\models\alpha.
\end{align*}
As a consequence, we also have
\begin{equation}\label{Team2K}
\alpha\vdash_\K\beta\iff\vdash_\K\alpha\to\beta\iff\models\alpha\to\beta\iff\alpha\models\beta.
\end{equation}

In view of the flatness property, lifting the single-possible world semantics to team semantics does not add essentially new features to classical modal formulas. 
Let us now extend classical modal logic by adding the \emph{dependence atoms} $\dep(\alpha_1,\dots,\alpha_n,\beta)$ and \emph{intuitionistic disjunction} $\vee$ that violate the flatness property. 
We call the resulting logic \emph{full modal downward closed team logic} (\MT) and its language is defined formaly as follows:
\[\phi::=\alpha\mid\dep(\alpha_1,\dots,\alpha_n,\beta)\mid\phi\wedge\phi\mid\phi\cmor\phi\mid\phi\vee\phi\mid\phi\to\phi\mid \Box \phi\mid \Diamond \phi,\]
where $\alpha,\alpha_1,\dots,\alpha_n,\beta$ are classical modal formulas. We write $\neg\phi$ for $\phi\to\bot$.

%For a set $\Phi$ of formulas and two possible worlds $w$ and $u$, we write $w\sim_\Phi u$ if $\mathfrak{M},w\models\phi \iff \mathfrak{M},u\models\phi$ for every $\phi\in \Phi$, and write $w\sim_\phi u$ for $w\sim_{\{\phi\}}u$.
\begin{definition}
The satisfaction relation for \MT-formulas is defined as for classical modal formulas and additionally:
\begin{itemize}
\item $\mathfrak{M},X\models\, \dep(\alpha_1,\cdots,\alpha_n,\beta)$ iff for any $w,u\in X$, $\mathfrak{M},w\models \alpha_i\Leftrightarrow \mathfrak{M},u\models\alpha_i$ for all $1\leq i\leq n$ implies $\mathfrak{M},w\models \beta\Leftrightarrow \mathfrak{M},u\models\beta$.
\item $\mathfrak{M},X\models  \phi\vee\psi$ iff $\mathfrak{M},X\models  \phi$ or $\mathfrak{M},X\models  \psi$
\end{itemize}
%For any Kripke models $\mathfrak{M}=(W,R,V)$, if $\mathfrak{M},X\models\phi$ holds for all $X\subseteq W$, then we say that $\phi$ is \emph{true} on $\mathfrak{M}$ and write $\mathfrak{M}\models\phi$. If $\mathfrak{M}\models\phi$ holds for all Kripke models $\mathfrak{M}$, then we write $\models_{\MT}\phi$ (or simply $\models\phi$) and say that $\phi$ is \emph{valid}.  
\end{definition}

Immediately from the semantics it follows that \MT-formulas  have the downward closure property and the empty team property defined below:
\begin{description}
\item[Downward Closure Property] $[\,\mathfrak{M},X\models\phi\text{ and }Y\subseteq X\,]\Longrightarrow \mathfrak{M},Y\models\phi$
\item[Empty Team Property] $\mathfrak{M},\emptyset\models\phi$
\end{description}
%From this and the flatness property of classical modal formulas, we know that \MT validates exactly the same classical formulas $\alpha$ as classical modal logic \K, because
%\begin{align*}
%\models_{\K}\alpha&\iff\mathfrak{M},w\models\alpha\text{ for any model $\mathfrak{M}$ and any possible world $w$ in $\mathfrak{M}$}\\
%&\iff\mathfrak{M},\{w\}\models\alpha\text{ for any model $\mathfrak{M}$ and any possible world $w$ in $\mathfrak{M}$}\\
%&\iff\mathfrak{M},X\models\alpha\text{ for any model $\mathfrak{M}$ and any team $X$ in $\mathfrak{M}$}\\
%&\iff\models_{\MT}\alpha.
%\end{align*}

We have discussed that the classical fragment of \MT (i.e., the fragment consisting of classical formals only) behaves exactly as classical modal logic. \MT also inherits from classical modal logic many other nice properties, such as the preservation property under taking disjoint unions. For any two Kripke models $\MM=(W,R,V)$ and $\MM'=(W',R',V')$, their \index{disjoint union}\emph{disjoint union} $\MM\uplus\MM'=(W_0,R_0,V_0)$ is defined as 
\[W_0=W\uplus W', ~R_0=R\uplus R'\text{ and }V_0(p)=V(p)\uplus V'(p)\text{ for all }p\in \mathop{\rm{Prop}},\] 
where $\uplus$ takes the disjoint union of two sets. It can be proved by a routine argument that for any collection $\{\mathfrak{M}_i=(W_i,R_i,V_i)\mid i\in I\}$ of Kripke models, for every $i\in I$ and every $X\subseteq W_i$, %it holds that
 \begin{equation}\label{truth-preserve-operators}
 \mathfrak{M}_i,X\models\phi \iff\displaystyle\biguplus_{j\in I}\mathfrak{M}_j,X\models\phi.
 \end{equation}
From this it follows that that \MT has the \emph{disjunction property} with respect to the intuitionistic disjunction $\ior$, as shown below.

\begin{theorem}[Disjunction Property]\label{disjunct_prop}
If $\models\phi\vee\psi$, then $\models \phi$ or $\models \psi$.
\end{theorem}
\begin{proof}
Suppose $\mathfrak{M}_0,X\not\models \phi$ and $\mathfrak{M}_1,Y\not\models \psi$ for some models $\mathfrak{M}_0$ and $\mathfrak{M}_1$ and teams $X$ and $Y$. Let $\mathfrak{M}=\mathfrak{M}_0\uplus\mathfrak{M}_1$ and $Z=X\cup Y$. By (\ref{truth-preserve-operators}), we have 
\(\mathfrak{M},X\not\models\phi\text{ and }\mathfrak{M},Y\not\models\psi.\)
Hence, by the downward closure property, 
$\mathfrak{M},Z\not\models\phi$ and $\mathfrak{M},Z\not\models\psi$,
implying $\mathfrak{M},Z\not\models\phi\vee\psi$.
\end{proof}

In this paper, we also study some interesting fragments  of \MT (referred to as \emph{modal dependence logics}) defined by restricting the language as follows: %Let $\alpha,\alpha_1,\dots,\alpha_n,\beta$ be arbitrary classical modal formulas. 

\begin{itemize}

\item The language of modal dependence logic (\MD): 
\[\phi::=\alpha\mid \dep(p_1,\dots,p_n,q)\mid\phi\wedge\phi\mid\phi\cmor\phi\mid \Box \phi\mid \Diamond \phi\]
where $\alpha$ is an arbitrary classical formula defined recursively as 
\[\alpha::=p\mid \bot\mid \neg\alpha\mid\alpha\wedge\alpha\mid \alpha\sor\alpha\mid\Box\alpha\mid\Diamond\alpha\]

\item The language of extended modal dependence logic (\MDe): 
\[\phi::=\alpha\mid\dep(\alpha_1,\dots,\alpha_n,\beta)\mid\phi\wedge\phi\mid\phi\cmor\phi\mid \Box \phi\mid \Diamond \phi\]
where $\alpha,\alpha_1,\dots,\alpha_n,\beta$ are classical formulas defined as in the case of \MD. 
\item The language of modal dependence logic with intuitionistic disjunction (\MDor): 
\[\phi::=\alpha\mid \dep(\alpha_1,\dots,\alpha_n,\beta)\mid\phi\wedge\phi\mid\phi\cmor\phi\mid\phi\vee\phi\mid \Box \phi\mid \Diamond \phi\]
where $\alpha,\alpha_1,\dots,\alpha_n,\beta$ are classical formulas defined as in the case of \MD. 

\item The language of modal intuitionistic dependence logic (\MID):
\[\phi::=\alpha\mid\bot\mid\dep(\alpha_1,\dots,\alpha_n,\beta)\mid\phi\wedge\phi\mid\phi\ior\phi\mid\phi\to\phi\mid \Box \phi\mid \Diamond \phi\]
where $\alpha,\alpha_1,\dots,\alpha_n,\beta$ are classical formulas defined recursively as 
\[\alpha::=p\mid \bot\mid \alpha\wedge\alpha\mid\alpha\to\alpha\mid\Box\alpha\mid\Diamond\alpha\]
\end{itemize}
%\MD with dependence atoms $\dep(\alpha_1,\dots,\alpha_n,\beta)$ with classical arguments is known in the literature as \emph{extended modal dependence logic} (see e.g.,\cite{EHMMVV2013}). 
%We do include dependence atoms with classical arguments in the other logics we consider, because

%\todo{expressive power of \MDe}

Note that negation is taken to be a defined connective,  i.e., $\neg\phi:=\phi\to\bot$, in the modal dependence logics that have implication in their languages (such as \MT and \MID), while in the other logics (i.e., \MD, \MDe, \MDor, etc.) negation is a primitive connective that applies to classical formulas only. %The negation of modal dependence logics behaves as the classical negation over classical formulas, as $\neg\neg\alpha\leftrightarrow\alpha$ and $\alpha\sor\neg\alpha$ hold in these logics for classical formulas $\alpha$.

We leave it for the reader to check that in the presence of intuitionistic connectives dependence atoms are definable, as
\begin{equation}\label{dep_def}
\begin{split}\dep(\alpha_1,\dots,\alpha_k,\beta)&\equiv\big((\alpha_1\vee\neg \alpha_1)\wedge \dots\wedge (\alpha_k\vee\neg \alpha_k)\to (\beta\vee\neg \beta)\big)\\
&\equiv\bigotimes_{v\in{2}^{\{1,\dots,n\}}}(\alpha_1^{v(1)}\wedge\dots\wedge \alpha_n^{v(n)}\wedge (\beta\vee\neg \beta))
\end{split}
\end{equation}
where $2=\{0,1\}$, $\gamma^1=\gamma$ and $\gamma^0=\neg\gamma$. %we will mainly concentrate on the dependence atom-free fragment of \MT and \MID. 
The modality-free fragments of \MD, \MDor and \MID are called \emph{propositional dependence logic}, \emph{propositional dependence logic with intuitionistic disjunction} and \emph{propositional intuitionistic dependence logic}, respectively. These propositional logics were studied in \cite{VY_PD}. The modality and dependence atom-free fragment of \MID is in fact \emph{inquisitive (propositional) logic}, which was introduced by Ciardelli and Roelofsen in \cite{InquiLog} and commented also in the context of dependence logic in  \cite{Ciardelli2015,VY_PD,Yang_dissertation}. Ciardelli  studied and axiomatized in  \cite{Ciardelli_PhD} various inquisitive modal logic obtained from inquisitive propositional logic by adding different modalities than the $\Box$ and $\Diamond$ modalities we consider here in this paper. 

% \emph{Inquisitive modal logic} is obtained from inquisitive propositional logic by adding two modalities that are different from the $\Box$ and $\Diamond$ modalities we consider here in this paper. This logic has been axiomatized in \cite{Ciardelli_PhD} and we will not consider this logic in this paper. 

The language of \MDe differs from that of \MD only in that the dependence atoms of the latter have only propositional arguments. In the literature, the terminology \emph{dependence atoms} is often used for dependence atoms $\dep(p_1,\dots,p_n,q)$ with propositional arguments only, while dependence atoms $\dep(\alpha_1,\dots,\alpha_n,\beta)$ with classical arguments are often referred to as \emph{extended dependence atoms}. It is proved in \cite{EHMMVV2013,HLSV14} that \MD is strictly less expressive than \MDe, and the latter has the same expressive power as \MDor.

%  the language of  \MDe differs from that of  \MD only in that 

An interesting feature of  modal dependence logics is that they are \emph{not} closed under \emph{Uniform Substitution}. For instance,  $\models\neg\neg p\to p$ and $\models p\sor\neg p$, whereas $\not\models\neg\neg(p\vee\neg p)\to(p\vee\neg p)$ and $\not\models(p\vee\neg p)\sor\neg(p\vee\neg p)$. For this reason, none of the deduction systems of modal dependence logics to be introduced in this paper admits the uniform substitution rule. In fact, for these logics, \emph{substitution} (being a mapping from the set of well-formed formulas to the set itself that commutes with the atoms, connectives and modalities) is not even a well-defined notion, because, for instance, dependence atoms $\dep(\phi_1,\dots,\phi_n,\psi)$ with arbitrary arguments are not necessarily well-formed formulas of the logics. For more details on substitution in dependence logics, we refer the reader  to \cite{Ciardelli_PhD,IemhoffYang15}.

The well-known  standard translation from the usual (single-world-based) modal logic into first-order logic provides interesting insights into the usual modal logic. 
In particular, the Compactness Theorem of the usual modal logic is an immediate consequence of the translation. 
Without going into  further details we point out that  a similar translation from modal dependence logics into first-order dependence logics can be defined as follows:
%\setlength{\columnsep}{0.005\linewidth}
%\begin{multicols}{2}
\begin{itemize}
\item $ST_x(p):=Px$
\item $ST_x(\bot):=\bot$
\item $ST_x(\neg\alpha):=\neg ST_x(\alpha)$
%\item $ST_x(\alpha\to\beta):=\neg ST_x(\alpha)\sor ST_x(\beta)$
\item $ST_x(\dep(\alpha_1,\dots,\alpha_k,\beta)):=\dep(ST_x(\alpha_1),\dots,ST_x(\alpha_k),ST_x(\beta))$
\item $ST_x(\phi\wedge\psi):=ST_x(\phi)\wedge ST_x(\psi)$
\item $ST_x(\phi\sor\psi):=ST_x(\phi)\sor ST_x(\psi)$
\item $ST_x(\phi\vee\psi):=ST_x(\phi)\vee ST_x(\psi)$
%\item If $\alpha,\beta$ are classical formulas, $ST_x(\alpha\to\beta):=\neg ST_x(\alpha)\sor ST_x(\beta)$
%\item If $\phi,\psi$ are non-classical formulas, 
\item $ST_x(\phi\to\psi):=ST_x(\phi)\to ST_x(\psi)$
\item $ST_x(\Box\phi):=\forall y(\neg xRy\sor ST_y(\phi))$
\item $ST_x(\Diamond\phi):=\exists y(xRy\wedge ST_y(\phi))$
\end{itemize}
%\end{multicols}
\noindent where the extended dependence atom $\dep(ST_x(\alpha_1),\dots,ST_x(\alpha_k),ST_x(\beta))$ and the intuitionistic disjunction $\vee$ can both be defined in first-order dependence logic in terms of the other atoms and connectives\footnote{The team semantics of an extended dependence atom $\dep(\alpha_1,\dots,\alpha_k,\beta)$ with $\alpha_1,\dots,\alpha_k,\beta$ first-order formulas is defined as $M\models_X\dep(\alpha_1,\dots,\alpha_k,\beta)$ iff
\[\text{for all }s,s'\in X:~M\models_s\alpha_i\Leftrightarrow M\models_{s'}\alpha_i\text{ for all }1\leq i\leq k\text{ implies }M\models_s\beta\Leftrightarrow M\models_{s'}\beta.\] 
This atom can be defined in first-order dependence logic as $\dep(\alpha_1,\dots,\alpha_k,\beta):=$
\[\forall x\forall y(x=y)\vee\exists w_1\dots \exists w_k\exists u\exists v_0\exists v_1\big(\dep(w_1,\dots,w_k,u)\wedge\dep(v_0)\wedge \dep(v_1)\wedge (v_0\neq v_1)\]
\[\quad\quad\quad\quad\quad\quad\quad\wedge\bigwedge_{i=1}^k \big(\theta(w_i,v_0,v_1)\wedge \delta(w_i,\alpha_i,v_0,v_1)\big) \wedge\theta(u,v_0,v_1)\wedge\delta(u,\beta,v_0,v_1)\big),\]
where \(\theta(v,v_0,v_1):=(v=v_0)\otimes(v=v_1)\) and \(\delta(v,\gamma,v_0,v_1):=\big(\neg\gamma\sor(v= v_1)\big)\wedge \big(\gamma\sor (v=v_0) \big).\) The defining formula states intuitively that ``either the model in question has only one element in its domain (in which case $\dep(\alpha_1,\dots,\alpha_k,\beta)$ is trivially satisfied), or the team in question satisfies $\dep(w_1,\dots,w_k,u)$, where $w_i$ simulates $\alpha_i$ and $u$ simulates $\beta$".

The intuitionistic disjunction can be defined in first-order dependence logic as
\begin{align*}
\phi\vee\psi:=\big(&\exists x\exists y(x\neq y)\otimes(\phi\otimes\psi)\big)\wedge\\
&\big(\forall x\forall y(x=y)\otimes \exists u\exists v\big(\dep(u)\wedge \dep(v)\wedge (u\neq v)\wedge ((u=v)\otimes \phi)\wedge ((u\neq v)\otimes\psi )\big)\big),
\end{align*}
where the first conjunct of the defining formula deals with the case when the model has cardinality $1$, and the second conjunct deals with the other cases. %when the model has cardinality greater than $1$.
}.
%is an extended dependence atom that is defined in general as 
%\[\dep(\alpha_1,\dots,\alpha_k,\beta):=\forall x\forall y(x=y)\vee\exists w_1\dots \exists w_k\exists u\exists v_0\exists v_1\Big(\dep(w_1,\dots,w_k,u)\wedge\dep(v_0)\wedge \dep(v_1)\]
%\[\quad\quad\quad\quad\wedge (v_0\neq v_1)\wedge\bigwedge_{i=1}^k \big(\theta(w_i,v_0,v_1)\wedge \delta(w_i,\alpha_i,v_0,v_1)\big) \wedge\theta(u,v_0,v_1)\wedge\delta(u,\beta,v_0,v_1)\Big),\]
%\(\theta(v,v_0,v_1):=(v=v_0)\otimes(v=v_1)\) and \(\delta(v,\gamma,v_0,v_1):=\big(\neg\gamma\sor(v= v_1)\big)\wedge \big(\gamma\sor (v=v_0) \big).\) \todo{intuitionistic disjunction is definable in  dependence logic.}
A modal Kripke model $\mathfrak{M}=(W,R,V)$ can be associated in the usual manner with a first-order model $M=(W,R,\{P^M\}_{p\in \PROP})$  by interpreting the unary predicate symbols $P$ as $P^M=V(p)$. A routine argument shows that for every \MT-formula $\phi$, for every Kripke model $\mathfrak{M}=(W,R,V)$ and team $X\subseteq W$,
\begin{equation}\label{ST_truth}
\mathfrak{M},X\models\phi\iff M\models_{X_x}ST_x(\phi),
\end{equation}
where $X_x=\{\{(x,w)\}\mid w\in X\}$ is a first-order team with domain $\{x\}$.

Note that the formula $ST_x(\phi)$ is in general (equivalent to) a formula in the language of first-order dependence logic extended with intuitionistic implication, which is known to  have the same expressive power as full second-order logic \cite{Yang2011}, and thus not compact. But if $\phi$ does not contain intuitionistic implication,  $ST_x(\phi)$ is (equivalent to) a formula in the language of first-order (classical) dependence logic, which has the same expressive power as the compact existential second-order logic \cite{Van07dl}. We end this section by deriving the Compactness Theorem for modal dependence logics without intuitionistic implication as a corollary of the standard translation.

%But  if $\alpha$ and $\beta$ are classical formulas, then $ST_x(\alpha\to\beta)=ST_x(\alpha)\to ST_x(\beta)\equiv \neg ST_x(\alpha)\otimes ST_x(\beta)$, which is 

% Thus, $\phi$ does not contain intuitionistic implication $\to$ applied to non-classical formulas, then 

\begin{theorem}[Compactness]\label{compactness}
For any set $\Gamma\cup\{\phi\}$ of formulas in the language of \MD, \MDe or \MDor, if $\Gamma\models\phi$, then there exists a finite set $\Delta\subseteq \Gamma$ such that $\Delta\models\phi$.
\end{theorem}
\begin{proof}
By (\ref{ST_truth}), we have
\(\Gamma\models\phi\) iff \(\{ST_x(\gamma)\mid \gamma\in \Gamma\}\models ST_x(\phi).\)
%Note that the formulas on the right-hand-side are formulas of first-order dependence logic. 
The theorem then follows from the fact that  first-order dependence logic is compact \cite{Van07dl}.
\end{proof}

\section{Axiomatizations}\label{sec:axioma}

In this section, we introduce deduction systems for the modal dependence logics defined in the previous section and prove the completeness theorems. Our systems extend both the system of the propositional base of these logics, defined in \cite{Ciardelli2015,VY_PD}, and the system of Fischer Servi's intuitionistic modal logic \IK  \cite{FS81_IK}. The systems to be introduced in this section complement the previous axiomatizations \cite{Hannula16entail,SanoVirtema2014} for modal dependence logics in two respects. First, we provide axiomatizations for modal dependence logics with intuitionistic implication that have not yet been axiomatized before. Second, we incorporate from \cite{VY_PD} simpler rules for dependence atoms, and our systems for the different logics exhibit more uniformity.

We introduce in Section \ref{sec:mt} Hilbert-style and natural deduction systems for the full logic \MT and  \MID. 
%All the other logics introduced in the previous section  do not have implication in their languages. For these logics, 
For  the other logics introduced in the previous section that do not have implication in their languages, namely, \MDor, \MD and \MDe,
we  only define natural deduction systems (which, unlike Hilbert-style systems, do not necessarily use implications in the presentation). In Section \ref{sec:mdor}, we show that the implication-free fragment of the natural deduction system of \MT together with some additional rules for negation and modalities form a complete system for \MDor. The proofs of the completeness theorems in \Cref{sec:mt,sec:mdor} make heavy use of a disjunctive normal form for these logics which can essentially be found in the literature. In \Cref{sec:nf}, by using the normal form we prove a characterization theorem for flat formulas, the Interpolation Theorem and the Finite Model Property of modal dependence logics. In \Cref{sec:md}, we introduce the systems of \MD and \MDe as fragments of the system of \MT together with some additional rules for dependence atoms.

Note that for a logic  that has an implication $\to$ in its language, if the Deduction Theorem (i.e., $\Gamma,\phi\models\psi\iff\Gamma\models\phi\to\psi$) holds for this implication, the entailment problem and the validity problem can be reduced to each other, because
\[\phi\models\psi\iff \models\phi\to\psi,\]
and as a consequence, assuming the compactness of the logic, the strong completeness (i.e., $\Gamma\models\phi\iff \Gamma\vdash\phi$) of a deduction system of \LL is equivalent to its weak completeness (i.e., $\models\phi\iff \vdash\phi$). However, this is not in general true for logics without an implication in their languages. To address this subtle point, we will present the Completeness Theorems as 
\[\phi\models\psi\iff\phi\vdash\psi,\]
especially for the implication-free logics \MDor, \MD and \MDe.
We will see that for the systems of \MD and \MDe (which are compact by \Cref{compactness}), there is indeed a difference between the strong and the weak completeness: the systems for which the weak completeness holds can have less rules than the ones for which the strong completeness holds. This subtle difference was noted also in \cite{Hannula16entail} in a different system for  \MDe.

%and not definable in the logic. For modal, propositional and first-order dependence logics, the intuitionistic implication clearly admits the Deduction Theorem. But it is shown in \cite{Yang15} that 

%some of the modal dependence logics introduced in the previous section do not implication in their languages, namely, \MDor, \MD and \MDe.

%resolve the axiomatization problem of all modal dependence logics
%
%
%Our axiomatizations completes the solution to the axiomatization 

%There are two subtle points in the study of the axiomatization problem of modal dependence logics.

%We will also introduce a disjunctive normal form for these logics and show some applications of the normal form (including deriving a characterization theorem for flat and classical formulas, the interpolation theorem and the finite model property of modal dependence logics).

\subsection{\MT and \MID}\label{sec:mt}

In this subsection, we introduce sound and complete Hilbert-style and natural  deduction systems for \MID and \MT. We first define the Hilbert-style systems and prove the completeness theorems by an argument that makes essential use of the disjunctive normal form of the logics. The natural deduction systems will be defined at the end of the section and  their completeness follows from a similar argument.

 %and prove the completeness theorems. Hilbert-style and natural

%We first treat the logic \MID, which is the fragment of \MT that does not have the disjunction $\sor$ in its language. 
By Expression (\ref{dep_def}), dependence atoms are definable in \MID. The dependence atom-free fragment of \MID turns out to have  \emph{inquisitive logic} (\Inql) \cite{InquiLog} as its propositional base.
%The propositional fragment of \MIDz (i.e, the modality-free fragment of \MIDz) turns out to be \emph{inquisitive logic} (\Inql) introduced by Ciardelli and Roelofsen \cite{InquiLog}. 
Below we recall the Hilbert system of \Inql defined in \cite{Ciardelli_PhD,InquiLog}. We refer the reader to \cite{Ciardelli2015,VY_PD} for further discussion on the connection between  inquisitive logic and dependence logics.

\begin{definition}
The Hilbert-style system of inquisitive logic \Inql is as follows: 
\begin{description}
\item[Axioms]\

%\vspace{-1.85\baselineskip}

\begin{enumerate}
%\addtolength{\itemindent}{0.2cm}
\item all axiom schemes of intuitionistic propositional logic (\IPC), namely
\begin{enumerate}%[(1)]
%\addtolength{\itemindent}{0.1cm}
\item $\phi\to(\psi\to\phi)$ 
\item $(\phi\to(\psi\to\chi))\to((\phi\to\psi)\to(\phi\to\chi))$
\item $\phi\wedge\psi\to\phi$, $\phi\wedge\psi\to\psi$
\item $\phi \to (\chi \to (\phi \wedge \chi ))$
\item $\phi\to\phi\ior\psi$, $\psi\to\phi\ior\psi$
\item $(\phi\to\chi)\to((\psi\to\chi)\to(\phi\ior\psi\to\chi))$
\item $\bot\to\phi$
\end{enumerate}

%\begin{tabular}{ll}
%(1) $\phi\to(\psi\to\phi)$&(4) $\phi\to\phi\ior\psi$, $\psi\to\phi\ior\psi$\\
%(2) $(\phi\to(\psi\to\chi))\to((\phi\to\psi)\to(\phi\to\chi))$&(5) $(\phi\to\chi)\to((\psi\to\chi)\to(\phi\ior\psi\to\chi))$\\
%(3)  $\phi\wedge\psi\to\phi$, $\phi\wedge\psi\to\psi$&(6) $\bot\to\phi$
%\end{tabular}

%\setlength{\columnsep}{0.01\linewidth}
%\begin{multicols}{2}
% \begin{enumerate}[(1)]
%\item $\phi\to(\psi\to\phi)$ 
%\item $(\phi\to(\psi\to\chi))\to((\phi\to\psi)\to(\phi\to\chi))$
%\item $\phi\wedge\psi\to\phi$, $\phi\wedge\psi\to\psi$
%\item $\phi\to\phi\ior\psi$, $\psi\to\phi\ior\psi$
%\item $(\phi\to\chi)\to((\psi\to\chi)\to(\phi\ior\psi\to\chi))$
%\item $\bot\to\phi$
%\end{enumerate}
%\end{multicols}
\item $(\alpha\to(\phi\ior\psi))\to(\alpha\to\phi)\ior(\alpha\to\psi)$ whenever $\alpha$ is a classical formula
\item $\neg\neg\alpha\to\alpha$ whenever $\alpha$ is a classical formula
\end{enumerate}

%\vspace{-0.5\baselineskip}

\item[Rule]  Modus Ponens: $\phi,\phi\to\psi/\psi$
\end{description}
\end{definition}

In the original presentation of the system of \Inql as given in \cite{InquiLog}, axiom 2  is formulated (equivalently) as any substitution instance of the KP axiom
\[(\neg p\to(q\vee r))\to (\neg p\to q)\vee(\neg p\to r)\] 
of the Kreisel-Putnam intermediate logic \KP \cite{KrsPutnamLog57}. The system of \Inql is then \KP without uniform substitution rule together with the Double Negation Law $\neg\neg\alpha\to\alpha$ for classical formulas. 
%It was proved in \cite{InquiLog} that the set of theorems of \Inql equals the theorems of the negative variant of any intermediate logic between Maksimova's logic \ND and Medvedev's logic \MLo.

% In \cite{InquiLog} axiom 2 was originally formulated as the KP axiom 
% \[(\neg p\to(q\vee r))\to (\neg p\to q)\vee(\neg p\to r)\] 
% of Kreisel-Putnam logic \KP \cite{KrsPutnamLog57}, and axiom 3 is formulated as the atomic double negation law $\neg\neg p\to p$. In other words, the system of \Inql is \KP together with the atomic double negation law $\neg\neg p\to p$ (being not closed under uniform substitution). It can be shown that our presentation above is equivalent to the original presentation in  \cite{InquiLog} (cf. \cite{Ciardelli2015}).

Our Hilbert-style systems of \MT and \MID will be extensions of both the system of \Inql and the system of Fischer Servi's intuitionistic modal logic (\IK) \cite{FS81_IK}. We refer the reader to \cite{SimpsonPhDthesis} for further discussion on intuitionistic modal logic, and we only remark that \IK has intuitionistic propositional logic \IPC as its propositional base and adding the Law of Excluded Middle (i.e., $\phi\vee\neg \phi$) or the Double Negation Law (i.e., $\neg\neg\phi\to\phi$) to the logic gives rise to \emph{classical} modal logic. In the literature there are a few (equivalent) variants of the system of Fischer Servi's intuitionistic modal logic \IK. For the convenience of our argument, we  use the system  defined by Plokin and Stirling \cite{IML_PlotkinStirling86}, which we recall below. %Let us now recall the Hilbert-style system of \IK. 

\begin{definition}
The Hilbert-style system of Fischer Servi's intuitionistic modal logic \IK  consists of the following axioms and rules:
\begin{description}
\item[Axioms]\

\vspace{-0.8\baselineskip}

\begin{multicols}{2}
\begin{enumerate}
%\addtolength{\itemindent}{0.2cm}
\item all axioms of \IPC 
\item $\Box( \phi\to \psi)\to(\Box \phi\to\Box \psi)$
\item $\Box( \phi\to \psi)\to(\Diamond \phi\to\Diamond \psi)$ 
\item $\neg\Diamond\bot$
\item $\Diamond(\phi\vee \psi)\to(\Diamond \phi\vee\Diamond \psi)$
\item $(\Diamond \phi\to\Box \psi)\to\Box( \phi\to \psi)$
\end{enumerate}
\end{multicols}

\vspace{-0.8\baselineskip}
%\begin{center}
%\begin{tabular}{C{0.05\linewidth}L{0.4\linewidth}C{0.05\linewidth}L{0.4\linewidth}}
%(i)& all axioms of \CPC & (iii)&  $\Diamond \phi\leftrightarrow\neg\Box\neg  \phi$\\
%(ii) &\textsf{K}: $\Box(\phi\to \psi)\to(\Box \phi\to\Box \psi)$& \\
%\end{tabular}
%\end{center}
\item[Rules]\

%\vspace{-1.85\baselineskip}

 \begin{enumerate}
\item Modus Ponens: $\phi,\phi\to\psi/\psi$
\item Necessitation: $\phi/\Box\phi$
\item Uniform Substitution: $\phi/\phi(\psi/p)$
\end{enumerate}
\end{description}
%\begin{center}
%\begin{tabular}{C{0.05\linewidth}L{0.4\linewidth}C{0.05\linewidth}L{0.4\linewidth}}
%(i)& all axioms of \IPC & (iv)& $\neg\Diamond\bot$\\
%(ii) &\textsf{K}: $\Box( \phi\to \psi)\to(\Box \phi\to\Box \psi)$& (v)& $\Diamond(\phi\vee \psi)\to(\Diamond \phi\vee\Diamond \psi)$\\
%(iii)&  $\Box( \phi\to \psi)\to(\Diamond \phi\to\Diamond \psi)$ & (vi)& $(\Diamond \phi\to\Box \psi)\to\Box( \phi\to \psi)$\\
%\end{tabular}
%\end{center}
\end{definition}

Now, we present our Hilbert-style systems of \MT and \MID.

\begin{definition}\label{Hilbert_MT}
\begin{itemize}
%\addtolength{\itemindent}{-0.4cm}
\item The Hilbert-style system of \MT is defined as follows:
\begin{description}
\addtolength{\itemindent}{-0.3cm}
\item[Axioms]\

%\vspace{-1.3\baselineskip}
\begin{enumerate}
\addtolength{\itemindent}{-0.5cm}
\item all  axioms of inquisitive logic \Inql
\item all  axiom schemes of  \IK 
\item $\dep(\alpha_1,\dots,\alpha_k,\beta)\leftrightarrow\big((\alpha_1\vee\neg \alpha_1)\wedge \dots\wedge (\alpha_k\vee\neg \alpha_k)\to (\beta\vee\neg \beta)\big)$ 
%\item  $(\alpha\otimes\beta)\leftrightarrow(\neg\alpha\to\beta)$ whenever $\alpha$ and $\beta$ are classical modal formulas
\item $\phi\to\phi\otimes\psi$ 
\item $(\phi\to\alpha)\to\big((\psi\to\alpha)\to(\phi\otimes\psi\to\alpha)\big)$ whenever $\alpha$ is a classical formula
\item $(\phi\to\chi)\to\big((\psi\to\theta)\to(\phi\otimes\psi\to\chi\otimes\theta)\big)$
%\item $\alpha\otimes\alpha\to\alpha$ whenever $\alpha$ is a classical modal formula
\item $\phi\otimes\psi\to\psi\otimes\phi$
\item $\phi\otimes(\psi\otimes\chi)\to(\phi\otimes\psi)\otimes\chi$
\item $\phi\otimes(\psi\vee\chi)\to(\phi\otimes\psi)\vee(\phi\otimes\chi)$
\item $\neg\Box\alpha\to\Diamond \neg\alpha$ whenever $\alpha$ is a classical formula
\item $\Box( \phi\vee \psi)\to(\Box \phi\vee\Box \psi)$
\end{enumerate}

%\vspace{0.2\baselineskip}

\item[Rules]\

%\vspace{-1.3\baselineskip}

 \begin{enumerate}
 \addtolength{\itemindent}{-0.5cm}
\item Modus Ponens: $\phi,\phi\to\psi/\psi$
\item Necessitation: $\phi/\Box\phi$
\end{enumerate}
\end{description}

\item The Hilbert-style system of \MID consists of all of the axioms and rules of the above system of \MT except the axioms that involve $\otimes$ (i.e., axioms 4-9).
\end{itemize}

\end{definition}

%We write $\Gamma\vdash_{\LL}\phi$ if $\phi$ follows from the assumptions in the set $\Gamma$ of formulas in the system of \LL, and $\vdash_{\LL}\phi$ if $\emptyset\vdash\phi$. 
Hereafter within this section, we let \LL denote either \MT or \MID. The following proposition lists some interesting derivable clauses that will play a role in the sequel. 
\begin{proposition}\label{derivable_thm_MT_MID}
%Let $\LL\in\{\MT,\MID\}$. 
Let $\phi,\psi$ be  \LL-formulas, and $\alpha,\beta$ classical formulas.
\vspace{-0.5\baselineskip}
\begin{multicols}{2}
\begin{enumerate}
\item[(a)] $\phi\otimes(\psi\bor\chi)\dashv\vdash_{\MT}(\phi\otimes \psi)\bor(\phi\otimes \chi)$
\item[(b)] $\Box(\phi\vee\psi)\dashv\vdash_{\LL}\Box\phi\vee\Box\psi$
\item[(c)] $\Diamond(\phi\vee\psi)\dashv\vdash_{\LL}\Diamond\phi\vee\Diamond\psi$
\item[(d)] $\neg\neg\alpha\dashv\vdash_{\LL}\alpha$
%\item[(e)] $\neg\Box\alpha\dashv\vdash_{\LL}\Diamond \neg\alpha$
\item[(e)] $\Diamond\alpha\dashv\vdash_{\LL}\neg\Box\neg\alpha$
%\item[(e)]\label{tensor_comu} $\alpha\cmor\beta\vdash_\MT\beta\cmor\alpha$
%\item[(f)]\label{tensor_contr} $\alpha\cmor\alpha\vdash_\MT\alpha$ 
%\item[(g)]\label{tensor_sub} $\alpha\cmor\beta,\beta\to\gamma\vdash_\MT\alpha\cmor\gamma$
\item[(f)] $\vdash_{\MT}\alpha\cmor\neg\alpha$
%\item[(f)] $\alpha\otimes\beta\dashv\vdash_{\MT}\neg\alpha\to\beta$
%\item[(g)]\label{tensor_eli} $\vdash_\MT(\phi\to\alpha)\to\big((\psi\to\alpha)\to(\phi\otimes\psi\to\alpha)\big)$
%\item[(g)] $\alpha\cmor\alpha\vdash_\MT\alpha$ 
\end{enumerate}
\end{multicols}
\end{proposition}
\begin{proof}
A routine proof. In particular, in order to derive item (e), one may first derive from the \IK axioms $\Box( \phi\to \psi)\to(\Diamond \phi\to\Diamond \psi)$ and $\neg\Diamond\bot$ that $\vdash_{\IK}\Box\neg\phi\to\neg\Diamond\phi$ for arbitrary formulas $\phi$.
\end{proof}

%In fact, in the presence of the axiom $\neg\neg\alpha\to\alpha$, the axiom $\neg\Diamond\bot$ turns out to be derivable in turns of the other axioms:
%\[\begin{array}{lll}
%(1)&\neg\bot&(\IPC\text{ theorem})\\
%(2)&\Box\neg\bot&\text{(Necessitation)}\\
%(3)&\neg\neg\Box\neg\bot&\text{(\Cref{derivable_thm_MT_MID}(d))}\\
%(4)&\neg\Diamond\bot&\text{(\Cref{derivable_thm_MT_MID}(f))}
%\end{array}
%\]

Next, we prove the Soundness  Theorem for the systems of \MT and \MID.

\begin{theorem}[Soundness]
For any \LL-formulas $\phi$ and $\psi$, $\phi\vdash_{\LL}\psi$ $\Longrightarrow$ $\phi\models\psi$.
\end{theorem}
\begin{proof}
It suffices to show that all the axioms of \MT are valid and all the rules are sound. We only verify the validity of axiom 2 of \Inql and axiom 10. %and 7.

Axiom 2 of \Inql: We prove a slightly more general fact that 
\begin{equation}\label{KP_valid}
\theta\to(\phi\vee\psi)\models(\theta\to\phi)\vee(\theta\to\psi)\text{ whenever $\theta$ is flat}.
\end{equation}
Suppose $\mathfrak{M},X\not\models(\theta\to\phi)\vee(\theta\to\psi)$. Then $\mathfrak{M},X\not\models\theta\to\phi$ and $\mathfrak{M},X\not\models\theta\to\psi$. Thus, there exist $Y,Z\subseteq X$ such that 
\[\mathfrak{M},Y\models\theta,\quad \mathfrak{M},Z\models\theta, \quad\mathfrak{M},Y\not\models\phi\quad\text{and}\quad\mathfrak{M},Z\not\models\psi.\]
Since $\theta$ is flat and \LL has the downward closure property, we have 
\[\mathfrak{M},Y\cup Z\models\theta,\quad\mathfrak{M},Y\cup Z\not\models\phi\quad\text{and}\quad\mathfrak{M},Y\cup Z\not\models\psi.\]
Hence, $\mathfrak{M},X\not\models\theta\to(\phi\vee\psi)$.

Axiom 10: Suppose $\mathfrak{M},X\models \neg\Box\alpha$, where $\mathfrak{M}=(W,R,V)$ and $\alpha$ is a classical formula. Then, for any $w\in X$, we have  $\mathfrak{M},\{w\}\not\models \Box\alpha$, i.e. $\mathfrak{M},R(w)\not\models \alpha$. Since $\alpha$ is flat, there exists $v_w\in R(w)$ such that $\mathfrak{M},\{v_w\}\not\models \alpha$. Define
\(Y=\{v_w\in R(X)\mid w\in X\}.\)
Clearly, $XRY$ and $\mathfrak{M},Y\models \neg\alpha$. Hence, $\mathfrak{M},X\models \Diamond\neg\alpha$.
%Axiom 7: Suppose $\mathfrak{M},X\models \Box(\phi\vee\psi)$. Then, $\mathfrak{M},R(X)\models \phi\vee\psi$, which implies that $\mathfrak{M},R(X)\models \phi$ or $\mathfrak{M},R(X)\models \psi$. Thus, $\mathfrak{M},X\models\Box\phi$ or $\mathfrak{M},X\models\Box\psi$, implying $\mathfrak{M},X\models\Box\phi\vee\Box\psi$.
\end{proof}

To prove the Completeness Theorem for \MT and \MID, we will transform every formula into a formula in \emph{disjunctive normal form}:  
\begin{equation}\label{dnf_mt}
\alpha_1\vee\dots\vee\alpha_n
\end{equation}
where each $\alpha_i$ is a classical formula.
This  normal form is a generalization of a similar normal form for \Inql defined in \cite{Ciardelli_PhD,InquiLog}, and  similar disjunctive normal forms for modal dependence logics without intuitionistic implication were discussed  in the literature  with a slightly different presentation (see e.g., \cite{HLSV14,lovo10}).
 Let us now introduce our disjunctive normal form for \MT as a recursively defined translation $\tau(\phi)$ for every formula $\phi$ of \MT:
\begin{description}
\item[Base case]\ 
\begin{itemize}
\item $\tau(\alpha)=\alpha$ when $\alpha$ is a classical formula
\item $\tau(\dep(\alpha_1,\dots,\alpha_k,\beta)):= \tau((\alpha_1\vee\neg\alpha_1)\wedge\dots\wedge (\alpha_k\vee\neg\alpha_k)\to (\beta\vee\neg\beta))$
\end{itemize}

\item[Induction step]\ 

Assume 
$\tau(\psi)=\alpha_1\vee\dots\vee\alpha_n$ and $\tau(\chi)=\beta_1\vee\dots\vee\beta_m$, where $\alpha_i$ and $\beta_j$ are classical formulas. 
\begin{itemize}
\item $\tau(\psi\vee\chi) := \tau(\psi)\vee\tau(\chi)$
\item $\tau(\psi\wedge\chi) := \bigvee\{\alpha_i\wedge\beta_j)\mid 1\leq i\leq n,~1\leq j\leq m\}$
\item $\tau(\psi\otimes\chi) := \bigvee\{\alpha_i\otimes\beta_j\mid 1\leq i\leq n,1\leq j\leq m\}$
\item $\tau(\psi\to\chi) = \bigvee\{\bigwedge_{i=1}^n(\alpha_{i}\to\beta_{f(i)})\mid f:\{1,\dots,n\}\to\{1,\dots,m\}\}$
\item $\tau(\Diamond\psi) := \bigvee\{\Diamond\alpha_i\mid 1\leq i\leq n\}$
\item $\tau(\Box\psi) := \bigvee\{\Box\alpha_i\mid 1\leq i\leq n\}$
\end{itemize}
\end{description}
The disjunctive normal form for \MID is defined the same way as above except that  \MID does not have the connective $\otimes$ in its language. In the next theorem, we show that every formula is provably equivalent to its disjunctive normal form.

\begin{theorem}[Normal Form]\label{dn_form_pid}
For any \LL-formula $\phi$, we have $\vdash_{\LL}\phi\leftrightarrow \tau(\phi)$.
\end{theorem}
\begin{proof}
We only give the proof for \MT, from which the \MID case follows. %immediately from the \MT case.

We proceed by induction on $\phi$. The base case is trivial. For the induction step, the cases $\phi=\psi\vee\chi$ and $\phi=\psi\wedge\chi$ follow immediately from the induction hypothesis and \IPC axioms. The cases $\phi=\Box\psi$ and $\phi=\Diamond\psi$ follow from \Cref{derivable_thm_MT_MID}(b)(c).
%the fact that $\vdash_{\MT}\Box(\psi\ior\chi)\leftrightarrow(\Box\psi\ior\Box\chi)$ and $\vdash_{\IK}\Diamond(\psi\ior\chi)\leftrightarrow(\Diamond\psi\ior\Diamond\chi)$.

If $\phi=\psi\cmor\chi$, then by the induction hypothesis, we derive in the system of \MT using axiom 6 and \Cref{derivable_thm_MT_MID}(a) that 
\[\psi\cmor\chi\dashv\vdash\big(\bigvee_{i=1}^n\alpha_i\big)\cmor\big(\bigvee_{j=1}^m\beta_j\big)\dashv\vdash\bigvee_{i=1}^n\bigvee_{j=1}^m(\alpha_i\cmor\beta_j).%\dashv\vdash\bigvee_{i=1}^n\bigvee_{j=1}^m(\neg\alpha_i\to\beta_j).
\]

If $\phi=\psi\to\chi$, then by the induction hypothesis, we derive in the system of \MT using axiom 2 of \Inql and \IPC axioms that 
\begin{align*}
\psi\to\chi&\dashv\vdash\big(\bigvee_{i=1}^n\alpha_i\big)\to\big(\bigvee_{j=1}^m\beta_j\big)\dashv\vdash\bigwedge_{i=1}^n\big(\alpha_i\to \bigvee_{j=1}^m\beta_j\big)\dashv\vdash\bigwedge_{i=1}^n\bigvee_{j=1}^m(\alpha_i\to\beta_j)\\
&\dashv\vdash \bigvee\{\bigwedge_{i=1}^n(\alpha_{i}\to\beta_{f(i)})\mid f:\{1,\dots,n\}\to\{1,\dots,m\}\}.
\end{align*}
\end{proof}

%We leave it for the reader to verify that the system is sound and the Deduction Theorem holds: 
%\begin{description}
%\item[(Deduction Theorem)] $\Gamma,\phi\vdash_\MT\psi$ iff $\Gamma\vdash_\MT\phi\to\psi$ 
%\end{description}
%The axioms $\neg\neg\alpha\to\alpha$ and $(\alpha\cmor\beta)\leftrightarrow(\neg\alpha\to\beta)$ imply that over classical formulas the disjunction $\cmor$ behaves as the classical disjunction. In particular, the following lemma holds.

Each disjunct $\alpha_i$ in the disjunctive normal form $\bigvee_{i\in I}\alpha_i$  is a classical formula. We now show that \MT and \MID derive the same entailment relation as \K does.

%prove exactly the same set of classical modal formulas as \K does.

\begin{lemma}\label{class_K_syntax_equiv}
If $\alpha$ and $\beta$ are classical formulas, then $\alpha\vdash_{\K}\beta\iff\alpha\vdash_{\LL}\beta$.
\end{lemma}
\begin{proof}
For the direction ``$\Longleftarrow$'', suppose $\alpha\vdash_{\LL}\beta$. By the Soundness Theorem, we have $\alpha\models\beta$, which by Expression (\ref{Team2K}) from \Cref{sec:prel} implies $\alpha\vdash_{\K}\beta$.

For the direction ``$\Longrightarrow$", by \Cref{derivable_thm_MT_MID}(e) and by inspecting the axioms and rules, we see easily that restricted to classical formulas, the systems of \MT and \MID admit all \K rules and axioms, including Uniform Substitution rule, the axiom $\Diamond\alpha\leftrightarrow\neg\Box\neg\alpha$,  the double negation axiom $\neg\neg\alpha\to\alpha$ and all classical axioms of disjunction with respect to tensor $\otimes$. This implies that any classical entailment relation $\alpha\vdash\beta$ that is derivable in \K is also derivable (by the same derivation) in \LL.
\end{proof}

We are now ready to prove the Completeness Theorem for the two systems.%\todo{Mention the difference between this and $\vdash$ completeness}

\begin{theorem}[Completeness]\label{completeness_MT}
For any \LL-formulas $\phi$ and $\psi$, $\phi\models\psi\Longrightarrow\phi\vdash_{\LL}\psi$.
\end{theorem}
\begin{proof} 
%It suffices to show the direction ``$\Longrightarrow$''. 
Suppose $\phi\models\psi$. By Lemma \ref{dn_form_pid},
$\phi\dashv\vdash_{\LL}\alpha_1\vee\dots\vee\alpha_k$ and $\psi\dashv\vdash_{\LL}\beta_1\vee\dots\vee\beta_m$, for some classical formulas $\alpha_i$ and $\beta_j$.
By the Soundness Theorem, we have $\alpha_1\vee\dots\vee\alpha_k\models\beta_1\vee\dots\vee\beta_m$, which implies $\alpha_i\models\beta_1\vee\dots\vee\beta_m$ for each $1\leq i\leq k$. Since $\alpha_i$ is flat, it follows from Expression (\ref{KP_valid}) that $\alpha_i\models\beta_{j_i}$ for some $1\leq j_i\leq m$. We then derive by applying Expression (\ref{Team2K}) from \Cref{sec:prel} that $\alpha_i\vdash_{\K}\beta_{j_i}$, which yields $\alpha_i\vdash_{\LL}\beta_{j_i}$ by Lemma \ref{class_K_syntax_equiv}. Hence, $\alpha_1\vee\dots\vee\alpha_k\vdash_{\LL}\beta_1\vee\dots\vee\beta_m$, which gives $\phi\vdash_\LL\psi$.
\end{proof}

%For the purpose of comparison, we end this section with presenting natural deduction systems for \MT and \MID. %We leave it for the reader to verify that the systems are sound and complete (the completeness theorem can be proved by a similar argument as above).

Having proved the Completeness Theorem for the Hilbert-style systems of \MT and \MID, we now present also natural deduction systems of the logics.

%We will discuss further in Section 4 that \MIDz has a close connection with Fischer Servi's intuitionistic modal logic \IK \cite{FS81_IK}.

\begin{table}[t]
\begin{center}
{\normalsize
\begin{tabular}{|C{\linewidth}|}
%\hline
%\multicolumn{1}{l}{Rules for Intuitionistic Implication:}\\
\multicolumn{1}{c}{\small(a)}\\
\hline
\AxiomC{}\noLine\UnaryInfC{$[\phi]$}\noLine \UnaryInfC{$\vdots$}\noLine \UnaryInfC{$\psi$}\RightLabel{\impi}\UnaryInfC{$\phi\to\psi$}\noLine\UnaryInfC{}\DisplayProof
\quad\AxiomC{$\phi$}\AxiomC{$\phi\to\psi$}\RightLabel{\impe}\BinaryInfC{$\psi$}\DisplayProof
\quad\AxiomC{$\alpha\to(\phi\vee\psi)$}\RightLabel{\textsf{Split}}\UnaryInfC{$(\alpha\to\phi)\vee(\alpha\to\psi)$}\DisplayProof\\%\hline
%\multicolumn{1}{c}{}\\
%\multicolumn{1}{l}{Rules of Interactions of Modalities and Intuitionistic Implication:}\\\hline%\hline
\AxiomC{}\noLine\UnaryInfC{$\Box(\phi\to\psi)$}
\UnaryInfC{$\Box\phi\to\Box\psi$}\noLine\UnaryInfC{}\DisplayProof 
\quad\AxiomC{}\noLine\UnaryInfC{$\Box(\phi\to\psi)$}
\UnaryInfC{$\Diamond\phi\to\Diamond\psi$}\noLine\UnaryInfC{}\DisplayProof 
\quad\AxiomC{}\noLine\UnaryInfC{$\Diamond\phi\to\Box\psi$}
\UnaryInfC{$\Box(\phi\to\psi)$}\noLine\UnaryInfC{}\DisplayProof 
\\\hline
\end{tabular}
%\caption{Rules for implication}
}
%\end{center}
%\caption{Rules for implication. Hereafter in the tables, Greek letters, such as $\alpha,\beta,\alpha_1,\dots,\alpha_k\dots$, are variables for arbitrary classical formulas.\vspace{-2\baselineskip}}
%\label{tab:imp}
%\end{table}%
%
%\begin{table}[t]
%\begin{center}
{\normalsize
\begin{tabular}{|C{0.67\linewidth}C{0.3\linewidth}|}
%\hline
%\multicolumn{1}{l}{Monotonicity Rule for Box}\\\hline%
%\hline
%\multicolumn{1}{l}{Rules for Conjunction:}\\
\multicolumn{2}{c}{\small(b)}\\
\hline
\multicolumn{2}{|c|}{\AxiomC{}\noLine\UnaryInfC{$\phi$} \AxiomC{$\psi$}\RightLabel{\ci}\BinaryInfC{$\phi\wedge\psi$}\noLine\UnaryInfC{}\DisplayProof 
~\AxiomC{}\noLine\UnaryInfC{$\phi\wedge\psi$} \RightLabel{\ce}\UnaryInfC{$\phi$}\noLine\UnaryInfC{} \DisplayProof
~\AxiomC{}\noLine\UnaryInfC{$\phi\wedge\psi$} \RightLabel{\ce}\UnaryInfC{$\psi$}\noLine\UnaryInfC{}\DisplayProof
~\AxiomC{}\noLine\UnaryInfC{$\neg\neg\alpha$}\RightLabel{\dnege}\UnaryInfC{$\alpha$}\DisplayProof
~\AxiomC{}\noLine\UnaryInfC{$\bot$}\RightLabel{\textsf{Ex falso}}\UnaryInfC{$\phi$}\noLine\UnaryInfC{}\DisplayProof}\\
%\hline
%\multicolumn{1}{c}{}\\
%\multicolumn{1}{l}{Rules for Negation and contradiction:}\\\hline%\hline
%\multicolumn{2}{|c|}{\AxiomC{}\noLine\UnaryInfC{$\neg\neg\alpha$}\RightLabel{\dnege}\UnaryInfC{$\alpha$}\DisplayProof\quad
%\AxiomC{}\noLine\UnaryInfC{$\bot$}\RightLabel{\textsf{Ex falso}}\UnaryInfC{$\phi$}\noLine\UnaryInfC{}\DisplayProof}
%\\
%\hline
%\multicolumn{1}{c}{}\\
%\multicolumn{1}{l}{Rule for modalities}\\\hline%\hline
\AxiomC{}\noLine\UnaryInfC{$[\phi_1]$}
 \AxiomC{$\dots$}
 \AxiomC{$[\phi_k]$}
\noLine
\TrinaryInfC{$D^\ast$}
\branchDeduce
\DeduceC{$\psi$}
\AxiomC{$\Box\phi_1$} \AxiomC{\!\!\!\!\!$\dots$\!\!\!\!\!}\AxiomC{$\Box\phi_k$}\RightLabel{\boxmon}\QuaternaryInfC{$\Box\psi$}\noLine\UnaryInfC{}\DisplayProof 
&\AxiomC{}\noLine\UnaryInfC{$\neg\Box\alpha$}\RightLabel{\boxdiainter}\UnaryInfC{$\Diamond\neg\alpha$}\noLine\UnaryInfC{}\DisplayProof\\
%\hline
%~\\
\multicolumn{2}{|l|}{($\ast$) The derivation $D$ is assumed to have no undischarged assumptions.}\\
\hline
\end{tabular}
}
%\end{center}
%\caption{Rules for conjunction, negation, contradiction and modalities}
%\label{tab:conj}
%\end{table}%
%
%\begin{table}[t]
%\begin{center}
{\normalsize
\begin{tabular}{|C{\linewidth}|}
%\hline
%\multicolumn{1}{l}{Rules for Intuitionistic Disjunction:}\\\hline%
\multicolumn{1}{c}{\small(c)}\\
\hline
\AxiomC{$\phi$} \RightLabel{\bori}\UnaryInfC{$\phi\bor\psi$} \DisplayProof 
\quad\AxiomC{$\psi$} \RightLabel{\bori}\UnaryInfC{$\phi\bor\psi$}  \DisplayProof\quad
 \AxiomC{}\noLine\UnaryInfC{$[\phi]$}\noLine\UnaryInfC{$\vdots$}\noLine\UnaryInfC{$\chi$} \AxiomC{}\noLine\UnaryInfC{$[\psi]$}\noLine\UnaryInfC{$\vdots$}\noLine\UnaryInfC{$\chi$} \AxiomC{$\phi\bor\psi$}\RightLabel{\bore}\TrinaryInfC{$\chi$}\noLine\UnaryInfC{} \DisplayProof\\
%\hline
%\multicolumn{1}{c}{}\\
%\multicolumn{1}{l}{Distributive Rules for Intuitionistic Disjunction:}\\\hline%\hline
\AxiomC{}\noLine\UnaryInfC{$\phi\otimes(\psi\bor\chi)$} \RightLabel{$\dstr\sor\bor$}\UnaryInfC{$(\phi\otimes \psi)\bor(\phi\otimes \chi)$}\noLine\UnaryInfC{} \DisplayProof\quad\AxiomC{}\noLine\UnaryInfC{$\Diamond(\phi\vee\psi)$} \RightLabel{\diaior}\UnaryInfC{$\Diamond\phi\vee\Diamond\psi$}\noLine\UnaryInfC{}\DisplayProof\quad\AxiomC{}\noLine\UnaryInfC{$\Box(\phi\vee\psi)$} \RightLabel{\boxior}\UnaryInfC{$\Box\phi\vee\Box\psi$}\noLine\UnaryInfC{}\DisplayProof\\
\hline
\end{tabular}
}
%\end{center}
%\caption{Rules for intuitionistic disjunction}
%\label{tab:ior}
%\end{table}%
%
%\begin{table}[b!]
%\begin{center}
{\normalsize
\begin{tabular}{|C{0.2\linewidth}C{0.77\linewidth}|}
%\hline
%\multicolumn{1}{l}{Rules for Tensor:}\\
\multicolumn{2}{c}{\small(d)}\\
\hline
%\AxiomC{$\phi\sor\bot$}\RightLabel{\bote}\UnaryInfC{$\phi$}\DisplayProof 
\AxiomC{}\noLine\UnaryInfC{$\phi$} \RightLabel{\sori}\UnaryInfC{$\phi\sor\psi$} \DisplayProof&
\multirow{3}{*}{\AxiomC{}\noLine\UnaryInfC{$[\phi]$}\noLine\UnaryInfC{$\vdots$}\noLine\UnaryInfC{$\alpha$} \AxiomC{}\noLine\UnaryInfC{$[\psi]$}\noLine\UnaryInfC{$\vdots$}\noLine\UnaryInfC{$\alpha$}\AxiomC{$\phi\sor\psi$}  \RightLabel{\sorwe}\TrinaryInfC{$\alpha$}\DisplayProof
\quad\AxiomC{}\noLine\UnaryInfC{$[\psi]$}\noLine\UnaryInfC{$\vdots$}\noLine\UnaryInfC{$\chi$} \AxiomC{$\phi\sor\psi$} \RightLabel{\sors}\BinaryInfC{$\phi\sor\chi$}\DisplayProof}\\
~~\AxiomC{}\noLine\UnaryInfC{}\noLine\UnaryInfC{$\phi\sor\psi$} \RightLabel{$\com\sor$}\UnaryInfC{$\psi\sor\phi$}\noLine\UnaryInfC{} \DisplayProof&\\
~\AxiomC{}\noLine\UnaryInfC{$\phi\sor(\psi\sor \chi)$} \RightLabel{$\ass\sor$}\UnaryInfC{$(\phi\sor\psi)\sor\chi$}\noLine\UnaryInfC{} \DisplayProof&\\
\hline
\end{tabular}
}
%\end{center}
%\caption{Rules for tensor}
%\label{tab:tensor}
%\end{table}%
%\begin{table}[t]
%\begin{center}
%{\normalsize
%\begin{tabular}{|C{\linewidth}|}
%%\hline
%%\multicolumn{1}{l}{Monotonicity Rule for Box}\\\hline%
%%\hline
%%\multicolumn{1}{l}{Rules for Negation and contradiction:}\\\hline%\hline
%%\AxiomC{}\noLine\UnaryInfC{$\neg\neg\alpha$}\RightLabel{\dnege}\UnaryInfC{$\alpha$}\DisplayProof\quad
%%\AxiomC{}\noLine\UnaryInfC{$\bot$}\RightLabel{\textsf{Ex falso}}\UnaryInfC{$\phi$}\noLine\UnaryInfC{}\DisplayProof
%%\\\hline
%\multicolumn{1}{l}{Monotonicity Rule for Box}\\\hline%\hline
%\AxiomC{}\noLine\UnaryInfC{$[\phi_1]$}
% \AxiomC{$\dots$}
% \AxiomC{$[\phi_k]$}
%\noLine
%\TrinaryInfC{$D$}
%\branchDeduce
%\DeduceC{$\psi$}
%\AxiomC{$\Box\phi_1$} \AxiomC{\!\!\!\!\!$\dots$\!\!\!\!\!}\AxiomC{$\Box\phi_k$}\RightLabel{\boxmon}\QuaternaryInfC{$\Box\psi$}\noLine\UnaryInfC{}\DisplayProof 
%\quad\AxiomC{}\noLine\UnaryInfC{$\neg\Box\alpha$}\RightLabel{\boxdiainter}\UnaryInfC{$\Diamond\neg\alpha$}\noLine\UnaryInfC{}\DisplayProof\\
%%\hline
%%~\\
%($\ast$) where the derivation $D$  has no undischarged assumptions\\
%\hline
%\end{tabular}
%}
%\end{center}
%\caption{Rules for modalities, negation and contradiction}
%\label{tab:imp}
%\end{table}%
%\begin{table}[t]
%\begin{center}
{\normalsize
\begin{tabular}{|C{\linewidth}|}
%\hline
%\multicolumn{1}{l}{Monotonicity Rule for Box}\\\hline%
\multicolumn{1}{c}{\small(e)}\\
\hline
\AxiomC{}\noLine\UnaryInfC{$\dep(\alpha_1,\dots,\alpha_k,\beta)$}\doubleLine \RightLabel{\depdf}\UnaryInfC{$\displaystyle\bigotimes_{v\in {2}^{\{1,\dots,k\}}}(\alpha_1^{v(1)}\wedge\dots\wedge \alpha_k^{v(k)}\wedge (\beta\vee\neg \beta))$}\noLine\UnaryInfC{}\DisplayProof\\\hline
\end{tabular}
}
\end{center}
\caption{Rules of the  system of \MT. Hereafter in the tables $\alpha,\beta,\alpha_1,\dots,\alpha_k\dots$ range over classical formulas, and $\alpha_i^1:=\alpha_i$ and $\alpha_i^0:=\neg\alpha_i$ for each $i$.\vspace{-1.6\baselineskip}}
%\caption{Rules for dependence atom. For each index $i$, $\alpha_i^1:=\alpha_i$ and $\alpha_i^0:=\neg\alpha_i$.}
\label{tab:MT}
%\label{tab:dep}
\end{table}%

\begin{definition}\label{MT_nds_df}
\begin{itemize}
\item The natural deduction system of \MT is defined as follows:
\begin{description}
\item[Axiom]~ 
\AxiomC{}\noLine\UnaryInfC{}\RightLabel{\textsf{Ax}}\UnaryInfC{$\neg\Diamond\bot$}\noLine\UnaryInfC{}\DisplayProof
\item[Rules] All rules in \Cref{tab:MT}.
\end{description}
\item The natural deduction system  \MID consists of all the axioms and rules of the system of \MT except for the rules that involve $\otimes$, i.e., the rules in \Cref{tab:MT}(d) and the first distributive rule in \Cref{tab:MT}(c).
\end{itemize}
\end{definition}

The rules for the propositional base of \MT and \MID are adapted from those introduced for \PDor in \cite{VY_PD} and for $\mathsf{QD}_{\mathsf{P}}$ in \cite{Ciardelli2015}, and the rules for modalities are obvious translations of the axioms for the logics in the Hilbert-style systems of \Cref{Hilbert_MT}. The rule \boxmon with empty assumption corresponds to the Hilbert-style Necessitation Rule. It is easy to verify that the systems defined in \Cref{MT_nds_df} are sound. The completeness of the systems can be proved by a similar argument (via the disjunctive normal form) to what we presented above. We will not provide the proof here. However, in the next section, we will prove that a sound and complete natural deduction system of the logic \MDor (a fragment of \MT) can be obtained by dropping the inapplicable rules in \Cref{MT_nds_df} and adding certain additional rules.

\subsection{\MDor}\label{sec:mdor}

In this section, we introduce a sound and complete natural deduction system for  $\MD^\vee$, the implication-free  fragment of \MT. %\todo{no implication, therefore move to natural deduction}%Our system is an extension of the system of propositional dependence logics defined in \cite{VY_PD}. 

\begin{definition}\label{Natrual_Deduct_MDor}%[A natural deduction system of \MD]\label{Natrual_Deduct_PDbor}\ 
The natural deduction system of \MDor consists of all rules in \Cref{tab:MT}(b)-(e), together with the additional rules in \Cref{tab:mdor_rule}. %\todo{also no axiom}\todo{The assumption of \diamon cannot be empty, but the sequence in \boxmon can be empty.}
%is defined as follows:
%\begin{description}
%\item[Axiom]~ 
%\AxiomC{}\noLine\UnaryInfC{}\noLine\UnaryInfC{}\noLine\UnaryInfC{}\noLine\UnaryInfC{} \RightLabel{\exclmid}\UnaryInfC{$\alpha\sor\neg \alpha$}\noLine\UnaryInfC{}\DisplayProof
%
%\item[Rules] All rules in \Cref{tab:MT_nd} except for those that involve intuitionistic implication, together with the additional rules for negation and diamond listed in \Cref{tab:mdor_rule}.
%\end{description}
\end{definition}

\begin{table}[t]
\begin{center}
{\normalsize
%\begin{tabular}{|C{\linewidth}|}
%%\hline
%\multicolumn{1}{c}{Excluded Middle}\\\hline\hline
%\AxiomC{}\noLine\UnaryInfC{}\noLine\UnaryInfC{}\noLine\UnaryInfC{}\noLine\UnaryInfC{} \RightLabel{\exclmid}\UnaryInfC{$\alpha\sor\neg \alpha$}\noLine\UnaryInfC{}\DisplayProof
%\\\hline
% \end{tabular}
%
%
\begin{tabular}{|C{\linewidth}|}
%\hline
%\multicolumn{1}{c}{Rules for Negation}\\\hline
\hline
\AxiomC{}\noLine\UnaryInfC{$[\alpha]$}\noLine \UnaryInfC{$\vdots$}\noLine \UnaryInfC{$\bot$}\RightLabel{\negi}\UnaryInfC{$\neg\alpha$}\noLine\UnaryInfC{}\noLine\UnaryInfC{}\DisplayProof\quad\AxiomC{}\noLine\UnaryInfC{$\alpha$}\AxiomC{$\neg\alpha$}\RightLabel{\boti}\BinaryInfC{$\bot$}\DisplayProof\quad \AxiomC{}\noLine\UnaryInfC{$[\phi]$}\noLine\UnaryInfC{$D^\ast$}\noLine\UnaryInfC{$\vdots$}\noLine\UnaryInfC{$\psi$}\AxiomC{$\Diamond\phi$} \RightLabel{\diamon}\BinaryInfC{$\Diamond\psi$}\noLine\UnaryInfC{} \DisplayProof\quad\AxiomC{}\noLine\UnaryInfC{$\Diamond\neg\alpha$}\RightLabel{\diaboxinter}\UnaryInfC{$\neg\Box\alpha$}\noLine\UnaryInfC{}\DisplayProof\\%\hline
% \end{tabular}
%%\end{center}
%
%\vspace{0.5\baselineskip}
%
%%\begin{center}
%\begin{tabular}{|C{0.67\linewidth}C{0.3\linewidth}|}
%%\hline
%\multicolumn{2}{c}{Rules for Modalities}\\\hline\hline
%\AxiomC{$\Diamond\phi$} \AxiomC{$[\phi]$}\noLine\UnaryInfC{$D^\ast$}\noLine\UnaryInfC{$\vdots$}\noLine\UnaryInfC{$\psi$} \RightLabel{\diamon}\BinaryInfC{$\Diamond\psi$}\noLine\UnaryInfC{} \DisplayProof
%\\
($\ast$) The derivation $D$  is assumed to have no undischarged assumptions.\\
\hline
\end{tabular}
}
\caption{Rules for negation and modalities}
\label{tab:mdor_rule}
\end{center}
\end{table}

%\vspace{\baselineskip}
The system of \MDor has all the rules of \MT that do not involve implication, together with some additional rules for negation and modalities.
The clauses in \Cref{derivable_thm_MT_MID} are derivable easily also in the system of \MDor. In particular, item (e) ($\Diamond\alpha\dashv\vdash_{\LL}\neg\Box\neg\alpha$) can be derived without applying the \IK axiom $\neg\Diamond\bot$ as follows:
\begin{equation}\label{box2diamond_der}
\AxiomC{$\Diamond\alpha$}
\AxiomC{[$\alpha$]}\AxiomC{[$\neg\alpha$]}\RightLabel{\boti}\BinaryInfC{$\bot$}\RightLabel{\negi}\UnaryInfC{$\neg\neg\alpha$}
\RightLabel{\diamon}\BinaryInfC{$\Diamond\neg\neg\alpha$}\RightLabel{\diaboxinter}\UnaryInfC{$\neg\Box\neg\alpha$} \DisplayProof\quad\quad
\AxiomC{$\neg\Box\neg\alpha$}\RightLabel{\boxdiainter}\UnaryInfC{$\Diamond\neg\neg\alpha$}
\AxiomC{[$\neg\neg\alpha$]}\RightLabel{\dnege}\UnaryInfC{$\alpha$}
\RightLabel{\diamon}\BinaryInfC{$\Diamond\alpha$} \DisplayProof
\end{equation}
Note that we did not include the \IK axiom $\neg\Diamond\bot$ in our system of \MDor, because this classical formula is derivable in the system:
\begin{center}
\AxiomC{[$\bot$]}\RightLabel{\negi}\UnaryInfC{$\neg\bot$}\RightLabel{\boxmon}\UnaryInfC{$\Box\neg\bot$} \RightLabel{\Cref{derivable_thm_MT_MID}(d)}\UnaryInfC{$\neg\neg\Box\neg\bot$}\RightLabel{\Cref{derivable_thm_MT_MID}(e)}\UnaryInfC{$\neg\Diamond\bot$}\DisplayProof
\end{center}

To prove the Completeness Theorem of the system, we adopt a very similar argument to that in the previous section. We first show that every \MDor-formula is provably equivalent to a formula in disjunctive normal form.

\begin{lemma}[Normal Form]\label{nf_MDor}
For any \MDor-formula $\phi$, $\phi\dashv\vdash\bigvee_{i\in I}\alpha_i$ for some set $\{\alpha_i\mid i\in I\}$ of classical formulas.
\end{lemma}
\begin{proof}
We prove the lemma  by induction on $\phi$. 
%Note that 
%\begin{equation}\label{distrsorborsor}
%(\theta\otimes \eta)\bor(\theta\otimes \delta)\vdash\theta\otimes(\eta\bor\delta)
%\end{equation} holds for all formulas $\theta,\eta,\delta$.
If $\phi$ is a classical formula, then the lemma holds trivially. If $\phi=\dep(\alpha_1,\dots,\alpha_n,\beta)$, then we derive
\[
\begin{array}{rl}
\dep(\alpha_1,\dots,\alpha_n,\beta)\!\!&\displaystyle\dashv\vdash\bigotimes_{v\in {2}^{\{1,\dots,n\}}}(\alpha_1^{v(1)}\wedge\dots\wedge \alpha_n^{v(n)}\wedge (\beta\vee\neg \beta))\quad\text{(by \depdf)}\\
&\displaystyle\dashv\vdash\bigvee_{f\in 2^{2^{\{1,\dots,n\}}}}\bigotimes_{v\in {2}^{\{1,\dots,n\}}}(\alpha_1^{v(1)}\wedge\dots\wedge \alpha_n^{v(n)}\wedge \beta^{f(v)})\\
&\hfill\text{(by \Cref{derivable_thm_MT_MID}(a))}.
\end{array}
\]

The induction steps are proved by applying \Cref{derivable_thm_MT_MID}(a)-(c) and the induction hypothesis (cf. the proof of \Cref{dn_form_pid}).
%The cases $\phi=\psi\wedge\chi$ and $\phi=\psi\vee\chi$ follow easily from the induction hypothesis. The case $\phi=\psi\sor\chi$ follows from the induction hypothesis and  \Cref{derivable_thm_MT_MID}(a). The cases $\phi=\Box\psi$ and $\phi=\Diamond\psi$ follow from the induction hypothesis and \Cref{derivable_thm_MT_MID}(b)(c).
\end{proof}

\begin{lemma}\label{class_K_equiv_MDor}
If $\alpha$ and $\beta$ are classical formulas, then $\alpha\vdash_{\K}\beta\iff\alpha\vdash_{\MDor}\beta$.
\end{lemma}
\begin{proof}
%It is easy to verify that restricted to classical formulas, the axiom and rules of \MDor are sound in the sense of the usual single-world semantics of classical modal logic. 
The direction ``$\Longleftarrow$'' then follows from the Soundness Theorem and Expression (\ref{Team2K}) from \Cref{sec:prel}.  

For the direction ``$\Longrightarrow$'', it suffices to show that restricted to classical formulas, the system of \MDor admits all axioms and rules of the Hilbert-style system of \K. 

For the rules of \K, by inspecting the rules of the system of \MDor, we see that restricted to classical formulas, the system admits Uniform Substitution rule, and Necessitation rule is a special case of the rule \boxmon when there is no undischarged assumption in the rule. The Modus Ponens rule is interpreted as $\alpha,\neg\alpha\sor\beta/\beta$ in the language of \MDor and it can be derived as follows:
\begin{center}
\AxiomC{$\alpha$}\AxiomC{$[\neg\alpha]$}\RightLabel{\boti} \BinaryInfC{$\bot$}\RightLabel{\textsf{Ex falso}}\UnaryInfC{$\beta$}\AxiomC{$\neg\alpha\sor\beta$}\RightLabel{\sors}\BinaryInfC{$\beta\sor\beta$}\RightLabel{\sore}\UnaryInfC{$\beta$}\DisplayProof
\end{center}

For the propositional axioms of \K, it is easy to see that restricted to classical formulas, the system of \MDor contains all the rules for the classical propositional connectives conjunction, disjunction with respect to $\otimes$, negation and falsum $\bot$. Therefore all axioms of classical propositional logic are derivable in the system of \MDor.

For the axioms of \K that involve modalities, the validity of (an equivalent form of) the \K axiom is stated in the language of \MDor as $\Box(\alpha\wedge\beta)\dashv\vdash\Box\alpha\wedge\Box\beta$, which can be derived easily by applying \boxmon. Finally, we derived the inter-definability of $\Box$ and $\Diamond$, i.e., $\Diamond\alpha\dashv\vdash\neg\Box\neg\alpha$,  already  in (\ref{box2diamond_der}).
\end{proof}

\begin{theorem}[Completeness]\label{completeness_MDor}
For any \MDor-formulas $\phi$ and $\psi$,  $\phi\models\psi\iff\phi\vdash_{\MDor}\psi$.
\end{theorem}
\begin{proof}
By a similar argument to that of the proof of \Cref{completeness_MT}, where we apply \Cref{nf_MDor,class_K_equiv_MDor} instead.
%%We leave it for the readers to check that all the axiom and rules of \MDor are sound. 
%It suffices to show the direction ``$\Longrightarrow$''. Suppose $\models\phi$. 
%%Since \MDor is compact (\Cref{compactness}), we may assume that $\Gamma$ is finite and let $\psi=\bigwedge\Gamma$. 
%By \Cref{nf_MDor}, we have $\phi\dashv\vdash\alpha_1\vee\dots\vee\alpha_n$ for some classical formulas $\alpha_1,\dots,\alpha_n$. By the direction ``$\Longleftarrow$", we obtain $\models\alpha_1\vee\dots\vee\alpha_n$.
%
%
%which implies that for each $1\leq i\leq n$, there exists $1\leq j_i\leq k$ such that $\alpha_i\models\beta_{j_i}$, thereby $\models\alpha_i\to\beta_{j_i}$. Since $\alpha_i\to\beta_{j_i}$ is a classical formula, by Expression (\ref{Team2K}), we obtain $\alpha_i\vdash_{\K}\beta_{j_i}$, which by \Cref{class_K_equiv_MDor} implies that $\alpha_i\vdash_{\MDor}\beta_{j_i}$. Finally, we derive $\alpha\vdash_{\MDor}\beta$ by \bori and \bore.
%%The direction $\Longrightarrow$ follows from the same argument as that of the proof of \Cref{completeness_MT}.
\end{proof}

Since  \MDor is compact (by \Cref{compactness}), we obtain also the Strong Completeness Theorem as a corollary.

\begin{corollary}[Strong Completeness]
For any set $\Gamma\cup\{\phi\}$ of \MDor-formulas, $\Gamma\models\phi\iff \Gamma\vdash_{\MDor}\phi$.
\end{corollary}
%\begin{proof}
%By \Cref{completeness_MDor} and the compactness of \MDor (\Cref{compactness}).
%\end{proof}

%\todo{Strong completeness as a corollary}

\subsection{Applications of the disjunctive normal form}\label{sec:nf}

We devote this section to three interesting applications of the disjunctive normal form (\ref{dnf_mt}) of modal dependence logics. %Throughout this section, let \LL stand for an arbitrary logic among modal dependence logics.

In the context of propositional logics of dependence, flat formulas admit a certain characterization theorem; see \cite{Ciardelli_PhD,IemhoffYang15} for the proof. We now generalize this  characterization result to the modal case by using the disjunctive normal form.

%We first prove a characterization theorem for flat and classical formulas. A similar characterization theorem is true also for %Items (c) and (d) in the proposition below are applicable only to those logics that have the relevant formulas in their languages.

\begin{theorem}\label{flat_char}
%Let $\phi$ be an \MT-formula. 
The following are equivalent.

\vspace{0.5\baselineskip}
%\setlength{\columnsep}{0.01\linewidth}
%\begin{multicols}{2}
\begin{enumerate}%[(a)]
\begin{minipage}{0.6\linewidth}
\item[(a)] $\phi$ is flat
\item[(b)] $\phi\equiv\alpha$ for some classical formula $\alpha$
\end{minipage}
\begin{minipage}{0.4\linewidth}
\item[(c)] $\neg\neg\phi\equiv\phi$
\item[(d)] $\models\phi\sor\neg\phi$
\end{minipage}
\end{enumerate}
%\end{multicols}
\end{theorem}
\begin{proof}
We only give the detailed proof for (a)$\Rightarrow$(b). Assume (a). We have $\phi\equiv\bigvee_{i\in I}\alpha_i$ for some set $\{\alpha_i\mid i\in I\}$ of classical formulas, and in particular $\models\phi\to\bigvee_{i\in I}\alpha_i$. Since $\phi$ is flat,  it follows from  Expression (\ref{KP_valid}) from \Cref{sec:mt} and the Disjunction Property (\Cref{disjunct_prop}) that  there exists $j\in I$ such that $\phi\models\alpha_{j}$. On the other hand, $\alpha_{j}\models\bigvee_{i\in I}\alpha_i$. Hence, $\phi\equiv\alpha_{j}$.
\end{proof}

We write $\phi(\vec{p})$ to indicate that the propositional variables occurring in $\phi$ are among $\vec{p}=p_1\dots p_n$. Next, we prove Craig's Interpolation Theorem for modal dependence logics that have intuitionistic disjunction in their languages.

\begin{theorem}[Interpolation]
Let \LL be a modal dependence logic that has intuitionistic disjunction in its language. For any \LL-formulas $\phi(\vec{p},\vec{q})$ and $\psi(\vec{q},\vec{r})$, if $\phi\vdash_\LL\psi$, then there exists an \LL-formula $\theta(\vec{q})$ such that $\phi\vdash_\LL\theta$ and $\theta\vdash_\LL\psi$.
\end{theorem}
\begin{proof}%\footnote{The author would like to thank Katsuhiko Sano for pointing out to me that my same argument for the proof of the Interpolation Theorem of propositional dependence logic applies also in the modal case.}
Suppose $\phi(\vec{p},\vec{q})\vdash_\LL\psi(\vec{q},\vec{r})$. Then $\bigvee_{i\in I}\alpha_i\dashv\vdash\phi\vdash\psi\dashv\vdash\bigvee_{j\in J}\beta_j$ for some sets $\{\alpha_i(\vec{p},\vec{q})\mid i\in I\}$ and $\{\beta_j(\vec{q},\vec{r})\mid j\in J\}$ of classical formulas. Then, for each $i\in I$, there exists $j_i\in J$ such that $\alpha_i\vdash_\LL\beta_{j_i}$. Since $\alpha_i$ and $\beta_{j_i}$ are classical formulas, $\alpha_i(\vec{p},\vec{q})\vdash_\K\beta_{j_i}(\vec{q},\vec{r})$. Now, by the Interpolation Theorem of \K, there exists a classical formula $\theta_i(\vec{q})$ such that $\alpha_i(\vec{p},\vec{q})\vdash_\K\theta_i(\vec{q})$ and $\theta_i(\vec{q})\vdash_\K\beta_{j_i}(\vec{q},\vec{r})$. Thus, $\alpha_i(\vec{p},\vec{q})\vdash_\LL\theta_i(\vec{q})$ and $\theta_i(\vec{q})\vdash_\LL\beta_{j_i}(\vec{q},\vec{r})$. The formula $\bigvee_{i\in I}\theta_i(\vec{q})$ is in the language of \LL, and clearly $\bigvee_{i\in I}\alpha_i(\vec{p},\vec{q})\vdash_\LL\bigvee_{i\in I}\theta_i(\vec{q})$ and $\bigvee_{i\in I}\theta_i(\vec{q})\vdash_\LL\bigvee_{j\in J}\beta_j$. 
\end{proof}

Lastly, we prove that modal dependence logics have the finite model property. %The argument below applies also to the logics \MD and \MDe that do not have intuitionistic disjunction $\vee$ in their languages.

\begin{theorem}[Finite Model Property]\label{fmp}
If $\not\models\phi$, then there exists a finite Kripke model $\mathfrak{M}$ and finite team $X$ such that $\mathfrak{M},X\not\models\phi$.
\end{theorem}
\begin{proof}
For any formula $\phi$, we have $\phi\equiv\bigvee_{i\in I}\alpha_i$ for some finite set $\{\alpha_i\mid i\in I\}$ of classical formulas. If $\not\models\phi$, then $\not\models\alpha_i$ for all $i\in I$. Since each $\alpha_i$ is a classical formula, by Expression (\ref{Team2K}) from \Cref{sec:prel}, $\not\vdash_{\K}\alpha_i$ for each $i\in I$. By the finite model property of \K, for each $i\in I$, there exists a finite Kripke model $\mathfrak{M}_i$ and $w_i$ such that $\mathfrak{M}_i,w_i\not\models\alpha_i$. It follows that $\mathfrak{M}_i,\{w_i\}\not\models\alpha_i$ in the sense of team semantics. Consider the finite model $\mathfrak{M}=\biguplus_{i\in I}\mathfrak{M}_i$ and the finite team $X=\{w_i\mid i\in I\}$. By Expression (\ref{truth-preserve-operators}) from \Cref{sec:prel}, we obtain that for each $i\in I$, $\mathfrak{M},\{w_i\}\not\models\alpha_i$, which implies $\mathfrak{M},X\not\models\alpha_i$ by the downward closure property. Hence, we conclude that $\mathfrak{M},X\not\models\bigvee_{i\in I}\alpha_i$, thereby $\mathfrak{M},X\not\models\phi$.
\end{proof}

%\todo{strong completeness via normal form}

\subsection{\MD and \MDe}\label{sec:md}

In this section, we define natural deduction systems for \MD and \MDe. As we pointed out in the introduction of \Cref{sec:axioma}, in these implication-free logics (which are compact by \Cref{compactness}) there is a subtle difference between the weak and the strong completeness. We  first introduce the systems for the two logics for which the strong completeness holds, and then point out that the systems with two rules less already admit the weak completeness.
%, which do not have intuitionistic disjunction $\vee$ in their languages. 
%These systems are generalizations of the system of propositional dependence logic defined in \cite{VY_PD}. \todo{Cite and comment on Mikka's paper}

The systems of \MD and \MDe have (essentially) the same rules for (extended) dependence atoms as introduced in \cite{VY_PD}. To define these rules, let us follow \cite{VY_PD} and first introduce some notations. A formula in the language of \MDe or \MD is a finite string of symbols. 
%A propositional variable $p$ is a symbol and the other symbols are $\wedge,\sor,\neg,=,(,),\Box,\Diamond$. 
We number the symbols in a formula with positive integers starting from the left, as in the following example:
\begin{center}
\begin{tabular}{cccccccccccccccc}
$=$& $($&$\Box$&$p$& $,$&$q$&$)$&$\sor$&$\Box$&$\dep$&$($&$\Box$&$p$&,&$q$&$)$\vspace{4pt}\\
1&2&3&4&5&6&7&8&9&10&11&12&13&14&15&16
\end{tabular}
\end{center}
%\begin{center}
%\begin{tikzpicture}
%\node at (0.15,0) {{\scriptsize$1~~~2~~~3~~~~~4~~~~5~~~6~~~7~~~8~~~9~~10~~11~~12~~13~~14$}};
%
%\node (0f)  at (0, 0.6) {$(~\neg~ p_1~\sor~  \dep~(~p_2~)~)~\wedge~\dep~(~p_2~)$};
%
%\end{tikzpicture}
%\end{center}
If the $m$th symbol of a formula $\phi$ starts a string $\psi$ that is a subformula of $\phi$, we denote the subformula by $[\psi,m]_\phi$, or simply $[\psi,m]$. When referring to an occurrence of a formula $\chi$ inside a subformula $\psi$ of $\phi$, we will be sloppy about the notations and use the same counting also for the subformula $\psi$.  We write $\phi(\beta/[\alpha,m])$ for the formula obtained from $\phi$ by replacing the occurrence of the subformula $[\alpha,m]$ with $\beta$. For example, for the  formula $\phi=\dep(\Box p,q)\sor\Box \dep(\Box p,q)$, we denote the second occurrence of the dependence atom $\dep(\Box p,q)$ by $[\dep(\Box p,q),10]$, and the same notation also designates the occurrence of $\dep(\Box p,q)$ inside the subformula $\Box \dep(\Box p,q)$. The notation $\phi(\beta/[\dep(\Box p,q),10])$ designates the formula $\dep(\Box p,q)\sor\Box \beta$.

\begin{table}[t]
\begin{center}
{\normalsize\begin{tabular}{|C{\linewidth}|}
\multicolumn{1}{c}{(a)}\\
%\multicolumn{2}{c}{Rules for Dependence Atoms}\\\hline\hline
\hline
\AxiomC{}\noLine\UnaryInfC{$\alpha$}\RightLabel{\depiz}\UnaryInfC{$\dep(\alpha)$}\DisplayProof
\quad\AxiomC{}\noLine\UnaryInfC{$\neg \alpha$}\RightLabel{\depiz}\UnaryInfC{$\dep(\alpha)$}\DisplayProof%\\
%\\
%\multicolumn{2}{|c|}{
%\def\defaultHypSeparation{\hskip .06in}
%\AxiomC{}\noLine\UnaryInfC{}\noLine\UnaryInfC{[$\phi(\alpha/[\dep(\alpha),m])$]}\noLine\UnaryInfC{$\vdots$} \noLine\UnaryInfC{$\theta$} \AxiomC{[$\phi(\neg \alpha/[\dep(\alpha),m])$]}\noLine\UnaryInfC{$\vdots$} \noLine\UnaryInfC{$\theta$}\AxiomC{$\phi$} \RightLabel{\depez}\TrinaryInfC{$\theta$}\noLine\UnaryInfC{}\DisplayProof
%%}
\\
~\\
\def\defaultHypSeparation{\hskip .09in}
 \AxiomC{}\noLine\UnaryInfC{$[\dep(\alpha_1)]$}
 \AxiomC{$\dots$}
 \AxiomC{}\noLine\UnaryInfC{$[\dep(\alpha_k)]$}
\noLine
\TrinaryInfC{}
\branchDeduce
\DeduceC{$\dep(\beta)$}
\RightLabel{\depik}\UnaryInfC{$\quad\dep(\alpha_1,\dots, \alpha_k,\beta)\quad$}\noLine\UnaryInfC{}\DisplayProof\quad
\AxiomC{$\dep(\alpha)$} \AxiomC{}\noLine\UnaryInfC{$[\alpha]$}\noLine\UnaryInfC{$\vdots$}\noLine\UnaryInfC{$\theta$} \AxiomC{$[\neg \alpha]$}\noLine\UnaryInfC{$\vdots$}\noLine\UnaryInfC{$\theta$} \RightLabel{\depez}\TrinaryInfC{$\theta$} \DisplayProof
\\\hline
\multicolumn{1}{c}{}\\
\end{tabular}
%\caption{Rules for Dependence Atoms (I)}
%\label{tab:dep_rule}
}
%\end{center}
%\end{table}%
%
%\begin{table}[t]
%\begin{center}
{\normalsize\begin{tabular}{|C{\linewidth}|}
%\multicolumn{2}{c}{Rules for Dependence Atoms}\\\hline\hline
\multicolumn{1}{c}{(b)}\\
\hline
%\AxiomC{}\noLine\UnaryInfC{$\alpha$}\RightLabel{\depiz}\UnaryInfC{$\dep(\alpha)$}\DisplayProof\AxiomC{}\noLine\UnaryInfC{$\neg \alpha$}\RightLabel{\depiz}\UnaryInfC{$\dep(\alpha)$}\DisplayProof\\
%\AxiomC{$\dep(\alpha)$} \AxiomC{}\noLine\UnaryInfC{$[\alpha]$}\noLine\UnaryInfC{$\vdots$}\noLine\UnaryInfC{$\theta$} \AxiomC{$[\neg \alpha]$}\noLine\UnaryInfC{$\vdots$}\noLine\UnaryInfC{$\theta$} \RightLabel{\depez}\TrinaryInfC{$\theta$} \DisplayProof
%\\
%\multicolumn{2}{|c|}{
\!\!\!\!\!\!\!\!
\def\defaultHypSeparation{\hskip .03in}
\AxiomC{}\noLine\UnaryInfC{}\noLine\UnaryInfC{}\noLine\UnaryInfC{$\dep(\alpha_1,\dots,\alpha_k,\beta)$} \AxiomC{$\dep(\alpha_1)$} \AxiomC{\!\!$\dots$\!\!}\AxiomC{$\dep(\alpha_k)$} \RightLabel{\depek}\QuaternaryInfC{$\dep(\beta)$}\DisplayProof\\
\def\defaultHypSeparation{\hskip .06in}
\AxiomC{}\noLine\UnaryInfC{}\noLine\UnaryInfC{[$\phi(\alpha/[\dep(\alpha),m])$]}\noLine\UnaryInfC{$\vdots$} \noLine\UnaryInfC{$\theta$} \AxiomC{[$\phi(\neg \alpha/[\dep(\alpha),m])$]}\noLine\UnaryInfC{$\vdots$} \noLine\UnaryInfC{$\theta$}\AxiomC{$\phi$} \RightLabel{\se}\TrinaryInfC{$\theta$}\noLine\UnaryInfC{}\DisplayProof
%}
\\\hline
%~\\
%\def\defaultHypSeparation{\hskip .09in}
%\!\!\!\!\!\!\!\!\!\!\!\!\!\!
% \AxiomC{}\noLine\UnaryInfC{$[\dep(\alpha_1)]$}
% \AxiomC{$\dots$}
% \AxiomC{$[\dep(\alpha_k)]$}
%\noLine
%\TrinaryInfC{}
%\branchDeduce
%\DeduceC{$\dep(\beta)$}
%\RightLabel{\depik}\UnaryInfC{$\quad\dep(\alpha_1,\dots, \alpha_k,\beta)\quad$}\noLine\UnaryInfC{}\DisplayProof
%\!\!\!\!\!\!\!\!
%\def\defaultHypSeparation{\hskip .03in}
%\AxiomC{}\noLine\UnaryInfC{}\noLine\UnaryInfC{}\noLine\UnaryInfC{$\dep(\alpha_1,\dots,\alpha_k,\beta)$} \AxiomC{$\dep(\alpha_1)$} \AxiomC{\!\!$\dots$\!\!}\AxiomC{$\dep(\alpha_k)$} \RightLabel{\depek}\QuaternaryInfC{$\dep(\beta)$}\DisplayProof
%\\\hline
\end{tabular}
%\caption{Rules for Dependence Atoms (II)}
%\label{tab:dep_se}
}
\caption{Rules for Dependence Atoms}
\label{tab:dep_rule}
\end{center}
\end{table}%

\begin{definition}\label{Natrual_Deduct_MD}
\begin{itemize}
\item The natural deduction system of \MDe consists of the rules in \Cref{tab:MT}(b)(d),  together with the rules  in \Cref{tab:dep_rule}.
\item The natural deduction system of \MD is the same as that of \MDe except that the dependence atoms can only have propositional variables as arguments.
\end{itemize}
\end{definition}

%\vspace{-1.5\baselineskip}

%\vspace{\baselineskip}

In the above systems, the rules \depiz, \depez and \se for dependence atoms simulate the equivalence $\dep(\alpha)\equiv \alpha\vee\neg\alpha$, and the rules \depik and \depek simulate the equivalence $\dep(\alpha_1,\dots,\alpha_k,\beta)\equiv \dep(\alpha_1)\wedge\dots\wedge\dep(\alpha_k)\to\dep(\beta)$ (see also Expression (\ref{dep_def}) in \Cref{sec:prel}). Clearly,  \depez is a special case of \se, but we present both rules in \Cref{tab:dep_rule} for reasons that will become clear in the sequel. We refer the reader to \cite{VY_PD} for further discussion on these rules. 

For simplicity, we  only give the proof of the Completeness Theorem for the system \MDe, from which  the Completeness Theorem for the system \MD follows. We follow the argument in \cite{VY_PD} for propositional dependence logic, and first define \emph{realizations} of formulas, a crucial notion of the argument.
Let $d=\dep(\alpha_1,\dots,\alpha_k,\beta)$ be a dependence atom. A function $f:2^{\{1,\dots,k\}}\to 2$ is called a \emph{realizing function} for $d$, where we stipulate $2^\emptyset=\{\emptyset\}$,  and the formula
\[d_{f}^\ast:=\displaystyle\bigotimes_{v\in 2^{\{1,\dots,k\}}}\left(\alpha_1^{v(1)}\wedge \dots \wedge \alpha_k^{v(k)}\wedge \beta^{f(v)}\right)\]
is called a \emph{realization} of the dependence atom $d$ over $f$. Let 
\(o=\langle [d_1,m_1],\,\dots,\,[d_c,m_c]\rangle\)
be the sequence of all occurrences of dependence atoms in  $\phi$.
%where 
%\[\{j^1_0,\dots,j^1_{k_1},\dots,{j^c_0},\dots,j^c_{k_c}\}\subseteq \{i_1,\dots,i_n\},\]
%\[j^1_{k_1}=i_{k^1},\dots,j^c_{k_c}=i_{k^c}.\]
A \emph{realizing sequence of $\phi$} is a sequence $\Omega=\langle f_1,\dots,f_{c}\rangle$ such that each $f_i$ is a realizing function for $d_i$. We call the classical formula $\phi_\Omega^\ast$ defined %inductively 
as follows a \emph{realization of $\phi$}: 
\[\phi_{\langle f_1,\dots,f_{c}\rangle}^\ast:=\phi((d_1)_{f_1}^\ast/[d_1,m_1],\dots,(d_c)_{f_c}^\ast/[d_c,m_c]).\]
For example, consider the formula $\phi=\dep(\Box p,q)\sor\Box\dep(\Box p,q)$ that we discussed earlier. Consider two realizing functions $f,g:2^{\{1\}}\to 2$ for $\dep(\Box p,q)$, defined as
\[f(\mathbf{1})=1=f(\mathbf{0}),~g(\mathbf{1})=0\text{ and }g(\mathbf{0})=1,\]
where $\mathbf{1}(1)=1$ and $\mathbf{0}(1)=0$. Both $\langle f, g\rangle$ and $\langle g,f\rangle$ are realizing sequences of $\phi$ giving rise to two realizations 
\[\big(\dep(\Box p,q)\big)^\ast_f\sor\Box\big(\dep(\Box p,q))\big)^\ast_g=\big((\Box p\wedge q)\sor(\neg\Box p\wedge q)\big)\sor\Box\big((\Box p\wedge \neg q)\sor(\neg\Box p\wedge q)\big)\] and $\big(\dep(\Box p,q)\big)^\ast_g\sor\Box\big(\dep(\Box p,q))\big)^\ast_f$ 
of $\phi$. 
%Note that a realization $\phi^\ast_\Omega$ is always a classical formula, and the realizing sequence of a formula that does not contain dependence atoms (i.e., a classical formula) is the (unique) empty sequence $\langle\rangle$.

%We can prove in the system that $\phi$ and $\bigvee_{\Omega\in \Lambda}\phi^\ast_\Lambda$ are equivalent in effect, where $\Lambda$ is the set of all realizing sequences of $\phi$. For the proof of the Completeness Theorem, we only need to derive that each $\phi^\ast_\Omega$ implies $\phi$.

The next lemma states the crucial properties of realizations that will be applied in the proof of the Completeness Theorem.

\begin{lemma}\label{realization_lm}
Let $\phi$ be a formula, and $\Lambda$ the set of all realizing sequences of $\phi$. 
\begin{enumerate}
\item[(a)] $\phi^\ast_\Omega\vdash_{\MDe}\phi$ for any $\Omega\in\Lambda$.
\item[(b)]  If $\phi^\ast_\Omega\vdash_{\MDe}\psi$ for all  $\Omega\in\Lambda$, then $\phi\vdash_{\MDe}\psi$.
\item[(c)] $\phi\equiv\bigvee_{\Omega\in \Lambda}\phi^\ast_\Omega$.
\end{enumerate}
\end{lemma}

To prove item (a) of the above lemma, we first prove a technical lemma.
\begin{lemma}\label{replacement_lm}
Let $\psi$ be a subformula of $\phi$ that is not inside the scope of a dependence atom or a negation. If $\delta\vdash_{\MDe}\theta$, then $\phi(\delta/[\psi,m])\vdash_{\MDe}\phi(\theta/[\psi,m])$. %\todofy{monotonicity}
\end{lemma}
\begin{proof}
We prove the lemma by induction on the subformulas $\chi$ of $\phi$. 

The case when $\chi$ is an atom is trivial. If $\chi=\chi_0\sor\chi_1$ and without loss of generality we assume that the formula $[\psi,m]$ occurs in the subformula $\chi_0$. By the induction hypothesis, $\chi_0(\delta/[\psi,m])\vdash\chi_0(\theta/[\psi,m])$, which by \sors implies $\chi_0(\delta/[\psi,m])\sor\chi_1\vdash\chi_0(\theta/[\psi,m])\sor\chi_1$. The case $\chi=\chi_0\wedge\chi_1$ is proved analogously.

The case $\chi=\Diamond\chi_0$ follows from the induction hypothesis and \diamon, and the case $\chi=\Box\chi_0$ follows from the induction hypothesis and \boxmon.
\end{proof}

\begin{proof}[Proof of \Cref{realization_lm}]
(a) We prove the item by induction on the complexity of $\phi$. 

If $\phi$ does not contain any occurrences of dependence atoms, then the property holds trivially. If $\phi=\dep(\alpha_1,\dots,\alpha_k,\beta)$ is a dependence atom and $f:2^{\{1,\dots,k\}}\to 2$ a realizing function of $d=\phi$, by \depik, to show $d^\ast_f\vdash \dep(\alpha_1,\dots,\alpha_k,\beta)$ it suffices to derive $d^\ast_f,\dep(\alpha_1),\dots,\dep(\alpha_k)\vdash\dep(\beta)$. This is proved by a similar argument to that of the proof of Lemma 4.15 in \cite{VY_PD} which makes use of the rules \depiz, \depez and \depik.

If $\phi$ is a complex formula with $c$ occurrences of dependence atoms, and $\phi_{\langle f_1,\dots,f_{c}\rangle}^\ast:=\phi((d_1)_{f_1}^\ast/[d_1,m_1],\dots,(d_c)_{f_c}^\ast/[d_c,m_c])$, where $\Omega=\langle f_1,\dots,f_{c}\rangle$. Then, since $(d_i)^\ast_{f_i}\vdash d_i$ for each $1\leq i\leq c$, we derive 
\[\phi((d_1)_{f_1}^\ast/[d_1,m_1],\dots,(d_c)_{f_c}^\ast/[d_c,m_c])\vdash\phi(d_1/[d_1,m_1],\dots,(d_c/[d_c,m_c])\] 
by applying \Cref{replacement_lm} repeatedly.

(b) This item is a special case of the statement of Lemma 4.18 in \cite{VY_PD}, and can be proved by essentially the same argument that makes use of \depiz, \depek, \se and other rules of the system of \MDe.

(c) The direction $\bigvee_{\Omega\in \Lambda}\phi^\ast_\Omega\models\phi$ follows from item (a) and the Soundness Theorem. We now prove the other direction $\phi\models\bigvee_{\Omega\in \Lambda}\phi^\ast_\Omega$ by induction on $\phi$.

The case when $\phi$ is a dependence atom can be easily checked using Expression (\ref{dep_def}) from \Cref{sec:prel}. The other propositional cases can be proved easily by the same argument as in Lemma 4.16 in \cite{VY_PD}. The case when $\phi=\Box\psi$ or $\phi=\Diamond\psi$ follows from the fact that $\Box(A\vee B)\models \Box A\vee \Box B$ and $\Diamond(A\vee B)\models \Diamond A\vee\Diamond B$.
\end{proof}

\begin{theorem}[Completeness]\label{mde_comp_thm}
For any \MDe-formula $\phi$ and $\psi$, $\phi\models\psi\iff\phi\vdash_{\MDe}\psi$.
\end{theorem}
\begin{proof}
Suppose $\phi\models\psi$. By \Cref{realization_lm}(c), we have 
\begin{equation*}%\label{comp_deduc_eq1}
\phi\equiv\bigvee_{\Omega\in \Lambda}\phi^\ast_\Omega\models\bigvee_{\Delta\in \Lambda'}\psi^\ast_\Delta\equiv
\psi
\end{equation*}
where $\Lambda$ and $\Lambda'$ are the (nonempty) sets of all realizing sequences of $\phi$ and $\psi$, respectively. Since each $\phi_\Omega^\ast$ and $\psi^\ast_\Delta$ are classical formulas, by (\ref{KP_valid}) from Section 2 we obtain that for each $\Omega\in\Lambda$, there is  $\Delta\in \Lambda'$ such that $\phi_\Omega^\ast\models\psi^\ast_{\Delta}$. From (\ref{Team2K}) from Section 1 we know that $\phi_\Omega^\ast\vdash_\K\psi^\ast_{\Delta}$, which implies $\phi_\Omega^\ast\vdash_{\MDe}\psi^\ast_{\Delta}$ by a similar argument to those in the previous sections (Cf. \Cref{class_K_equiv_MDor}). Now, by \Cref{realization_lm}(a) we derive $\phi_\Omega^\ast\vdash_{\MDe}\psi$. Finally, by \Cref{realization_lm}(b) we conclude that $\phi\vdash_{\MDe}\psi$.
\end{proof}

Since \MDe is compact (by \Cref{compactness}), we obtain the strong completeness as a corollary.

\begin{corollary}[Strong Completeness]
For any set $\Gamma\cup\{\phi\}$ of \MDe-formulas, $\Gamma\models\phi\iff \Gamma\vdash_{\MDor}\phi$.
\end{corollary}
%\begin{proof}
%By \Cref{mde_comp_thm} and the compactness of \MDe (\Cref{compactness}).
%\end{proof}

Finally, it is interesting to note that the \emph{theoremhood} or \emph{validity problem} of the logic \MDe can actually be axiomatized by a slightly weaker system that contains less rules than the one defined in \Cref{Natrual_Deduct_MD} for the \emph{entailment problem}. In the system of \Cref{Natrual_Deduct_MD}, if we drop the rules in \Cref{{tab:dep_rule}}(b) and write $\vdash_{\MDe}^0\phi$ if $\phi$ is a theorem (i.e., a formula derivable from  the empty assumption) in the resulting system, items (a) and (b) of \Cref{realization_lm} are still true. By a very similar argument to the proof of \Cref{mde_comp_thm} (namely, simply discard the arguments that involve $\phi$), one can prove the following weaker form of Completeness Theorem  for this weaker system without applying \Cref{realization_lm}(b). 
%It is worth noting that since \MDe does not have intuitionistic implication in its language, the Strong Completeness Theorem does not follow immediately from compactness and the Weak Completeness Theorem. A similar subtle difference between the systems for entailment problem and validity problem was noted also in \cite{Hannula16entail} in a different system for modal dependence logic.

\begin{theorem}[(Weak) Completeness]
For any \MDe-formula $\phi$, $\models\phi\iff\vdash_{\MDe}^0\phi$.
\end{theorem}

%Noting that $\phi_{\langle f_1,\dots,f_{c}\rangle}^\ast:=\phi((d_1)_{f_1}^\ast/[d_1,m_1],\dots,(d_c)_{f_c}^\ast/[d_c,m_c])$, where $\Omega=\langle f_1,\dots,f_{c}\rangle$,  apply \Cref{dep_atom_realization_lm} and \Cref{replacement_lm} repeatedly we obtain the following desired corollary.
%
%\begin{corollary}\label{realization_lm}
%For each realizing sequence $\Omega$ of $\phi$, $\phi^\ast_\Omega\vdash\phi$.
%\end{corollary}

%\begin{theorem}[Completeness]
%For any \MDe-formula $\phi$, $\models\phi\iff\vdash_{\MDe}\phi$.
%\end{theorem}
%\begin{proof}
%We only show the direction ``$\Longrightarrow$".
%Suppose $\models\phi$. By \Cref{realization_lm}(b), we have $\phi\equiv \bigvee_{\Omega\in \Lambda}\phi^\ast_\Omega$, where $\Lambda$  is the (nonempty) set of all realizing sequences of $\phi$.  By a similar argument to that in the proof of \Cref{completeness_MT}, we obtain $\vdash_{\MDe}\phi^\ast_\Omega$ for some $\Omega\in \Lambda$. Hence, we conclude $\vdash_{\MDe}\phi$ by \Cref{realization_lm}(a). 
%\end{proof}

%\begin{lemma}\label{realization_str_lm} 
% If $\phi^\ast_\Omega\vdash_{\MDe}^\ast\psi$ for all realizing sequences $\Omega$ of $\phi$, then $\phi\vdash_{\MDe}\psi$.
%\end{lemma}
%\begin{proof}
%The lemma is a special case of the statement of Lemma 4.18 in \cite{VY_PD}, and can be proved by essentially the same argument that makes use of \depek, \se and other rules of the system of \MDe.
%\end{proof}

%Below we sketch the proof of the completeness of consequence relation with a slightly different argument.

\section{Interpreting team semantics in single-world semantics}

%\todo{intermediate modal logic}

In the previous section, we have defined the systems of modal dependence logics as extensions of Fischer Servi's intuitionistic modal logic \IK and inquisitive logic \Inql (which is a variant of the  Kreisel-Putnam  intermediate logic \KP).
%The dependence atom-free fragment of \MID, denoted by \MIDz, has the same syntax as \IK
In this section, we explore the connection between  the single-world-based  intuitionistic modal logic and intermediate  logics and modal dependence logics from the model-theoretic point of view. 
%Especially, we will interpret the modal team semantics in the single-world semantics setting. 
We  first prove that the team semantics of modal dependence logics over a usual (modal) Kripke model $\mathfrak{M}$ coincides with the usual single-world semantics over an intuitionistic Kripke model $\mathfrak{M}'$, whose domain consists of the teams of $\mathfrak{M}$. 
%This proof confirms the common intuition that team semantics is a natural generalization of the usual single-world  semantics. 
For simplicity, we  only perform this construction for the dependence atom-free fragment of \MID and \MT, denoted \MIDz and \MTz, subsequently in \Cref{sec:mid_sin} and \Cref{sec:tensor_dia}. Depending on whether tensor $\otimes$ is present in the language of the logic, the domain of a model $\mathfrak{M}'$ for interpreting the team semantics will consist of either the full powerset of the domain of $\mathfrak{M}$ or the same powerset excluding the empty set. The tensor $\otimes$  (which corresponds to \emph{multiplicative conjunction}, as discussed in  \cite{AbVan09}, but is often understood as a \emph{disjunction}) will be naturally interpreted as a binary diamond modality in this framework.

Furthermore, we  generalize the properties of the specific powerset models we built for interpreting team semantics to establish the connection on a general level. In \Cref{sec:mid_sin}   we identify a class of bi-relation intuitionistic Kripke models that enjoy the abstract properties of the powerset models for \MIDz, and we  show that the system of \MIDz  defined in the previous section\footnote{The deduction system of \MIDz is obtained (natrually) from the system of \MID, as defined in \Cref{Hilbert_MT} or in \Cref{MT_nds_df}, by simply dropping all the rules that involve dependence atoms. Similarly for the deduction system of \MTz.} is complete with respect to this class of models in the single-world semantics sense. Similar result for \MTz will be obtained in \Cref{sec:tensor_dia} with respect to a class of tri-relation intuitionistic Kripke models with an extra ternary relation corresponding to the binary diamond $\otimes$. 
%These completeness results also link the team semantics and the single-world semantics in a concrete sense, which we hope will inspire further research. 
This approach is based on a similar construction for inquisitive logic given in \cite{ivano_msc}.

%A model $\mathfrak{M}^\circledcirc$ for interpreting the team semantics of \MID will have the set of all nonempty teams of $\mathfrak{M}$ as its domain, and a model $\mathfrak{M}^\circledcirc$ for interpreting the team semantics of \MT instead has the set of all teams of $\mathfrak{M}$ as its domain

%We will perform this construction first for \MIDz in \Cref{sec:mid_sin}, and then for \MTz \Cref{sec:tensor_dia}
%
%
%Depending on whether we include the tensor $\sor$ in the language of the logic we consider, the domain of the induced intuitionistic Kripke model will consist of either the nonempty teams 
%
%For simplicity, in this approach we will only consider 
%
%
%
%We have seen in the previous sections that team-based modal dependence logic inherit many nice properties from the single-world based usual modal logic. It is then natural to think that the team semantics of modal dependence logics actually has a single-world semantics interpretation. In this section, we make this connection precise. In particular, we will prove that (the dependence atom-free fragment of) modal intuitionistic dependence logic can be viewed as a type of (single-world based) intermediate modal logic over a certain class of bi-relation intuitionistic Kripke models, and the tensor (disjunction) $\sor$ of team semantics can be interpreted naturally as a binary diamond modality. Our approach is inspired by the similar interpretation for inquisitive logic given in \cite{ivano_msc}. \todo{Let us denote the dependence atom-free fragment of \MID by \MIDz}

\subsection{A single-world semantics for \MIDz}\label{sec:mid_sin}
In this section, we define a single-world semantics for the system of the dependence atom-free fragment of \MID (\MIDz). 
%and prove that \MIDz is complete (in the single-world semantics sense) to a class of bi-relation intuitionistic Kripke models.
Observe from \Cref{Hilbert_MT} that the set of theorems of \MIDz  includes all theorems of \IK and is included in the set of theorems of \K. In other words, \MIDz can be understood as an intermediate modal logic\footnote{In the literature intermediate modal logics are often obtained by adding to intuitionistic modal logic (either Fischer Servi's \IK or some other versions) extra modal axioms, such as $\mathbf{S4}$, $\mathbf{S5}$ axioms. The approach we take in this paper is, roughly, to add to \IK an extra propositional axiom, the \KP axiom.} that is not closed under uniform substitution (also called an intermediate modal theory). To be more precise, \MIDz can be viewed as the fusion of \IK and \KP together with one axiom stating that box distributes over disjunction (i.e., $\Box(\phi\vee\psi)\to(\Box\phi\vee\Box\psi)$), and two axioms describing the classical behavior of disjunction-free formulas (i.e., $\neg\Box\alpha\to\Diamond\neg\alpha$ and $\neg\neg\alpha\to\alpha$). % In addition, instances of axioms of \KP in the modal language of \MIDz are also theorems of \MIDz.
 Recall that \IK is complete with respect to bi-relation intuitionistic Kripke frames (see e.g. \cite{SimpsonPhDthesis}) and \KP is complete with respect to \KP-frames (see e.g. \cite{ZhaCha_ml}). We will show in this section that \MIDz is complete (in the single-world semantics sense) with respect to a class of Kripke models whose frames are both bi-relation intuitionistic Kripke frames and \KP-frames. 
 
 Let us first recall  relevant definitions for \IK.
%The system of \MIDz is also an extension of \KP.

%In this section, we show that \MIDz can be viewed as a type of (single-world based) intermediate modal logics (i.e., modal logics between intuitionistic modal logic and classical modal logic).

% interpret the team semantics of  the dependence atom-free fragment of modal intuitionistic dependence logic \MIDz

%the connection between the dependence atom-free fragment of modal intuitionistic dependence logic \MIDz and (the single-world based) intermediate modal logics (i.e., modal logics between intuitionistic modal logic and classical modal logic).

%We have defined the Hilbert system of \MID (and therefore the system of \MIDz, namely the system of \MID without the axioms for dependence atoms) as an extension of both Fischer Servi's intuitionistic modal logic \IK and Kreisel-Putnam logic \KP (see \Cref{Hilbert_MT}). Recall that \IK is complete with respect to bi-relation intuitionistic Kripke frames (see e.g. \cite{SimpsonPhDthesis}) and \KP is complete with respect to \KP-frames (see e.g. \cite{ZhaCha_ml}). We will show in this section that \MIDz is complete (in the single-world semantics sense) with respect to a class of Kripke models whose frames are both bi-relation intuitionistic Kripke frames and \KP-frames. Let us first recall  relevant definitions for \IK.

 \begin{figure}[t]
 \begin{center}
\begin{tikzpicture}[scale=.8, transform shape] 
\node  (1)  at (0, 0) {\normalsize$w$};
\node  (3) at +(0,2) {\normalsize$w'$};
\node (4) at +(2,2) {\normalsize$v'$};
\node  (2) at +(2,0) {\normalsize$v$};
\draw[line width=1pt, dashed] (2) -- (4);
\draw[line width=1pt] (1) -- (3);
\draw [dashed, ->] (3) -- (4);
\draw [->] (1) -- (2);
\node at (1,-0.5) {\normalsize(F1)};
\node  (11)  at (5, 0) {\normalsize$w$};
\node  (13) at +(5,2) {\normalsize$w'$};
\node (14) at +(7,2) {\normalsize$v'$};
\node  (12) at +(7,0) {\normalsize$v$};
\draw[line width=1pt] (12) -- (14);
\draw[line width=1pt, dashed] (11) -- (13);
\draw [dashed, ->] (13) -- (14);
\draw [->] (11) -- (12);
\node at (6,-0.5) {\normalsize(F2)};

%\node at (3.5,-1.4) {(Arrows represent $R$ relation and undirected lines represent $\geq$ relation)};
\end{tikzpicture}
 \caption{Frame conditions. The directed lines represent the $R$ relation and the undirected lines represent the $\geq$ relation with the nodes positioned above being accessible from the ones positioned below.}\label{Fig_F1F2}
\end{center}
  \end{figure}
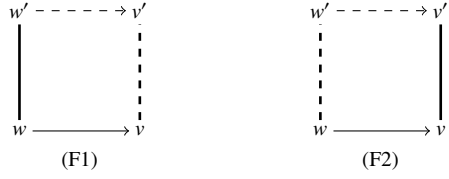

\begin{definition}\label{bi-relation_ml}
A \emph{bi-relation intuitionistic Kripke frame}  is a triple \allowbreak$\FF=(W,\geq,R)$, where
\begin{itemize}
%\addtolength{\itemsep}{-0.3\baselineskip}
\item $W$ is a nonempty set
\item $\geq$ is a partial ordering and $R$ is a binary relation on $W$
%\item $V:\rm{Prop}\to\wp(W)$ is a function (a \emph{valuation}) satisfying \emph{monotonicity} with respect to $\geq$, that is, $w\in V(p)$ and $w\geq v$ imply $v\in V(p)$;
\item $R$ and $\geq$ satisfy the following two conditions (F1) and (F2) (see  Figure \ref{Fig_F1F2}):
\begin{description}
  \item[F1] If $w\geq w'$ and $wRv$, then there exists $v'\in W$ such that $v\geq v'$ and $w'Rv'$.
  \item[F2] If $wR v$ and $v\geq v'$, then there exists $w'\in W$ such that $w\geq w'$ and $w'R v'$.
\end{description}
\end{itemize}
A \emph{bi-relation intuitionistic Kripke model} is a quadruple $\MM=(W,\geq,R,V)$ such that $(W,\geq,R)$ is a bi-relation intuitionistic Kripke frame and $V:\PROP\to\wp(W)$ is a \emph{valuation} satisfying \emph{monotonicity} with respect to $\geq$, that is, $w\in V(p)$ and $w\geq v$ imply $v\in V(p)$.
%\[[\,w\in V(p)\text{ and }w\geq v\,]~\Longrightarrow v\in V(p).\vspace{-2pt}\]
\end{definition}

%  Let us also recall the (single-world-based) Kripke semantics of \IK. %(\cite{FS81_IK},\cite{IK_FS84},\cite{IML_PlotkinStirling86}, see also \cite{SimpsonPhDthesis}).

\begin{definition}\label{sat_df_bir}
The  \emph{satisfaction} relation $\MM,w\Vdash\phi$ between a bi-relation intuitionistic Kripke model $\MM=\mathop{(W,\geq,R,V)}$, a node $w\in W$ and a formula $\phi$ in the language of \IK is defined inductively as follows:
%We define inductively the notion of a formula $\phi$ in the language of \IK being \emph{satisfied} in a bi-relation intuitionistic Kripke model $\MM=\mathop{(W,\geq,R,S,V)}$ on a node $w\in W$, denoted $\MM,w\Vdash\phi$, as follows:
\begin{itemize}
%\addtolength{\itemsep}{-0.3\baselineskip}
\item $\MM,w\Vdash p$ iff $w\in V(p)$
\item $\MM,w\nVdash\bot$
\item $\MM,w\Vdash\phi\wedge\psi$ iff $\MM,w\Vdash\phi$ and $\MM,w\Vdash\psi$
\item $\MM,w\Vdash\phi\vee\psi$ iff $\MM,w\Vdash\phi$ or $\MM,w\Vdash\psi$
\item $\MM,w\Vdash\phi\to\psi$ iff for all $v\in W$ such that $w\geq v$, if $\MM,v\Vdash\phi$, then $\MM,v\Vdash\psi$
%\[\MM,v\Vdash\phi\,\Longrightarrow\,\MM,v\Vdash\psi\]
\item $\MM,w\Vdash\Diamond\phi$ iff there exists $v\in W$ such that $wRv$ and $\MM,v\Vdash\phi$
\item $\MM,w\Vdash\Box\phi$ iff for all $u,v\in W$ such that $w\geq u$ and $uRv$, it holds that $\MM,v\Vdash\phi$
%\item $\MM,w\Vdash\otimes(\phi,\psi)$ iff there exist $v,u\in W$ such that $S(w,v,u)$, $\MM,v\Vdash\phi$ and $\MM,u\Vdash\psi$
\end{itemize}
%If $\MM,w\Vdash\phi$ for all nodes $w$ in a model \MM, then we say that $\phi$ is \emph{true} on the model \MM and write $\MM\models\phi$. If $(\FF,V)\models\phi$ for any model $(\FF,V)$ on a frame \FF, then we say that $\phi$ is \emph{valid} on the frame \FF and write $\FF\models\phi$.
\end{definition}

%In this section, we reveal the relationship between \MID and intuitionistic modal logic. Let us first recall some basic results on the intuitionistic modal logic \IK defined in the previous section.

It is easy to show that the $\geq$-monotonicity extends to arbitrary formulas $\phi$, that is, 
\begin{description}
\item[Monotonicity] [$\MM,w\Vdash \phi$ and $w\geq v$] $\Longrightarrow\,\MM,v\Vdash \phi$.
\end{description}

%In the remainder of this section, we will interpret the team semantics of \MIDz as single-world based semantics over a certain set $\mathsf{M}$ of bi-relation intuitionistic Kripke models and show that the system of \MIDz is complete with respect to $\mathsf{M}$.

%We will now interpret the team semantics of \MIDz over classical modal Kripke models as the single-world  semantics over certain bi-relation intuitionistic Kripke models. A crucial observation is that 
%every classical modal Kripke model induces a bi-relation intuitionistic Kripke model which we shall call \emph{powerset model}.

Every classical modal Kripke model induces a bi-relation intuitionistic Kripke model which we shall call \emph{powerset model}.

%Next, we show that the team semantics of \MIDz over classical modal Kripke models coincides with the single-world semantics over a bi-relation intuitionistic Kripke model 

\begin{definition}\label{pwset_ml_df}
Let $\MM=(W,R,V)$ be a classical modal Kripke model. The \emph{powerset  model $\MM^\circ$ induced by \MM} is a quadruple $\MM^\circ=(W^\circ,\supseteq,R^\circ,V^\circ)$, where
\begin{itemize}
\item $W^\circ=\wp(W)\setminus\{\emptyset\}$, i.e. $W^\circ$ consists of all nonempty teams $X\subseteq W$
\item $\supseteq$ is the superset relation 
\item $XR^\circ Y$ iff $XRY$ iff $Y\subseteq R(X)$ and $Y\cap R(w)\neq\emptyset$ for every $w\in X$%\footnote{Note that the $X$ and $Y$ on the left-hand side stand for single nodes in a model, while the same symbols on the right-hand side stand for teams (i.e., sets of nodes).}
%\item $R_\sor$ is a ternary relation defined as $R_\sor(X,Y,Z)$ iff $X=Y\cup Z$
\item $X\in V^\circ(p)$ iff $X\subseteq V(p)$
\end{itemize}
%All powerset Kripke models are defined in the above way, in other words, each powerset Kripke models $\MM^\circ$ is induced by a unique modal Kripke model \MM.
\end{definition}

%In a sense, the powerset model $\MM^\circ$ carries the information of teams of its associated modal Kripke model \MM, and the relation $R^\circ$ resembles the successor team relation $R$ on $W$.%, but $R^\circ\neq R$ as the successor team relation $R$ is a relation between teams (sets of worlds), whereas $R^\circ$ is a relation between single worlds of $W^\circ$. 

 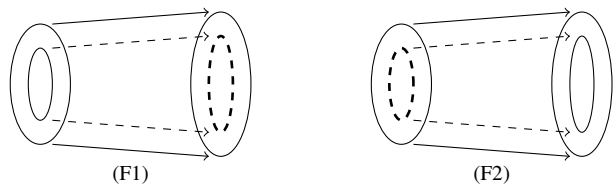
\begin{figure}[t]
 \begin{center}
\begin{tikzpicture}[scale=.8, transform shape] 
\draw (0, 1) ellipse (0.5 and 1);
\draw (0, 1) ellipse (0.2 and 0.6);
\draw (3, 1) ellipse (0.5 and 1.2);
\draw[line width=1pt, dashed]  (3, 1) ellipse (0.2 and 0.8);
%\node at (-0.8,1.2) {\normalsize$X$};
%\node at (3.8,1.2) {\normalsize$Y$};
%\node at (1.5,2.4) {\normalsize$R$};
\draw [->] (0.2,0) -- (2.8,-0.2);
\draw [dashed, ->] (0.2,0.4) -- (2.8,0.2);
\draw [->] (0.2,2) -- (2.8,2.2);
\draw [dashed, ->] (0.2,1.6) -- (2.8,1.8);
\node at (1.5,-0.5) {\normalsize(F1)};
\draw (6, 1) ellipse (0.5 and 1);
\draw[line width=1pt, dashed] (6, 1) ellipse (0.2 and 0.6);
\draw (9, 1) ellipse (0.5 and 1.2);
\draw  (9, 1) ellipse (0.2 and 0.8);
\draw [->] (6.2,0) -- (8.8,-0.2);
\draw [dashed, ->] (6.2,0.4) -- (8.8,0.2);
\draw [->] (6.2,2) -- (8.8,2.2);
\draw [dashed, ->] (6.2,1.6) -- (8.8,1.8);
\node at (7.5,-0.5) {\normalsize(F2)};
\end{tikzpicture}
\end{center}
 \caption{Conditions (F1) and (F2) in powerset models}\label{Fig_F1F2_MID}
\end{figure}

To see that $\MM^\circ$ is indeed a bi-relation intuitionistic Kripke model, note that the superset relation $\supseteq$ is a partial ordering, and the monotonicity of $V^\circ$ is immediate. To verify  condition (F1), for any $X,X',Y\in W^\circ$  such that $X\supseteq X'$ and $XR^\circ Y$, letting $Y'=R(X')\cap Y$, it is easy to show that $Y\supseteq Y'$ and $X'R^\circ Y'$ (see also Figure \ref{Fig_F1F2_MID}). Similarly, to verify condition (F2),  for any $X,Y,Y'\in W^\circ$ such that $XR^\circ Y$ and $Y\supseteq Y'$, letting $X'=R^{-1}(Y')\cap X$, clearly $X\supseteq X'$ and $X'R^\circ Y'$ (see also Figure \ref{Fig_F1F2_MID}).

Next, we show that the team-based satisfaction relation with respect to classical modal Kripke models is equivalent to the single-world-based satisfaction relation with respect to the associated powerset  models. 
\begin{lemma}\label{comodel_model}
Let $\MM=(W,R,V)$ be a classical modal Kripke model and $X\subseteq W$ a nonempty team. For any \MIDz-formula $\phi$, $\MM,X\models \phi\iff \MM^\circ,X\Vdash\phi$.
\end{lemma}
%\footnotetext{Note that the symbol ``$X$'' on the left-hand side stands for a team (a set of nodes), while the ``$X$'' on the right-hand side stands for a single node.}

\begin{proof}
We prove the lemma by induction on $\phi$. The only interesting case is when $\phi=\Box\psi$. If $\MM^\circ,X\Vdash \Box\psi$, then $\MM^\circ,R(X)\Vdash \psi$, since $X\supseteq X$ and $XR^\circ R(X)$. The induction hypothesis implies that $\MM,R(X)\models \psi$.  Hence $\MM,X\models \Box\psi$.

Conversely, if $\MM,X\models \Box\psi$, then $\MM,R(X)\models \psi$.  For all $Y,Z\in W^\circ$ such that $X\supseteq Y$ and $YR^\circ Z$, since $Z\subseteq R(X)$, the downward closure property implies that $\MM,Z\models \psi$ yielding $\MM^\circ,Z\Vdash \psi$ by the induction hypothesis. Hence $\MM^\circ,X\Vdash \Box\psi$.
% By the downward closure property of \MIDz, it follows that for all $Y,Z\in W^\circ$, such that $X\supseteq Y$ and $YR^\circ Z$, we have $\MM,Z\models \psi$, which implies $\MM^\circ,Z\Vdash \psi$ by the induction hypothesis. Hence $\MM^\circ,X\Vdash \Box\psi$.
%We have
%\begin{align*}
%\MM^\circ,X\Vdash \Box\psi&\Longrightarrow  \MM^\circ,R(X)\Vdash \psi\text{ (since $X\supseteq X$ and $XR^\circ R(X)$)}\\
%&\Longrightarrow \MM,R(X)\models \psi\text{ (by the induction hypothesis)}\\
%&\Longrightarrow \MM,X\models \Box\psi.
%\end{align*}
%and 
%\begin{align*}
%\MM,X\models \Box\psi&\Longrightarrow \MM,R(X)\models \psi\\
%&\Longrightarrow \text{for all nonempty }Y,Z\subseteq W\text{ s.t. }X\supseteq Y\text{ and }YRZ,~ \MM,Z\models \psi\\
%&~~~~~~~\text{ (since $Z\subseteq R(X)$ and $\models$ is downward closed)}\\\
%&\Longrightarrow \text{for all }Y,Z\in W^\circ\text{ s.t. }X\supseteq Y\text{ and }YR^\circ Z,~ \MM^\circ,Z\Vdash \psi\\
%&~~~~~~~\text{ (by the induction hypthesis)}\\\
%&\Longrightarrow \MM^\circ,X\Vdash \Box\psi.
%\end{align*}
\end{proof}

Inquisitive logic (being the propositional fragment of \MIDz) is shown in \cite{ivano_msc} to be complete with respect to negative saturated (single-relation) intuitionistic Kripke models, which are also negative \KP-models. Let us now give the corresponding definitions in the context of bi-relation intuitionistic Kripke models. %\todo{Mention or not mention the KP axioms}

%\begin{definition}
A point $w$ in a bi-relation intuitionistic Kripke model $\MM=\mathop{(W,\geq,R,V)}$ is called an \emph{$\geq$-endpoint} iff there is no point $v\neq w$ such that $w\geq v$. Denote by $E_w$ the set of all $\geq$-endpoints seen from $w$, i.e., 
\[E_w=\{v\in W\mid w\geq v\text{ and }v\text{ is an $\geq$-endpoint}\}.\]
A bi-relation intuitionistic Kripke frame $\FF=(W,\geq,R)$ is said to be  \emph{$\geq$-saturated} if for every $w\in W$, $E_w\neq \emptyset$, and for every nonempty subset $E\subseteq E_w$, there exists a point $v\in W$ such that $w\geq v$ and $E_v=E$.
%\begin{itemize}
%\item $E_w\neq \emptyset$;
%\item for every nonempty subset $E\subseteq E_w$, there exists $v\leq w$ such that $E_v=E$.
%\end{itemize}
A model $\MM$ is called  \emph{negative} if $\MM,w\Vdash p\iff \MM,w\Vdash \neg \neg p$. %A model $(\FF,V)$ is called a \emph{negative} model if $V$ is negative.
%\end{definition}
It is easy to verify that a powerset model $\MM^\circ=(W^\circ,\supseteq,R^\circ,V^\circ)$ is a negative $\geq$-saturated model, and in particular, $\supseteq$-endpoints in $\MM^\circ$ are singletons $\{w\}$ of elements $w$ in $W$.

%As in \cite{ivano_msc}, and as we shall see in the sequel, in a saturated bi-relation intuitionistic Kripke model, $\geq$-endpoints behave as singleton teams of a usual modal Kripke model.

The system of \MIDz (see \Cref{Hilbert_MT}) extends the systems of \IK and \Inql with two extra axioms: $\Box(\phi\vee\psi)\to(\Box\phi\vee\Box\psi)$ and $\neg\Box\alpha\to\Diamond\neg\alpha$, where $\alpha$ is any classical formula. The latter axiom is equivalent to the axiom $\neg\Box\neg p\to\Diamond\neg\neg p$, because a classical formula $\alpha$ is always equivalent to a (double) negation $\neg\neg\alpha$ (by  \Cref{flat_char}).
%By  \Cref{flat_char}, a classical formula $\alpha$ is equivalent to a (double) negation $\neg\neg\alpha$. Thus, we can rewrite the axiom $\neg\Box\alpha\to\Diamond\neg\alpha$ as  $\neg\Box\neg p\to\Diamond\neg\neg p$. 
%Also, the axiom $\Box(\phi\vee\psi)\to(\Box\phi\vee\Box\psi)$ can be written as $\Box(p\vee q)\to(\Box p\vee\Box q)$. 
In what follows we show that the axioms $\Box(p\vee q)\to(\Box p\vee\Box q)$ and $\neg\Box\neg p\to\Diamond\neg\neg p$ both characterize certain frame condition. We write $R_1\circ R_2$ for the \emph{composition} of the two binary relations $R_1$ and $R_2$ on a set $W$, defined as
\((x,y)\in R_1\circ R_2~\text{ iff }~\exists z\in W(xR_1z\wedge zR_2y).\)
%\todo{equivalence of axioms}

\begin{figure}[t]
\begin{center}
\begin{tikzpicture}[scale=.8, transform shape] 
\node  (1)  at (-7, -0.5) {\normalsize$w$};
%\node  (3) at +(-7,1) {$t'$};
\node  (3) at +(-7,1) {$\bullet$};
\node (4) at +(-7,2.5) {\normalsize$t$};
%\node  (2) at +(-8,1.2) {$u'$};
\node (2) at +(-8,1.2) {$\bullet$};
\node  (5) at +(-9,3) {\normalsize$u$};
%\node  (6) at +(-6,1.2) {$v'$};
\node  (6) at +(-6,1.2) {$\bullet$};
\node  (7) at +(-5,3) {\normalsize$v$};
\draw[line width=1pt] (1) -- (2);
\draw[line width=1pt] (1) -- (6);
\draw[line width=1pt, dashed] (4) -- (5);
\draw[line width=1pt, dashed] (4) -- (7);
\draw[line width=1pt, dashed] (1) -- (3);
\draw [dashed, ->] (3) -- (4);
\draw [->] (2) -- (5);
\draw [->] (6) -- (7);

\node at (-7,-1.5) {\normalsize(a)};

\draw (0, 0) ellipse (2.3 and 0.7);
\draw[fill=red, color=gray] (-1, 0) ellipse (0.55 and 0.3);
\draw[fill=blue, color=gray] (1, 0) ellipse (0.55 and 0.3);
\draw[dashed] (-1.7, -0.4) rectangle (1.7, 0.4);
\draw[fill=red, color=gray] (-1, 2.5) ellipse (0.8 and 0.3);
\draw[fill=blue, color=gray] (1, 2.5) ellipse (0.8 and 0.3);
\draw[dashed] (-1.9, 2) rectangle (1.9, 2.9);

\node at (2.2,2.3) {\normalsize$t$};

\node at (-1,3.1) {\normalsize$u$};
\node at (1,3.1) {\normalsize$v$};

\node at (2.3,-0.4) {\normalsize$w$};

\draw[line width=2pt, color=gray] [->] (-1,0) -- (-1,2.3);
%\draw[color=red] [->] (-0.45,0) -- (-0.2,1.9);
\draw[line width=2pt, color=gray] [->] (1,0) -- (1,2.3);
\draw[line width=2pt, dashed] [->] (0,0) -- (0,2.3);
%\draw[color=blue] [->] (1.55,0) -- (1.8,1.9);

\node at (0,-1.5) {\normalsize(b)};

\end{tikzpicture}
\caption{Condition (G1')}\label{fig:G1}
\end{center}
\end{figure}
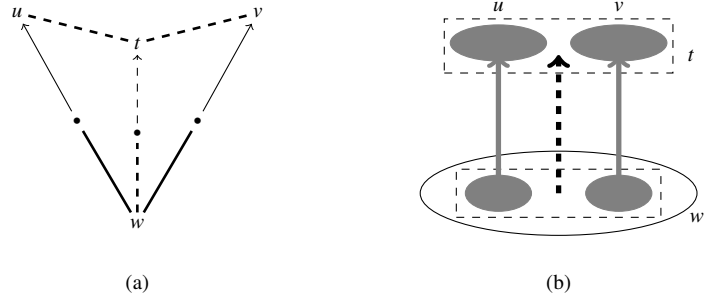

\begin{lemma}\label{sat_G1}
Let $\FF=(W,\geq,R)$ be a bi-relation intuitionistic Kripke frame. Then,
\(\FF\models\Box(p\vee q)\to(\Box p\vee\Box q)\iff \FF \text{ satisfies condition (G1') defined below}:\)
%where (G1') is defined as follows:
\begin{description}
\item[G1'] For all $w,u,v\in W$, if $u,v\in (\geq\circ R)(w)$%\footnote{We write $R_1\circ R_2$ for the \emph{composition} of the two binary relations $R_1$ and $R_2$ on a set $W$, defined as
%\((x,y)\in R_1\circ R_2~\text{ iff }~\exists z\in W(xR_1z\wedge zR_2y).\)}
, then there exists $t\in W$ such that $w(\geq\circ R)t$, $t\geq u$ and $t\geq v$. (See \Cref{fig:G1}(a)) %This condition is abstracted from the property of teams in a usual modal Kripke frame as illustrates in the right figure below.)
\end{description}
%In case \FF is finite, (G1') is equivalent to (G1) defined as follows:
%\begin{description}
%\item[(G1)]For any $w\in W$ and any nonempty $X\subseteq (\geq\circ R)(w)$, there exists a node $u\in (\geq\circ R)(w)$ such that $u\geq v$ for all $v\in X$.
%\end{description}
\end{lemma}

Before we give the proof of the lemma, let us first check that the  underlying frames of powerset models satisfy (G1'). First note that for any points $X,Y$ in a powerset model $\MM^\circ=(W^\circ,\supseteq,R^\circ,V^\circ)$, $X(\supseteq \circ R^\circ)Y$ iff $Y\subseteq R(X)$. Now, for any three points $w,u,v$ in $W^\circ$ such that $u,v\in\supseteq\circ R^\circ(w)$, we have $u,v\subseteq R(w)$ implying $u\cup v\subseteq R(w)$. Clearly, $t=u\cup v$ is a nonempty subset of $W$ such that $w(\supseteq\circ R^\circ) t$, $t\supseteq u$ and $t\supseteq v$ (see \Cref{fig:G1}(b)).
 %Indeed, for any three points $w,u,v$ in a powerset model $\MM^\circ=(W^\circ,\supseteq,R^\circ,V^\circ)$ such that $u,v\in\supseteq\circ R^\circ(w)$, clearly, $t=u\cup v$ is a nonempty subset of $W$ such that $w(\supseteq\circ R) t$, $t\supseteq u$ and $t\supseteq v$ (see \Cref{fig:G1}(b)). 
 As a powerset model $\MM^\circ$ carries the information of teams in the model $\MM$, condition (G1') can be viewed as a property that is abstracted from the corresponding property of teams of the usual classical modal Kripke frames.

In the sequel, we will also work with the following equivalent form (G1) of (G1'):
\begin{description}
\item[G1]For any $w\in W$ and any nonempty finite set $X\subseteq (\geq\circ R)(w)$, there exists a node $u\in (\geq\circ R)(w)$ such that $u\geq v$ for all $v\in X$.
\end{description}

\begin{proof}[Proof of \Cref{sat_G1}]
Suppose $\mathfrak{F}$ satisfies (G1') and $(\mathfrak{F},V),w\nVdash \Box p\vee\Box q$ for some valuation $V$ and some $w\in W$. Then there exist $u,v\in W$ such that $w(\geq\circ R)u$, $w(\geq\circ R)v$,
\[(\mathfrak{F},V),u\nVdash p\text{ and }(\mathfrak{F},V),v\nVdash  q.\]
Let $t\in W$ be the point given by (G1'). Then by the $\geq$-monotonicity, we have $(\mathfrak{F},V),t\nVdash p\vee q$,
which implies that $(\mathfrak{F},V),w\nVdash \Box(p\vee q)$.

Conversely, suppose that $\mathfrak{F}$ does not satisfy (G1'). Then there exist $w,u,v\in W$ such that $w(\geq\circ R)u$, $w(\geq\circ R)v$ and for all $t\in W$ such that $w\geq \circ Rt$, either $t\ngeq u$ or $t\ngeq v$. %In particular, $u\ngeq v$ and $v\ngeq u$. 
Clearly, we can find a $\geq$-monotone valuation $V$ such that
%Let $V$ be any valuation such that
\[ V(p)=W\setminus\geq^{-1} (v)\text{ and }V(q)=W\setminus\geq^{-1} (u).\]
%It is easy to see that $V$ is monotone.

For each $t\in W$ such that $w(\geq\circ R)t$, either $t\notin \geq^{-1} (v)$ or $t\notin \geq^{-1} (u)$. Thus $(\mathfrak{F},V),t\Vdash p\vee q$, thereby $(\mathfrak{F},V),w\Vdash \Box(p\vee q)$.
On the other hand,   $(\mathfrak{F},V),u\nVdash q$ and $(\mathfrak{F},V),v\nVdash p$. Hence $(\mathfrak{F},V),w\nVdash \Box p\vee \Box q$.
\end{proof}

%Note that condition (G1') is equivalent to (G1) defined as follows:

\begin{figure}[t]
\begin{center}
\begin{tikzpicture}[scale=.8, transform shape] 
\draw (0, 0)[fill=gray, color=gray] ellipse (2 and 0.6);
\draw (0, 2.5) ellipse (2.5 and 0.6);
\draw[color=gray!60, fill=gray!40, line width=1.5pt, dashed] (0.4, 2.5) ellipse (1.6 and 0.45);
%\tikzstyle{every node} = [circle, fill=black]
%\draw[color=white, fill=red!50] (0, 2.5) ellipse (2 and 0.5);

\node  (1)  at (-1, 0) {$\bullet$};
\node  (2)  at (0, 0) {$\bullet$};
\node  (3)  at (1, 0) {$\bullet$};
\node  (4)  at (-1, 2.5) {$\bullet$};
\node  (5)  at (-1.5, 2.5) {$\bullet$};
\node  (6)  at (-0.3, 2.5) {$\bullet$};
\node  (7)  at (0.4, 2.5) {$\bullet$};
\node  (8)  at (1.5, 2.5) {$\bullet$};
\draw [->] (1) -- (4);
\draw [->] (1) -- (5);
\draw [->] (1) -- (6);
\draw [->] (2) -- (6);
\draw [->] (2) -- (7);
\draw [->] (3) -- (8);

\node at (1,3.27) {\normalsize$\bigcup E=t$};
\node at (1,-0.7) {\normalsize$w$};

\node at (-7,-1.7) {\normalsize(a)};

\node at (0,-1.7) {\normalsize(b)};

\draw[line width=2pt, dashed, color=gray] [->] (0.65,0.2) -- (0.65,2.4);

\node (r) at (-7,-1) {\normalsize$w$};
\draw (-7, 0.2) ellipse (2 and 0.6);
\node  (l1)  at (-8, 0.2) {$\bullet$};
\node  (l2)  at (-7, 0.2) {$\bullet$};
\node  (l3)  at (-6, 0.2) {$\bullet$};
\draw [-] (r) -- (l1);
\draw [-] (r) -- (l2);
\draw [-] (r) -- (l3);

\node at (-9.5,0.2) {\normalsize$E_w$};

\draw (-7, 2.5) ellipse (2.5 and 0.6);

\node  (l4)  at (-8, 2.5) {$\bullet$};
\node  (l5)  at (-8.5, 2.5) {$\bullet$};
\node  (l6)  at (-7.3, 2.5) {$\bullet$};
\node  (l7)  at (-6.6, 2.5) {$\bullet$};
\node  (l8)  at (-5.5, 2.5) {$\bullet$};
\draw [->] (l1) -- (l4);
\draw [->] (l1) -- (l5);
\draw [->] (l1) -- (l6);
\draw [->] (l2) -- (l6);
\draw [->] (l2) -- (l7);
\draw [->] (l3) -- (l8);

\node at (-10.2,2.5) {\normalsize$R(E_w)$};
\node at (-6,3.25) {\normalsize$E$};
\draw[color=white, gray, line width=1.5pt] (-6.6, 2.5) ellipse (1.6 and 0.45);

\node(t) at (-4.6,0.9) {\normalsize$t$};
\draw [->,dashed] (r) -- (t);

\draw [-,dashed] (t) -- (l6);
\draw [-,dashed] (t) -- (l7);
\draw [-,dashed] (t) -- (l8);

\end{tikzpicture}
\caption{Condition (G2)}
\label{fig:G2}
\end{center}
\end{figure}
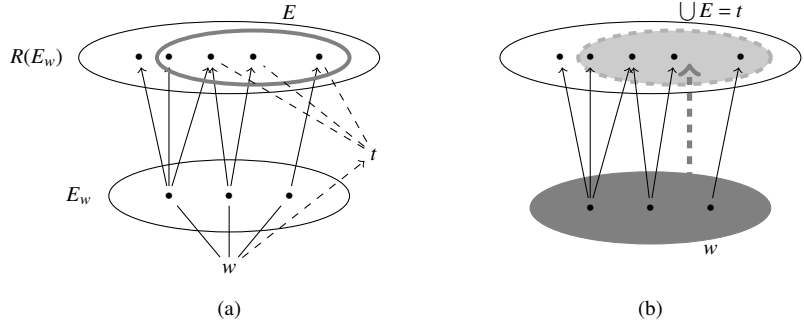

\begin{lemma}\label{sat_G2}
Let $\mathfrak{F}=(W,\geq,R)$ be a $\geq$-saturated bi-relation intuitionistic Kripke frame. Then,
\(\mathfrak{F}\models\neg\Box\neg p\to\Diamond\neg\neg p\iff \mathfrak{F} \text{ satisfies condition (G2) defined below}:\)
%where (G2) is defined as follows:
\begin{description}
\item[G2] Let $w\in W$ be an arbitrary point and $E$ a set of $\geq$-endpoints such that $E\subseteq R(E_w)$ and $E\cap R(v)\neq\emptyset$ for every $v\in E_w$.
%, there exists $u_v\in E$ with $v R u_v$. 
Then, there exists $t\in W$ such that
\(wRt\text{ and }E_t\subseteq E.\) (see \Cref{fig:G2}(a)) %This condition is abstracted from the property of teams in a usual modal Kripke frame as illustrates in the right figure below.)
\end{description}
\end{lemma}

Before we give the proof of the lemma, let us first check that the underlying frames of powerset models satisfy (G2). Indeed, for any point $w$ in a powerset model $\MM^\circ=(W^\circ,\supseteq,R^\circ,V^\circ)$ and any set $E$ of $\supseteq$-endpoints (i.e., a set of singletons of elements in $W$) such that $E\subseteq R^\circ(E_w)$ and $E\cap R^\circ (v)$ for every $v\in E_w$, it is easy to see that $t=\bigcup E$ is a nonempty subset of $W$ such that $wR^\circ t$ and $E_t= E$ (see \Cref{fig:G2}(b)).

\begin{proof}[Proof of \Cref{sat_G2}]
Suppose $\mathfrak{F}$ satisfies (G2) and $(\mathfrak{F},V),w\Vdash \neg\Box\neg p$ for some valuation $V$ and some $w\in W$. Then, $(\mathfrak{F},V),v\nVdash\Box\neg p$ for each $\geq$-endpoint $v\geq w$, i.e., each $v\in E_w$. It follows that there exists $u_v'$ such that $vRu_v'$ and $(\mathfrak{F},V),u_v'\nVdash \neg p$, which implies $(\mathfrak{F},V),u_v\Vdash p$ for some $u_v\leq u_v'$. Since $\mathfrak{F}$ is $\geq$-saturated, $u_v$ sees an $\geq$-endpoint, and thus we may w.l.o.g. assume that $u_v$ is itself an $\geq$-endpoint. 

Consider the set 
\(E=\{u_v\mid v\in E_w\}\). For each $u_v\in E$, by the construction we have $v\in E_w$ and $vRu_v'\geq u_v$, which by (F2) implies that  there exists $v'\in W$ such that $v\geq v'Ru_v$. But as $v$ is an $\geq$-endpoint, we must have  $v=v'$ and $vRu_v$. Thus, we have proved that the set $E$ satisfies the condition $E\subseteq R(E_w)$ in (G2). 
On the other hand, for every $v\in E_w$, we have $u_v\in E$ by definition, and the same argument as above shows that $u_v\in R(v)$. Hence, $u_v\in E\cap R(v)\neq \emptyset$, namely $E$ also satisfies the other condition  in (G2). Then, (G2) applies to the set $E$ and the point $w$, and therefore there exists a point $t\in W$ such that
\(wRt\text{ and }E_t\subseteq E.\)

Now, since $E_t\subseteq E$, every $\geq$-endpoint that $t$ can see is a $u_v\in E$ with $(\mathfrak{F},V),u_v\Vdash p$ for some $v\in E_w$. This means $(\mathfrak{F},V),t\Vdash \neg\neg p$, which gives $(\mathfrak{F},V),w\Vdash \Diamond\neg\neg p$ as $wRt$.

% Since $v$ is an $\geq$-endpoint and $\mathfrak{F}$ is $\geq$-saturated, there exists an $\geq$-endpoint $u_v$ such that $vR\circ \geq u_v$ and $(\mathfrak{F},V),u_v\Vdash p$.
%By (F2), there exists $v'\in W$ such that $v\geq v'$ and $v'Ru_v$. But as $v$ is a $\geq$-endpoint, we must have that $v=v'$ and $vRu_v$. Thus, the set
%\(E=\{u_v\mid v\in E_w\}\)
%satisfies the condition in (G2).
%By (G2), there exists a point $t\in W$ such that
%\(wRt\text{ and }E_t\subseteq E.\)
%Hence, $(\mathfrak{F},V),w\Vdash \Diamond\neg\neg p$.

Conversely, suppose $\mathfrak{F}$ does not satisfy (G2). Then there exists $w\in W$ and a set $E$ of $\geq$-endpoints satisfying $E\subseteq R(E_w)$ and $E\cap R(v)\neq \emptyset$ for every $v\in E_w$ such that  
%for each $v\in E_w$, there exists $u_v\in E$ with $v R u_v$, 
 for all $t\in W$,
$wRt$ implies $E_t\nsubseteq E$. Since $E$ is a set of $\geq$-endpoints, one can find a $\geq$-monotone valuation $V$ such that $V(p)=E$. We will show that $(\mathfrak{F},V),w\nVdash \neg\Box\neg p\to \Diamond\neg\neg p$.

For every $v\in E_w$, there exists $u\in E\cap R(v)\neq \emptyset$ with $(\mathfrak{F},V),u\Vdash p$. Since $(\mathfrak{F},V),u\nVdash \neg p$ and $vRu$, we obtain $(\mathfrak{F},V),v\nVdash\Box\neg p$ for every $v\in E_w$. Hence, $(\mathfrak{F},V),w\Vdash\neg\Box\neg p$. 

On the other hand, for every $t\in R(w)$, by the assumption there exists $s\in E_t$ such that $s\notin E$ meaning  $(\mathfrak{F},V),s\nVdash p$. It follows that $(\mathfrak{F},V),t\nVdash \neg\neg p$ for every $t\in R(w)$. Hence $(\mathfrak{F},V),w\nVdash \Diamond\neg\neg p$.
\end{proof}

%Now, we are ready to define the class $\mathsf{K}$ of generalizations of all powerset Kripke models.
%\begin{definition}

%\end{definition}

Let $\mathsf{M}$ be the class of all finite negative $\geq$-saturated bi-relation intuitionistic Kripke models  satisfying (G1) and (G2).
In the remainder of this section, we show that the system of \MIDz is complete with respect to $\mathsf{M}$, that is, we will prove the following theorem. The idea of the proof is inspired by that of Theorem 3.2.18 in \cite{ivano_msc}. Note that since \MIDz is not closed under uniform substitution, one can only obtain the completeness theorem in the sense of the theorem below for a class  $\mathsf{M}$ of \emph{models} (with restricted valuations) instead of a class of \emph{frames} (with arbitrary valuations).

%\todo{The logic is not closed under uniform substitution, it is not complete with respect to frames, but models with negative valuations}

\begin{theorem}\label{completeness_mid_K}
For any \MIDz-formula $\phi$, \(\vdash_{\MIDz}\phi\iff \mathsf{M}\Vdash\phi.\)
%$\models \phi$ iff for all finite bi-relation intuitionistic Kripke model $\MM\in \mathsf{M}$, $\MM\Vdash\phi$.
\end{theorem}
\begin{proof}[Proof of ``$\Longleftarrow$'']\let\qed\relax
We have checked that each (finite) powerset  model is  in $\mathsf{M}$. Then,
\begin{align*}
\mathsf{M}\Vdash\phi\Longrightarrow~~ &\MM^\circ\Vdash\phi\text{ for all finite powerset models }\MM^\circ\\
\Longrightarrow~~ &\MM\models\phi\text{ for all finite classical modal Kripke models }\MM~~\text{(by Lemma \ref{comodel_model})}\\
\Longrightarrow~~ &\vdash_{\MIDz}\phi~~\text{(by the finite model  property (\Cref{fmp}) and}\\
&\quad\quad\quad\quad\quad\quad\quad\quad\quad\quad\quad\text{the Completeness Theorem  of \MIDz)}.
\end{align*}
\hfill$\varclubsuit$
\end{proof}

To prove the other direction ``$\Longrightarrow$'' of the above theorem, we first show that every model in $\mathsf{M}$ can be mapped via a p-morphism into a finite powerset Kripke model. As p-morphisms are truth-preserving, the required result will then follow. Now, we recall the definition of p-morphisms of bi-relation intuitionistic Kripke models given by Wolter and Zakharyaschev in  \cite{WolterZakhIML99}. %As usual, such defined p-morphisms are truth preserving.

\begin{definition}\index{p-morphism}
Let $\MM_1=(W_1,\geq_1,R_1,V_1)$ and $\MM_2=(W_2,\geq_2,R_2,V_2)$ be bi-relation intuitionistic Kripke models. A function $f:W_1\to W_2$ is called a \emph{p-morphism} iff
\begin{description}
\item[P1] $w\in V_1(p)\Longleftrightarrow f(w)\in V_2(p)$ for all propositional variables $p$
\item[P2] $w\geq_1v$ $\Longrightarrow$ $f(w)\geq_2f(v)$
\item[P3] $wR_1v$ $\Longrightarrow$ $f(w)R_2f(v)$
\item[P4] $f(w)\geq_2v'$ $\Longrightarrow$ $\exists v\in W_1$ s.t.  $f(v)=v'$ and $w\geq_1v$
\item[P5] $f(w)R_2v'$ $\Longrightarrow$ $\exists v\in W_1$ s.t.  $v'\geq_2 f(v)$ and $wR_1v$
\item[P6] $f(w)(\geq_2 \circ R_2)v'$ $\Longrightarrow$ $\exists v\in W_1$ s.t. $w\geq_1\circ R_1 v$ and $f(v) \geq_2 v'$
\end{description}
\end{definition}
 
%The above defined p-morphisms is truth preserving, that is, 
\begin{theorem}[see  \cite{WolterZakhIML99}]\label{p_map_truth_bimod}
If $f:\MM_1\to\MM_2$ is a p-morphism between two bi-relation intuitionistic Kripke models $\MM_1$ and $\MM_2$, then $\MM_1,w\Vdash\phi\Longleftrightarrow \MM_2,f(w)\Vdash\phi$.
%\[\MM_1,w\Vdash\phi\Longleftrightarrow \MM_2,f(w)\Vdash\phi.\]
\end{theorem}
%\begin{proof}
%By induction on $\phi$. 
%\end{proof}

%Next, we prove the crucial lemma for the proof of the direction ``$\Longrightarrow$'' of Theorem \ref{completeness_mid_K}.

\begin{lemma}\label{bi_mod_p_map_comod}
For every finite bi-relation intuitionistic Kripke model $\MM=\mathop{(W,\geq,R,V)}$ in $\mathsf{M}$, there exists a  finite classical modal Kripke model $\mathfrak{N}$ such that there exists a  p-morphism $f$ of \MM into the powerset model $\mathfrak{N}^\circ$ induced by $\mathfrak{N}$.
\end{lemma}
\begin{proof}
Define a modal Kripke model $\mathfrak{N}=(W_0,R_0,V_0)$ as follows:
\begin{itemize}
\item $W_0$ is the set of all $\geq$-endpoints of $W$,
\item $R_0=R\upharpoonright W_0$ and $V_0=V\upharpoonright W_0$.
\end{itemize}
Now, consider the powerset Kripke model $\mathfrak{N}^\circ=(W_0^\circ,\supseteq,R_0^\circ,V_0^\circ)$ associated with $\mathfrak{N}$.  Define a function $f:W\to W^\circ_0$ by taking
\[f(w)=E_w\text{ for all }w\in W.\]
Since \MM is saturated, $E_w\neq\emptyset$ for all $w\in W$. Thus $E_w\in W_0^\circ$ and $f$ is well-defined.

%Before we continue the proof, let us ponder over the above construction. %Singletons of $\mathfrak{N}^\circ$ are singletons of $\geq$-endpoints of $W$. 
%As defined, 
Note that an $\geq$-endpoint $e$ of \MM is mapped through $f$ to the singleton $\{e\}=E_e$. Intuitively,  $\geq$-endpoints of \MM are simulated in our argument by singletons of $\mathfrak{N}^\circ$, and  it may be helpful for the reader to think of a node $w$ of \MM as the team formed by all $\geq$-endpoints seen from $w$, namely the set $E_w$.

% As it is defined, elements of powerset Kripke models resemble teams of the model 

%Intuitively, $\geq$-endpoints of $W$ are singleton elements in $W_0^\circ$.

 Now, we proceed to show that $f$ is a p-morphism, i.e., $f$ satisfies (P1)-(P6).

(P1). It suffices to show that $\MM,w\Vdash p\iff \mathfrak{N}^\circ,E_w\Vdash p$. The direction ``$\Longrightarrow$'' follows from the $\geq$-monotonicity of $V$. For the direction ``$\Longleftarrow$'', if $\MM,w\nVdash p$, then since $V$ is negative, $\MM,w\nVdash \neg \neg p$. Thus, there exists $v\in E_w$ such that $\MM,v\nVdash  p$, which implies that  $\mathfrak{N},\{v\}\not\models  p$, thereby $\mathfrak{N}^\circ,E_w\nVdash p$.

(P2). Clearly, if $w\geq v$, then $E_w\supseteq E_v$, i.e. $f(w)\supseteq f(v)$.

(P3). Assume $wRv$, we show that $E_wR_0^\circ E_v$, namely $E_wR_0E_v$. For any $s\in E_w$, by (F1) of \MM, there exists $t\in W$ such that
\(v\geq t\text{ and }sRt.\)
For each $t'\in E_v$ such that $t\geq t'$, by (F2), there exists $s'\in W$ such that
\(s\geq s'\text{ and }s'Rt'.\)
As $s$ is an $\geq$-endpoint, we must have $s=s'$ and $sRt'$. 

On the other hand, for any $t\in E_v$, consecutively applying (F2) and (F1) of \MM, by a similar argument to the above, we can find an $s'\in E_w$ such that $s'Rt$. Hence, we conclude that $E_wR_0E_v$.

%by (F2) of \MM, there exists $s\in W$ such that 
%\[w\geq s\text{ and }sRt.\]
%For each $s'\in E_w$ such that $s\geq s'$, by (F1), there exists $t'\in W$ such that
%\[t\geq t'\text{ and }s'Rt'.\]
%As $t$ is a $\geq$-endpoint, it must be that $t=t'$, so $sRt'$.

(P4). If $E_w\supseteq v'$, then as \MM is $\geq$-saturated, there exists $v\in W$ such that $w\geq v$ and $E_v=v'$, as required.

(P5). If $E_wR_0^\circ v'$, then $E_wR_0v'$. Clearly, $v'$ is a set of $\geq$-endpoints such that $v'\subseteq R(E_w)$ and $v'\cap R(s)\neq\emptyset$ for every $s\in E_w$.
Thus, by (G2) of \MM, there exists $v\in W$ such that
\(wRv\text{ and }v'\supseteq E_v,\)
as required.

(P6) Suppose $E_w(\supseteq \circ R_0^\circ) u'$ and $u'\neq \emptyset$. Then
\(u'\subseteq\, (\geq\circ R) (w).\) Since \MM is finite, the set $u'$ must be finite and (G1) applies.
Thus, there exists $u\in W$ such that
\(w(\geq\circ R)u\text{ and }u\geq s\text{ for all }s\in u'.\)
Since $u'$ is a set of $\geq$-endpoints, the latter of the above implies that $f(u)=E_u\supseteq u'$.
\end{proof}

Finally, we complete the proof of Theorem \ref{completeness_mid_K} as follows.
\begin{proof}[Proof of Theorem \ref{completeness_mid_K}, the direction ``$\Longrightarrow$''] Suppose $\vdash_\MIDz\phi$. For each $\mathfrak{M}\in \mathsf{M}$, by  \Cref{bi_mod_p_map_comod}, there is a classical modal Kripke model $\mathfrak{N}$ and a p-morphism $f:\mathfrak{M}\to\mathfrak{N}^\circ$. By the assumption, we have $\mathfrak{N}\models\phi$, which implies $\mathfrak{N}^\circ\Vdash\phi$ by \Cref{comodel_model}. Finally, by \Cref{p_map_truth_bimod}, we conclude that $\MM\Vdash\phi$, as required.
\end{proof}

Note that the finiteness of the models in the class $\mathsf{M}$ is used in the proof of \Cref{completeness_mid_K} (only) for establishing condition (P6) in \Cref{bi_mod_p_map_comod} when quoting condition (G1), the frame condition that the axiom $\Box(p\vee q)\to(\Box p\vee\Box q)$ characterizes (see \Cref{sat_G1}). Consider a  stronger version of (G1) :
\begin{description}
\item[G1${}^+$]For any $w\in W$ and any nonempty  set $X\subseteq (\geq\circ R)(w)$, there exists a node $u\in (\geq\circ R)(w)$ such that $u\geq v$ for all $v\in X$.
\end{description}
 If one, instead, defines $\mathsf{M}$ as the class of all (possibly infinite) negative $\geq$-saturated bi-relation intuitionistic Kripke models  satisfying  (G1${}^+$) and (G2), \Cref{completeness_mid_K} will still hold. But we choose to adopt the current setting in this section, as it exhibits more interaction between the properties of the models and the axioms.

%\todo{Classical dependence logic has Kripke models with no accessibility relation. This construction for classical dependence logic generates the variant of \KP. }

\subsection{A single-world semantics for \MTz}\label{sec:tensor_dia}
In this section, we define a single-world semantics for the system of the dependence atom-free fragment of \MT (\MTz), and prove that \MTz is complete (in the single-world semantics sense) with respect to a class of tri-relation intuitionistic Kripke models.
  
The language of \MTz has one connective more than that of \MID. The team semantics of the one additional connective of \MTz, the tensor $\sor$,  is generalized naturally  from the usual single-world semantics of  the disjunction of classical logic. Yet, the tensor, being understood as a disjunction, has a few odd behaviors. For instance, it does not admit the usual elimination and distributive rules for disjunction, in particular, none of $\phi\sor\phi\models\phi$, $(\phi\sor\psi)\wedge (\phi\sor\chi)\models\phi\sor(\psi\wedge\chi)$ and $(\phi\wedge\psi)\sor(\phi\wedge\chi)\models\phi\wedge(\psi\sor\chi)$ is in general true.  Indeed, although the tensor behaves truly as a disjunction over classical formulas (Cf. \Cref{class_K_syntax_equiv}), Abramsky and V\"{a}\"{a}n\"{a}nen \cite{AbVan09} observed that the tensor should rather be understood as a \emph{multiplicative conjunction} (hence the notation $\sor$) as in linear logic (or, in fact, in bunched implication logic \cite{OHearnPym1999}). Multiplicative conjunction can often be read as a \emph{binary diamond} modality, and this is the interpretation that we will adopt  for tensor  in this section.

Following the approach of the previous section, we first show that the team semantics of \MTz over the usual classical modal Kripke models coincides with the single-world semantics of \MTz over the associated \emph{full  powerset models}, which are powerset models equipped also  with  a ternary relation $R_\otimes$ for the interpretation of the binary diamond $\otimes$, and have also the empty team in their domains in order to characterize the property that the constant $\bot$ is the neutral element of the tensor, i.e., $\bot\sor \phi\equiv\phi$.  Let us now define formally this stronger notion of powerset model.%, called \emph{full powerset model}.

%To interpret the binary diamond $\otimes$, the powerset models will also be equipped with  a ternary relation $R_\otimes$. In view of the fact that the constant $\bot$ is the neutral element of the tensor, i.e., $\bot\sor \phi\equiv\phi\equiv\phi\sor\bot$, the domain of a powerset model will also contain the empty team.

% the tensor of team semantics as the usual binary diamond of the single-world semantics over a modified version of powerset models.

%In the sense that shall become clear in the sequel, tensor is a binary diamond with 
%
%
%One solution is to treat $\sor$ as a binary diamond-like modality. One difficulty is to characterize the fact that the constant $\bot$ is a neutral element of the disjunction $\sor$, i.e., $\bot\sor \phi\equiv\phi$ for all formulas $\phi$. For this purpose, one may refine the powerset model notion to include the empty team in the model.

%\todo{Tensor is actually a binary diamond that happens to behave like a disjunction restricted to classical formulas. It preserves $\ior$: $\otimes(A,B\vee C)\equiv \otimes(A,B)\vee\otimes(A,C)$}

%To define the semantics of tensor as a binary modality, we will add a ternary relation into the powerset models. Moreover, to  characterize the peculiar behavior of tensor as also a type of disjunction, especially to characterize the fact that the constant $\bot$ is the neutral element of the tensor disjunction, i.e., $\bot\sor \phi\equiv\phi$ for all formulas $\phi$, we will need to include also the empty team  in the powerset models. Let us now define this stronger notion of powerset model.

\begin{definition}
Let $\MM=(W,R,V)$ be a classical modal Kripke model. The \emph{full powerset  model} $\MM^\bullet$ induced by $\MM$ is a quintuple $\MM^\bullet=(W^\bullet,\supseteq,R^\circ,R_\sor,V^\circ)$, where  
\begin{itemize}
\item $W^\bullet=\wp(W)$ i.e., $W$ consists of all teams $X\subseteq W$ including the empty team $\emptyset$
\item $R_\sor$ is a ternary relation defined as $R_\sor(X,Y,Z)$ iff $X=Y\cup Z$
\item and the other components are defined as in \Cref{pwset_ml_df}.
\end{itemize}
\end{definition}

Full powerset models are special cases of tri-relation intuitionistic Kripke models defined as follows.

\begin{definition}\label{trirel_df}
A \emph{tri-relation intuitionistic Kripke frame} is a quardruple \allowbreak$\FF=(W,\geq,R,S)$, where $(W,\geq,R)$ is a bi-relation intuitionistic Kripke frame and $S$ is a binary relation on $W$ satisfying condition (H1) defined below:
\begin{description}
%\item[(H1)] If $S(w,u,v)$, then $w\geq u$ and $w\geq v$.
\item[H1] If $S(w,u,v)$ and $w\geq w'$, then there exist $u',v'\in W$ such that $S(w',u',v')$, $u\geq u'$ and $v\geq v'$.
%\item[(H3)]  If $E_s\subseteq E_u$, $E_t\subseteq E_v$ and $E_s\cup E_t=E_u\cup E_v$, then $S(w,u,v)\Rightarrow S(w,s,t)$, where we redefine the set $E_x$ as the set of all $\geq$ second least points seen from $x$.  \todo{Is this necessary?}
\end{description}
A tri-relation intuitionistic Kripke model is a tuple $(\FF,V)$ such that $\FF$ is a tri-relation intuitionistic Kripke frame and $V$ satisfies $\geq$-monotonicity.
\end{definition}

%Full powerset models are clearly tri-relation intuitionistic Kripke models. %Next, we define the satisfaction relation for tri-relation intuitionistic Kripke models.

\begin{definition}
Let $\MM=\mathop{(W,\geq,R,S,V)}$ be a tri-relation intuitionistic Kripke model. We define the  \emph{satisfaction} relation $\MM,w\Vdash^\bullet\phi$ inductively as follows:
\begin{itemize}
%\addtolength{\itemsep}{-0.3\baselineskip}
\item $\MM,w\Vdash^\bullet\bot$ iff $w$ is an $\geq$-endpoint
\item $\MM,w\Vdash^\bullet\otimes(\phi,\psi)$ iff there exist $v,u\in W$ such that $S(w,v,u)$, $\MM,v\Vdash^\bullet\phi$ and $\MM,u\Vdash^\bullet\psi$
\item The other cases are defined the same way as the $\Vdash$ relation in \Cref{sat_df_bir}.
\end{itemize}
\end{definition}

Note that under the above definition the constant $\bot$ does not any more behave as the falsum of the logic (which is satisfied nowhere in any model) and the semantics of the negation $\neg\phi=\phi\to\bot$ is changed accordingly.  The real falsum, denoted by $\mathbf{0}$, will not be studied in the present paper.

%We will not consider the real falsum (denoted by $\mathbf{0}$) in this paper. 

%We now check that $\Vdash^\bullet$ is also $\geq$-monotone. 

Observe that condition (H1) in \Cref{trirel_df} characterizes a similar type of interaction between $S$ and $\geq$ to the interaction between $R$ and $\geq$ that is characterized by condition (F1) in \Cref{bi-relation_ml}. Condition (F1) ensures the monotonicity of the (unary) diamond $\Diamond$. We now prove that the binary diamond $\otimes$ as well as $\bot$ preserve monotonicity by applying (H1).

\begin{lemma}[Monontonicity]
If $\MM,w\Vdash^\bullet\phi$ and $w\geq u$, then $\MM,u\Vdash^\bullet\phi$.
\end{lemma}
\begin{proof}
We prove the lemma  by induction on $\phi$. We only check the interesting cases. 

If $\phi=\bot$, then $\MM,w\Vdash^\bullet\bot$ implies that $w$ is an $\geq$-endpoint. If $w\geq u$, then $w=u$, and so $\MM,u\Vdash^\bullet\bot$.

If $\phi=\psi\sor\chi$, then $\MM,w\Vdash^\bullet\psi\sor\chi$ implies that there exist $s,t\in W$ such that $S(w,s,t)$, $\MM,s\Vdash^\bullet\psi$ and $\MM,t\Vdash^\bullet\chi$. Since $w\geq u$, by (H1), there exist $s',t'\in W$ such that $S(u,s',t')$, $s\geq s'$ and $t\geq t'$. By the induction hypothesis, we have $\MM,s'\Vdash^\bullet\psi$ and $\MM,t'\Vdash^\bullet\chi$. Thus, $\MM,u\Vdash^\bullet\psi\sor\chi$.
\end{proof}

%We now show that  the team-based satisfaction relation with respect to classical modal Kripke models is equivalent to the single-world-based satisfaction relation with respect to full powerset  models. The reader may compare the lemma below with \Cref{comodel_model}.

\begin{lemma}\label{fulpwmodel_model}
Let $\MM=(W,R,V)$ be a classical modal Kripke model and $X\subseteq W$ a team. For any \MTz-formula $\phi$, $\MM,X\models \phi\iff \MM^\bullet,X\Vdash^\bullet\phi$.
\end{lemma}
\begin{proof}
We prove the lemma by induction on $\phi$.  If $\phi=\bot$, then $\MM,X\models\bot$ iff $X=\emptyset$ iff $\MM^\bullet,X\Vdash^\bullet\bot$, since $\emptyset$ is an $\supseteq$-endpoint in $\MM^\bullet$.

%If $\phi=p$, then $\MM,X\models p\iff \MM^\bullet,X\Vdash p$ by the definition of $V^\circ$. In particular, we have $\MM,\emptyset\models p$ and $\MM^\bullet,\emptyset\Vdash p$.

If $\phi=\psi\sor\chi$, then 
\begin{align*}
\MM,X\models\psi\sor\chi&\iff\exists Y,Z\text{ s.t. }X=Y\cup Z,~\MM,Y\models\psi\text{ and }\MM,Z\models\chi\\
&\iff\exists Y,Z\text{ s.t. }R_\sor(X,Y,Z),~\MM^\bullet,Y\Vdash^\bullet\psi\text{ and }\MM^\bullet,Z\Vdash^\bullet\chi\\
&\quad\quad\quad\text{ (by the induction hypothesis)}\\
&\iff \MM^\bullet,X\Vdash^\bullet\psi\sor\chi.
\end{align*}

The other cases follow from the same argument as in \Cref{comodel_model}.
\end{proof}

We remarked in the previous section that singletons in a powerset model correspond to $\geq$-endpoints in its associated bi-relation intuitionistic Kripke models. In a tri-relation intuitionistic Kripke model $\MM=\mathop{(W,\geq,R,S,V)}$, singletons are simulated by the \emph{$\geq$-second least points} instead, i.e., points $w$ in $W$ such that there is an $\geq$-endpoint $e$ such that $w>e$, and  for all $v\in W$, $w>v$ implies that  $v$ is an $\geq$-endpoint.
%$w\geq v$ implies that $v=w$ or $v$ is an $\geq$-endpoint. 
We denote by $E_w^\bullet$ the set of all $\geq$-second least points seen from $w$, i.e.,
\[E_w^\bullet=\{v\in W\mid w\geq v\text{ and $v$ is a $\geq$-second least point}\}.\] 
%We now  re-interpret the set $E_w$ defined in the previous section as the set of all $\geq$-second least points seen from $w$, i.e.,
%\[E_w=\{v\in W\mid w\geq v\text{ and $v$ is a $\geq$-second least point}\}.\] 
In particular, if $w$ is an $\geq$-endpoint, then $E_w^\bullet=\emptyset$; and if $w$ is itself a $\geq$-second least point, then $E_w^\bullet=\{w\}$.

We say that a tri-relation intuitionistic Kripke frame $\FF=(W,\geq,R,S)$ is \emph{weakly $\geq$-saturated} if $E_w^\bullet\neq\emptyset$ for every non-$\geq$-endpoint $w\in W$, and for every  subset $E\subseteq E_w^\bullet$, there exists a point $v\in W$ such that $w\geq v$ and $E_v^\bullet=E$. 
%For example, in \Cref{fig:tri-frame-exam} the frame $\FF_1$  is not weakly $\geq$-saturated, since $E_w\neq\emptyset$; while the frame $\FF_2$ is  weakly $\geq$-saturated where $w$ sees two $\geq$-endpoints. \todo{??} 
A model $\MM$ is called  \emph{weakly negative} if $\MM,w\Vdash^\bullet p\iff \MM,w\Vdash^\bullet \neg \neg p$, and for all $\geq$-endpoints $e$, $\MM,e\Vdash^\bullet p$. 

%\begin{figure}[t]
%\begin{center}
%\begin{tikzpicture}[scale=.8, transform shape] 
%\node  (1)  at (-7, 0) {\normalsize$w$};
%%\node  (2) at +(-8,1.2) {$u'$};
%\node (2) at +(-8,1) {$\bullet$};
%\node  (5) at +(-9,2) {$\bullet$};
%%\node  (6) at +(-6,1.2) {$v'$};
%\node  (6) at +(-6,1) {$\bullet$};
%\draw[line width=1pt] (1) -- (2);
%\draw[line width=1pt] (1) -- (6);
%\draw [line width=1pt] (2) -- (5);
%
%\node at (-7,-0.8) {\normalsize$\FF_1$};
%
%
%\node  (1)  at (0, 0.5) {\normalsize$w$};
%%\node  (2) at +(-8,1.2) {$u'$};
%\node (2) at +(-1,1.5) {$\bullet$};
%%\node  (5) at +(-2,2) {$\bullet$};
%%\node  (6) at +(-6,1.2) {$v'$};
%\node  (6) at +(1,1.5) {$\bullet$};
%\draw[line width=1pt] (1) -- (2);
%\draw[line width=1pt] (1) -- (6);
%%\draw [line width=1pt] (2) -- (5);
%
%\node at (0,-0.8) {\normalsize$\FF_2$};
%
%\end{tikzpicture}
%\caption{tri-relation intuitionistic Kripke frames}
%\label{fig:tri-frame-exam}
%\end{center}
%\end{figure}

Let $\mathsf{M}^\bullet$ be the class of all finite weakly negative and weakly $\geq$-saturated tri-relation intuitionistic Kripke models  satisfying (G1), (G2) with ``$\geq$-endpoints'' and ``$E_w$'' in the definition replaced by ``$\geq$-second least points'' and ``$E_w^\bullet$", respectively, and (H2), (H3) and (H4) defined as follows:
\begin{description}
\item[H2] $S(w,u,v)\iff E_w^\bullet= E_u^\bullet\cup E_v^\bullet$. 
\item[H3] For any $\geq$-endpoint $e$, $eRw$ or $wRe$ implies that $w$ is also an $\geq$-endpoint.
%If $wRv$, then $w$ is an $\geq$-endpoint if and only if $v$ is also an $\geq$-endpoint.%\todo{??}
%(b) $E_w=E_u\cup E_v$  $\Longrightarrow S(w,u,v)$ %\todo{Or just one inclusion?}
\item[H4] $eRe$ for all $\geq$-endpoints $e$
\end{description}
%It is easy to check that each finite full powerset  model is  in $\mathsf{M}^\bullet$. 

A full powerset model $\MM^\bullet$ has a unique $\supseteq$-endpoint, namely the empty set $\emptyset$. Since $\MM,\emptyset\models p$ for all $p$, by \Cref{comodel_model} we know that $\MM^\bullet,\emptyset\Vdash^\bullet p$ as well. %Together with what we obtained in the previous section, it is easy to see that $\MM^\bullet$ is a weakly negative model. 
We leave it for the reader to verify that all the other conditions of $\mathsf{M}^\bullet$ are satisfied by the full powerset model $\MM^\bullet$ of any finite modal Kripke model $\MM$, i.e., $\MM^\bullet\in\mathsf{M}^\bullet$.

The main result of this section is that the system of \MTz is complete %(in the sense of single-world semantics) 
with respect to the class $\mathsf{M}^\bullet$, namely, the following theorem holds.

\begin{theorem}\label{completeness_mt_M}
For any \MTz-formula $\phi$, \(\vdash_{\MTz}\phi\iff \mathsf{M}^\bullet\Vdash\phi.\)
%$\models \phi$ iff for all finite bi-relation intuitionistic Kripke model $\MM\in \mathsf{M}$, $\MM\Vdash\phi$.
\end{theorem}

The proof of the above theorem goes through a similar argument to that of \Cref{completeness_mid_K}. The direction ``$\Longleftarrow$" follows from \Cref{fulpwmodel_model} and the finite model property of \MTz  (\Cref{fmp}). The other direction ``$\Longrightarrow$" will follow from \Cref{p_map_truth_trimod} and \Cref{tri_p_map_bullet} to be stated and proved in the remainder of this section.

%Next, we prove the characterization property of double negation of propositional variables in models in $\mathsf{M}^\bullet$.

%we check that for models in $\mathsf{M}^\bullet$, double negation of propositional variables

We first prove a lemma concerning the behavior of the double negation under the satisfaction relation $\Vdash^\bullet$.

\begin{lemma}\label{dneg}
Let $\MM=(W,\geq,R,S,V)$ be a model in $\mathsf{M}^\bullet$ and $w$ is a non-$\geq$-endpoint in $\MM$. Then, $\MM,w\nVdash^\bullet\neg\neg p$ iff  $\MM,v\nVdash^\bullet p$ for some $v\in E_w^\bullet$.
\end{lemma}
\begin{proof}
``$\Longleftarrow$": Suppose $\MM,v\nVdash^\bullet p$ for some $v\in E_w^\bullet$. For any point $e$ such that $v\geq e$ and $v\neq e$, $e$ is an $\geq$-endpoint, implying $\MM,e\Vdash^\bullet  \bot$. It follows that $\MM,v\Vdash^\bullet  p\to \bot$. Since $\MM,v\nVdash^\bullet  \bot$, we conclude $\MM,w\nVdash^\bullet (p\to\bot)\to\bot$.

``$\Longrightarrow$'': Suppose $\MM,w\nVdash^\bullet \neg\neg p$. Then, there exists $u\leq w$ such that $\MM,u\Vdash^\bullet  p\to \bot$ and $\MM,u\nVdash^\bullet  \bot$. The latter implies that $u$ is not an $\geq$-endpoint. Since $\MM$ is weakly $\geq$-saturated, $E_u^\bullet\neq\emptyset$. Pick an element $v\in E_u^\bullet\subseteq E_w^\bullet$. By the monotonicity of $\geq$, we have $\MM,v\Vdash^\bullet  p\to \bot$, implying $\MM,v\nVdash^\bullet  p$, as required.
\end{proof} 

%\begin{theorem}\label{completeness_mt_M}
%\(\models_{\MT}\phi\iff \mathsf{M}^\bullet\models\phi.\)
%%$\models \phi$ iff for all finite bi-relation intuitionistic Kripke model $\MM\in \mathsf{M}$, $\MM\Vdash\phi$.
%\end{theorem}
%\begin{proof}[Proof of ``$\Longleftarrow$'']
%It is easy to check that each finite full powerset  model is  in $\mathsf{M}^\bullet$. Then,
%\begin{align*}
%\mathsf{M}^\bullet\models\phi\Longrightarrow~~ &\MM^\bullet\models\phi\text{ for all finite full powerset models }\MM^\bullet\\
%\Longrightarrow~~ &\MM\models\phi\text{ for all finite modal Kripke models }\MM\\
%\Longrightarrow~~ &\models_{\MT}\phi~~\text{(by the finite model property of \MT, \Cref{fmp})}.
%\end{align*}
%\end{proof}

The notion of p-morphism for tri-relation intuitionistic Kripke models is the p-morphism for bi-relation intuitionistic Kripke models parametrized by the standard clause for binary diamonds.

\begin{definition}\index{p-morphism_bullet}
Let $\MM_1=(W_1,\geq_1,R_1,S_1,V_1)$ and $\MM_2=(W_2,\geq_2,R_2,S_2,V_2)$ be tri-relation intuitionistic Kripke models. A function $f:W_1\to W_2$ is called a \emph{p-morphism} iff $f$ satisfies (P1)-(P6) and (Q1) and (Q2) defined below:
\begin{description}
\item[Q1] $S_1(w,u.v)$ $\Longrightarrow$ $S_2(f(w),f(u),f(v))$
\item[Q2] $S_2(f(w),u',v')$ $\Longrightarrow$ $\exists u,v\in W_1$ s.t.  $f(u)=u'$, $f(v)=v'$  and $S_1(w,u,v)$
\end{description}
\end{definition}

%It is easy to verify that p-morphisms $f$ are truth preserving, i.e., $\MM_1,w\Vdash\phi\Longleftrightarrow \MM_2,f(w)\Vdash\phi$.

%Next, we check that p-morphisms are truth-preserving. Note that since the constant $\bot$ has now a different semantics than the usual one, the conditions in the truth theorem of p-morphisms need to be strengthened. 

We call a p-morphism $f$ between two models $\MM_1=(W_1,\geq_1,R_1,V_1)$ and $\MM_2=(W_2,\geq_2,R_2,V_2)$ \emph{endpoint-preserving} if 
\begin{description}
\item[Q3] $e$ is an $\geq_1$-endpoint $\iff$ $f(e)$ is an $\geq_2$-endpoint.
\end{description}

\begin{theorem}\label{p_map_truth_trimod}
If $f$ is an endpoint-preserving p-morphism between  tri-relation intuitionistic Kripke models $\MM_1$ and $\MM_2$,
%\begin{equation}\label{p_map_truth_trimod_eq1}
%e\text{ is an $\geq_1$-endpoint }\iff f(e)\text{ is an $\geq_2$-endpoint},
%\end{equation}
 then $\MM_1,w\Vdash^\bullet\phi\Longleftrightarrow \MM_2,f(w)\Vdash^\bullet\phi$.
%\[\MM_1,w\Vdash\phi\Longleftrightarrow \MM_2,f(w)\Vdash\phi.\]
\end{theorem}
\begin{proof}
The theorem is proved by a routine argument.
%We prove the theorem by induction on $\phi$. We only check the interesting cases;  the other cases are proved as in \Cref{p_map_truth_bimod} or see  \cite{WolterZakhIML99}. The case $\phi=\bot$ follows from (Q3). The case $\phi=\psi\otimes\chi$ is proved by applying (Q1) and (Q2). 
%We prove the theorem by induction on $\phi$. If $\phi=\bot$, then by the semantics of $\bot$, it suffices to verify that $w$ is an $\geq_1$-endpoint if and only if $f(w)$ is an $\geq_2$-endpoint. 
%
%``$\Longrightarrow$": Suppose $w$ is an $\geq_1$-endpoint and $f(w)\geq_2v'$. We have to show that $f(w)=v'$. Now, by (P4), there exists $v\in W_1$ such that $f(v)=v'$ and $w\geq_1 v$. But since $w$ is an $\geq_1$-endpoint, we must have that $v=w$, and thus $f(w)=f(v)=v'$.
%
%``$\Longleftarrow$": Suppose $f(w)$ is an $\geq_2$-endpoint and $w\geq_1 v$. We have to show that $w=v$. Now, by (P2), we have $f(w)\geq_2f(v)$. But since $f(w)$ is an $\geq_2$-endpoint, we must have that $f(w)=f(v)$
\end{proof}

Finally, we prove the crucial lemma of this section, from which the direction ``$\Longrightarrow$" of \Cref{completeness_mt_M} will follow. The reader may compare this lemma with \Cref{bi_mod_p_map_comod}. %in the previous section.

\begin{lemma}\label{tri_p_map_bullet}
For every finite tri-relation intuitionistic Kripke model $\MM$ in $\mathsf{M}^\bullet$, there exists a  finite classical modal Kripke model $\mathfrak{N}$ such that there exists an endpoint preserving  p-morphism $f$ of \MM into the full powerset model $\mathfrak{N}^\bullet$ induced by $\mathfrak{N}$.
%\[e\text{ is an $\geq_1$-endpoint }\iff f(e)\text{ is an $\geq_2$-endpoint}.\]
\end{lemma}
\begin{proof}
Let $\MM=\mathop{(W,\geq,R,S,V)}$. Define a modal Kripke model $\mathfrak{N}=(W_0,R_0,V_0)$ as:
\begin{itemize}
\item $W_0$ is the set of all $\geq$-second least points of $W$,
\item $R_0=R\upharpoonright W_0$ and $V_0=V\upharpoonright W_0$.
\end{itemize}
Now, consider the full powerset Kripke model $\mathfrak{N}^\bullet=(W_0^\bullet,\supseteq,R_0^\circ,R_\sor,V_0^\circ)$ associated with $\mathfrak{N}$.  Define a function $f:W\to W^\bullet_0$ by taking
\[f(w)=E_w^\bullet\text{ for all }w\in W.\]

We show that $f$ is an endpoint-preserving p-morphism. Conditions (P2) and (P4) are verified by a similar argument to that in the proof of \Cref{bi_mod_p_map_comod} taking into account the fact that $E_w$ is replaced by $E_w^\bullet$ in this proof. We now give the proof for the other conditions.

(P1). We show that $\MM,w\Vdash^\bullet p\iff \mathfrak{N}^\bullet,E_w^\bullet\Vdash^\bullet p$. If $e$ is an $\geq$-endpoint, since $\MM$ is weakly negative, $\MM,e\Vdash^\bullet p$. We also have $E_e^\bullet=\emptyset$, and $\mathfrak{N}^\bullet,\emptyset\Vdash^\bullet p$ by the empty team property and \Cref{fulpwmodel_model}. If $w$ is not an $\geq$-endpoint, applying \Cref{dneg}, the condition is proved by a similar argument to that in the proof of \Cref{bi_mod_p_map_comod}.

(P3). Assume $wRv$, we show that $E_w^\bullet R_0^\circ E_v^\bullet $, namely $E_w^\bullet R_0E_v^\bullet $. If one of $w$ and $v$ is an $\geq$-endpoint, then by (H3), both $w$ and $v$ are $\geq$-endpoints, which implies that $E_w^\bullet=\emptyset=E_v^\bullet $. By definition, $\emptyset R_0\emptyset$, i.e., $E_w^\bullet R_0E_v^\bullet $.

Now, assume that both $w$ and $v$ are not $\geq$-endpoints. By a similar argument to that in \Cref{bi_mod_p_map_comod}, for any $s\in E_w^\bullet$, by (F1) of \MM, there exists $t\in W$ such that
\(v\geq t\text{ and }sRt.\)
For each $t'\in E_v^\bullet $ such that $t\geq t'$, by (F2), there exists $s'\in W$ such that
\(s\geq s'\text{ and }s'Rt'.\)
But since $s$ is a $\geq$-second least point, either $s'=s$ or $s'$ is an $\geq$-endpoint. In the latter case, we conclude from (H3) that $t'$ is also an $\geq$-endpoint, which is a contradiction. Thus, the former is the case, and $sRt'$. On the other hand, for any $t\in E_v^\bullet $, by a similar argument, we can find an $s'\in E_w^\bullet$ such that $s'Rt$. Hence we conclude that $E_w^\bullet R_0E_v^\bullet $.

(P5). Suppose $E_w^\bullet R_0^\circ v'$. If one of $E_w^\bullet$ and $v'$ is the empty set, then $E_w^\bullet=\emptyset=v'$ by the definition of $R_0^\circ$. Since $\MM$ is weakly $\geq$-saturated, $w$ must be an $\geq$-endpoint, which by (H4) implies that $wRw$. Also, clearly, $v'\supseteq \emptyset =E_w^\bullet=f(w)$, as required. If $E_w^\bullet,v'\neq\emptyset$, then the condition follows from a similar argument to that in the proof of \Cref{bi_mod_p_map_comod}.

(P6). Suppose $E_w^\bullet(\supseteq \circ R_0^\circ) v'$. %If $E_w^\bullet=\emptyset$, then $v'=\emptyset$ and $w$ is an $\geq$-endpoint. Thus, for $v=w$, by (H4) we have $w\geq\circ R v$ and $f(v)=f(w)=E_w^\bullet=\emptyset\supseteq v'$. 
If $v'=\emptyset$, then let $v$ be any $\geq$-endpoint such that $w\geq v$. Clearly, $f(v)=E_v^\bullet =\emptyset\supseteq v'$. By (H4), $vRv$, and thus $w(\geq\circ R)v$. If $v'\neq\emptyset$, then the condition follows from a similar argument to that in the proof of \Cref{bi_mod_p_map_comod}.

(Q1). Suppose $S(w,u,v)$. By (H2), $E_w^\bullet=E_u^\bullet \cup E_v^\bullet $, thereby $R_\sor(f(w),f(u),f(v))$.

(Q2). Suppose $R_\sor(f(w),u',v')$. Then $E_w^\bullet=u'\cup v'$. Since $\MM$ is weakly $\geq$-saturated, there exist $u,v\in W$ such that $w\geq u,v$, $E_u^\bullet =u'$ and $E_v^\bullet =v'$. It follows that $f(u)=u'$, $f(v)=v'$. Since $E_w^\bullet=E_u^\bullet \cup E_v^\bullet $, by (H2), we obtain $S(w,u,v)$.

(Q3). If $e$ is an $\geq$-endpoint, then $f(e)=E^\bullet_e=\emptyset$, which is the unique $\supseteq$-endpoint of $\mathfrak{N}^\bullet$. Conversely, if $w$ is not an $\geq$-endpoint, then since $\MM$ is weakly $\geq$-saturated, we have $E_w^\bullet\neq \emptyset$, i.e., $f(w)$ is not an $\supseteq$-endpoint.
\end{proof}

%\begin{theorem}\label{completeness_mt_M}
%\(\vdash_{\MT}\phi\iff \mathsf{M}^\bullet\models\phi.\)
%%$\models \phi$ iff for all finite bi-relation intuitionistic Kripke model $\MM\in \mathsf{M}$, $\MM\Vdash\phi$.
%\end{theorem}
%\begin{proof}
%By a similar argument to that of the proof of \Cref{completeness_mid_K}.
%\end{proof}

\section{Concluding remarks}

In this paper, we have studied the axiomatization problem and some model-theoretic properties of the major modal dependence logics considered in the literature, namely \MT, \MID, \MDor, \MD and \MDe. In the first part of the paper, we introduced sound and complete Hilbert-style or natural deduction systems for all these logics, among which those logics with intuitionistic implication have not been axiomatized before.   We presented the system of \MT as an extension of Fischer Servi's intuitionistic modal logic \IK and the Kreisel-Putnam intermediate logic \KP, and the systems of all the other modal dependence logics are its fragments and variants. We showed that formulas of all these modal dependence logics (essentially) enjoy a same disjunctive normal form that is essentially already known in the literature. We also derived some metalogical properties of the logics, such as Craig's Interpolation Theorem and the Finite Model Property, as immediate corollaries of the normal form.

%Normal forms are in general useful for logics, and we showed also some examples of deriving important properties of the logics, such as Craig's Interpolation Theorem and the Finite Model Property, by simple applications of the normal form.

%In the second part of the paper, we developed single-world semantics for the team-based logics, and showed that \MIDz and \MTz and sound and complete (in the single-world semantics sense) with respect to certain classes of intuitionistic Kripke models. 
First-order teams are essentially relations in first-order models, and first-order dependence logic is expressively equivalent to existential second-order logic \cite{Van07dl,KontVan09}. 
In a similar fashion,  in the second part of the paper
we interpreted modal teams as possible worlds in powerset models in \Cref{comodel_model,fulpwmodel_model}, and on the basis of this we showed in \Cref{completeness_mid_K,completeness_mt_M} that \MIDz and \MTz can be understood as intermediate modal logics also from the model-theoretic perspective in the sense that they are  complete (in the single-world semantics sense) with respect to certain classes of intuitionistic Kripke models.
It is worth pointing out that although \Cref{comodel_model,fulpwmodel_model} are not explicitly found in the literature, their intuitive idea seems to be folklore in the field or have in some sense already been used as a guideline in some research. For instance, the perfect information semantic set game introduced in \cite{VaMDL08} for modal dependence logic played over modal Kripke models can actually be viewed as a standard perfect information semantic game played over the associated full powerset models, and the correctness of the set game is essentially justified by \Cref{fulpwmodel_model}. It is the author's hope that the connections established in the paper between team semantics and single-world semantics, and between modal dependence logics and intermediate modal logics can provide a pointer for a deeper understanding of team semantics and team-based logics. 

As acknowledged in the corresponding sections, many results of this paper are built on or inspired by the literature of inquisitive logic. Inquisitive modal  logic (see e.g. \cite{Ciardelli_PhD}) can be viewed as a variant of model dependence logic with different modalities. It is interesting to see whether inquisitive modal logic can also be given a single-world semantics in a similar manner, and to compare inquisitive modal logic with intermediate modal logics.

\begin{acks}
The author would like to thank Giuseppe Greco, Dick de Jongh, Tadeusz Litak, Alessandra Palmigiano and Katsuhiko Sano for useful discussions related to this paper. The author is also grateful to an anonymous referee for helpful and stimulating comments concerning both the presentation and technical details of the paper.
\end{acks}

\bibliographystyle{jflnat_igpl} 
%\bibliography{../fan_delft}

\end{document}